\documentclass{amsart}
\usepackage{amssymb}

\usepackage{a4}

\theoremstyle{plain}
\newtheorem{thm}{Theorem}[section]
\newtheorem{cor}[thm]{Corollary}
\newtheorem{lemma}[thm]{Lemma}

\newtheorem{prop}[thm]{Proposition}

\theoremstyle{definition}
\newtheorem{definition}[thm]{Definition}
\newtheorem{notation}[thm]{Notation}

\theoremstyle{remark}
\newtheorem{rem}[thm]{Remark}

\theoremstyle{plain}

\numberwithin{equation}{section}

\begin{document}
\title[Path integral formulas on manifolds]{Finite dimensional approximations to Wiener measure and path integral
formulas on manifolds}
\author[Andersson]{Lars Andersson\footnotemark {$^*$}}
\thanks{\footnotemark {$^*$}Supported in part by NFR, contract no. F-FU 4873-307.}
\address{Department of Mathematics \\
Royal Institute of Technology \\
S-100 44 Stockholm, Sweden}
\email{larsa{\char'100}math.kth.se}
\author[Driver]{Bruce K. Driver\footnotemark {$^\dagger$}}
\thanks{\footnotemark {$^\dagger$}This research was partially supported by NSF Grant
DMS 96-12651.}
\address{Department of Mathematics, 0112\\
University of California, San Diego \\
La Jolla, CA 92093-0112 }
\email{driver{@}euclid.ucsd.edu}
\date{\today \ \emph{File:\jobname{.tex}}}

\begin{abstract}
Certain natural geometric approximation schemes are developed for Wiener
measure on a compact Riemannian manifold. These approximations closely mimic
the informal path integral formulas used in the physics literature for
representing the heat semi-group on Riemannian manifolds. The path space is
approximated by finite dimensional manifolds $\mathrm{H}_{\mathcal{P}}(M)$
consisting of piecewise geodesic paths adapted to partitions $\mathcal{P}$
of $[0,1].$ The finite dimensional manifolds $\mathrm{H}_{\mathcal{P}}(M)$
carry both an $H^{1}$ and a $L^{2}$ type Riemannian structures, $G_{\mathcal{%
P}}^{1}$ and $G_{\mathcal{P}}^{0}$ respectively. It is proved that 
\begin{equation*}
\frac{1}{Z_{\mathcal{P}}^{i}}e^{-\frac{1}{2}E(\sigma )}d\mathrm{Vol}_{G_{%
\mathcal{P}}^{i}}(\sigma )\rightarrow \rho _{i}(\sigma )d\nu (\sigma )\text{
as mesh}(\mathcal{P)}\rightarrow 0,
\end{equation*}
where $E(\sigma )$ is the energy of the piecewise geodesic path $\sigma \in 
\mathrm{H}_{\mathcal{P}}(M),$ and for $i=0$ and $1,$ $Z_{\mathcal{P}}^{i}$
is a ``normalization'' constant, $\mathrm{Vol}_{G_{\mathcal{P}}^{i}} $ is
the Riemannian volume form relative $G_{\mathcal{P}}^{i},$ and $\nu $ is
Wiener measure on paths on $M.$ Here $\rho _{1}(\sigma )\equiv 1$ and $\rho
_{0}(\sigma )=\exp \left( -\frac{1}{6}\int_{0}^{1}\mathrm{Scal} (\sigma
(s))ds\right) $ where $\mathrm{Scal}$ is the scalar curvature of $M.$ These
results are also shown to imply the well know integration by parts formula
for the Wiener measure.
\end{abstract}

\maketitle
\setcounter{tocdepth}{1}
\tableofcontents

\section{Introduction}

\label{s.1}Let $(M,g,o)$ be a Riemannian manifold $M$ of dimension $d$, with
Riemannian metric $g$ (we will also use $\langle \cdot ,\cdot \rangle $ to
denote the metric) and a given base point $o\in M$. Let $\nabla $ be the
Levi-Civita covariant derivative, $\Delta ={\text{tr}}\nabla ^{2}$ denote
the Laplacian acting on $C^{\infty }(M)$ \ and $p_{s}(x,y)$ be the
fundamental solution to the heat equation, $\partial u/\partial s=\frac{1}{2}%
\Delta u.$ More explicitly, $p_{s}(x,y)$ \ is the integral kernel of the
operator $e^{\frac{s}{2}\Delta }$ acting on $L^{2}(M,dx),$ where $dx$
denotes the Riemannian volume measure.

For simplicity we will restrict our attention to the case where $M$ is
either compact or $M$ is $\mathbb{R}^{d}.$ If $M=\mathbb{R}^{d}$, we will
always take $o=0$ and $\langle \cdot ,\cdot \rangle $ to be the standard
inner product on $\mathbb{R}^{d}.$ \ In either of these cases $M$ is
stochastically complete, i.e. $\int_{M}p_{s}(x,y)dy=1$ for all $s>0$ and $%
x\in M.$ Recall, for $s$ small and $x$ and $y$ close in $M,$ that 
\begin{equation}
p_{s}(x,y)\thickapprox (\frac{1}{2\pi s})^{d/2}e^{-\frac{1}{2s}d(x,y)^{2}},
\label{e.1.1}
\end{equation}
where $d(x,y)$ is the Riemannian distance between $x$ and $y.$ Moreover if $%
M=\mathbb{R}^{d},$ then $\Delta =\sum_{i=1}^{d}\partial ^{2}/\partial
x_{i}^{2},$ $d(x,y)=|x-y|$ and equation (\ref{e.1.1}) is exact.

\begin{definition}
\label{d.1.1}The \textbf{Wiener space} $\mathrm{W}([0,T];M),$ $T>0$ is the
path space 
\begin{equation}
\mathrm{W}([0,T];M)=\{\sigma :[0,T]\rightarrow M:\sigma (0)=o\text{ and }%
\sigma \text{ is continuous}\}.  \label{e.1.2}
\end{equation}
The \textbf{Wiener measure} $\nu $ associated to $(M,\langle \cdot ,\cdot
\rangle ,o)$ is the unique probability measure on $\mathrm{W}([0,T];M)$ such
that 
\begin{equation}
\int_{\mathrm{W}([0,T];M)}f(\sigma )d\nu _{T}(\sigma
)=\int_{M^{n}}F(x_{1},\ldots ,x_{n})\prod_{i=1}^{n}p_{\Delta
_{i}s}(x_{i-1},x_{i})dx_{1}\cdots dx_{n},  \label{e.1.3}
\end{equation}
for all functions $f$ of the form $f(\sigma )=F(\sigma (s_{1}),...,\sigma
(s_{n})),$ where $\mathcal{P}:=\{0=s_{0}<s_{1}<s_{2}<...<s_{n}=T\}$ is a
partition of $I:=[0,T],$ $\Delta _{i}s:=s_{i}-s_{i-1},$ and $%
F:M^{n}\rightarrow \mathbb{R}$ is a bounded measurable function. In equation
(\ref{e.1.3}), $dx$ denotes the Riemann volume measure on $M$ and by
convention $x_{0}:=o.$ For convenience we will usually take $T=1$ and write $%
\mathrm{W}(M)$ for $\mathrm{W}([0,1];M)$ and $\nu $ for $\nu _{1}.$\qed
\end{definition}

As is well known, there exists a unique probability measure $\nu _{T}$ on $%
\mathrm{W}([0,T];M)$ satisfying (\ref{e.1.3}). The measure $\nu _{T}$ is
concentrated on continuous but \emph{nowhere} differentiable paths. In
particular we get the following path integral representation for the heat
semi-group in terms of the measure $\nu _{T}$, 
\begin{equation}
e^{\frac{s}{2}\Delta }f(o)=\int_{\mathrm{W}([0,T];M)}f(\sigma (s))d\nu
_{T}(\sigma ),  \label{e.1.4}
\end{equation}
where $f$ is a continuous function on $M$ and $0\leq s\leq T.$

\begin{notation}
\label{n.1.2}When $M=\mathbb{R}^{d}$, $\langle \cdot ,\cdot \rangle $ is the
usual dot product and $o=0,$ the measure $\nu $ defined in Definition \ref
{d.1.1} is standard Wiener measure on $\mathrm{W}(\mathbb{R}^{d}).$ We will
denote this standard Wiener measure by $\mu $ rather than $\nu .$ We will
also let $B(s):\mathrm{W}(\mathbb{R}^{d})\rightarrow \mathbb{R}^{d}$ be the
coordinate map $B(s)(\sigma ):=B(s,\sigma ):=\sigma (s)$ for all $\sigma \in 
\mathrm{W}(\mathbb{R}^{d}).$ \qed
\end{notation}

\begin{rem}[Brownian Motion]
\label{r.1.3}The process $\left\{ B(s)\right\} _{s\in \lbrack 0,1]}$ is a
standard $\mathbb{R}^{d}$--valued Brownian motion on the probability space $(%
\mathrm{W}(\mathbb{R}^{d}),\mu ).$
\end{rem}

\subsection{A Heuristic Expression for Wiener Measure\label{s.1.1}}

Given a partition $\mathcal{P}:=\{0<s_{1}<s_{2}<...<s_{n}=1\}$ of $[0,1]$
and $\mathbf{x}:=(x_{1},\ldots ,x_{n})\in M^{n},$ let $\sigma _{\mathbf{x}}$
denote a path in $\mathrm{W}(M)$ such that $\sigma _{\mathbf{x}%
}(s_{i})=x_{i} $ and such that $\sigma _{\mathbf{x}}|_{[s_{i-1},s_{i}]}$ is
a geodesic path of shortest length for $i=1,2,\ldots n.$ (As above, $%
x_{0}:=o\in M.$) \ With this notation and the asymptotics for $p_{s}(x,y)$
in equation (\ref{e.1.1}), \ we find 
\begin{eqnarray*}
\prod_{i=1}^{n}p_{\Delta _{i}s}(x_{i-1},x_{i}) &\thickapprox
&\prod_{i=1}^{n}(\frac{1}{2\pi \Delta _{i}s})^{d/2}\exp \{-\frac{1}{2\Delta
_{i}s}d(x_{i-1},x_{i})^{2}\} \\
&=&\frac{1}{Z_{\mathcal{P}}}\exp \{-\frac{1}{2}\int_{0}^{1}\mid \sigma _{%
\mathbf{x}}^{\prime }(s)\mid ^{2}ds\},
\end{eqnarray*}
where $\sigma _{\mathbf{x}}^{\prime }(s):=\frac{d}{ds}\sigma _{\mathbf{x}%
}(s) $ for $s\notin \mathcal{P}$ and $Z_{\mathcal{P}}:=\prod_{i=1}^{n}(2\pi
\Delta _{i}s)^{d/2}.$\ Using this last expression in equation (\ref{e.1.3})
and letting the mesh of the partition $\mathcal{P}$ tend to zero we are lead
to the following \emph{heuristic} expression: 
\begin{equation}
d\nu (\sigma )\text{``}=\text{''}\frac{1}{Z}e^{-\frac{1}{2}E(\sigma )}%
\mathcal{D}\sigma ,  \label{e.1.5}
\end{equation}
where 
\begin{equation}
E(\sigma ):=\int_{0}^{1}\langle \sigma ^{\prime }(s),\sigma ^{\prime
}(s)\rangle ds  \label{e.1.6}
\end{equation}
is the \textbf{energy} of $\sigma ,$ $\mathcal{D}\sigma $ denotes a
``Lebesgue'' like measure on $\mathrm{W}(M)\ $and $Z$ is a ``normalization
constant'' chosen so as to make $\nu $ a probability measure.

Let $V$ be a continuous function on $M$. Then equation (\ref{e.1.5}) and
Trotter's product formula leads to the following heuristic path integral
formula for the parabolic heat kernel of the Schr\"{o}dinger operator $\frac{%
1}{2}\Delta -V$, 
\begin{equation}
e^{s(\frac{1}{2}\Delta -V)}f(o)\text{``}=\text{''}\frac{1}{Z}\int_{\mathrm{W}%
(M)}f(\sigma (1))e^{-(\frac{1}{2s}E(\sigma )+s\int_{0}^{1}V(\sigma (r))dr)}%
\mathcal{D}(\sigma )  \label{e.1.7}
\end{equation}

Equation (\ref{e.1.7}) can be interpreted as a prescription for the path
integral quantization of the Hamiltonian $\frac{1}{2}g^{ij}p_{i}p_{j}+V$.
The use of ``path integrals'' in physics including heuristic expressions
like those in equations (\ref{e.1.5}) and (\ref{e.1.7}) started with Feynman
in \cite{Feynman48} with very early beginnings being traced back to Dirac 
\cite{Dirac33}. See Gross \cite{Gross84} for a brief survey of the role of
path integrals in constructive quantum field theory and Glimm and Jaffe \cite
{Glimm87} for a more detailed account.

The heuristic interpretation of the ``measure'' $\mathcal{D}\sigma $ is
somewhat ambiguous in the literature. Some authors, for example \cite
{Dewitt72,Dewitt75,Dewitt79b,Dewitt81} tend to view $\mathrm{W}(M)$ as the
infinite product space $M^{I}$ and $\mathcal{D}\sigma $ as an infinite
product of Riemann volume measures on this product space. This is the
interpretation which is suggested by the ``derivation'' of equation (\ref
{e.1.5}) which we have given above.

Other authors, \cite{Atiyah85,Bismut85} interpret $\mathcal{D}\sigma $ as a
Riemannian ``volume form'' on $\mathrm{W}(M).$ We prefer this second point
of view. One reason for our bias towards the volume measure interpretation
is the fact that the path space $\mathrm{W}(M)$ is topologically trivial
whereas the product space $M^{I}$ is not. \ This fact is reflected in the
ambiguity (which we have glossed over) in assigning a path $\sigma _{\mathbf{%
x}}$ to a point $\mathbf{x}=(x_{1},\ldots ,x_{n})\in M^{n}$ as above in the
case when there are multiple distinct shortest geodesics joining some pair $%
(x_{i-1},x_{i}).$ However, from the purely measure theoretic considerations
in this paper we shall see that the two interpretations of $\mathcal{D}%
\sigma $ are commensurate.

Of course equations (\ref{e.1.5}) and (\ref{e.1.7}) are meaningless as they
stand because: 1) infinite dimensional Lebesgue measures do not exist and 2)
Wiener measure $\nu $ concentrates on nowhere differentiable paths which
renders the exponent in (\ref{e.1.5}) meaningless. Nevertheless, in Theorem 
\ref{t.1.8} we will give two precise interpretations of equation (\ref{e.1.5}%
).

\subsection{Volume elements on path space}

To make the above discussion more precise, let $\mathrm{H}(M)\subset \mathrm{%
W}(M)$ be the Hilbert manifold modeled on the space $\mathrm{H}(\mathbb{R}%
^{d})$ of finite energy paths: 
\begin{equation}
\mathrm{H}(M)=\{\sigma \in \mathrm{W}(M):\sigma \text{ is absolutely
continuous and }E(\sigma )<\infty \}.  \label{e.1.8}
\end{equation}
Recall that $\sigma \in \mathrm{W}(M)$ is said to be absolutely continuous
if $f\circ \sigma $ is absolutely continuous for all $f\in C^{\infty }(M).$
(It is easily checked that the space $\mathrm{H}(M)$ is independent of the
choice of Riemannian metric on $M.$)\vspace{1pt}\ The tangent space $%
T_{\sigma }\mathrm{H}(M)$ to $\mathrm{H}(M)$ at $\sigma $ may be naturally
identified with the space of absolutely continuous vector fields $%
X:[0,1]\rightarrow TM$ along $\sigma $ (i.e. $X(s)\in T_{\sigma (s)}M$ for
all $s$) such that $X(0)=0$ and $G^{1}(X,X)<\infty ,$ where 
\begin{equation}
G^{1}(X,X):=\int_{0}^{1}\left\langle \frac{\nabla X(s)}{ds},\frac{\nabla X(s)%
}{ds}\right\rangle ds,  \label{e.1.9}
\end{equation}
\begin{equation}
\frac{\nabla X(s)}{ds}:=//_{s}(\sigma )\frac{d}{ds}\{//_{s}(\sigma
)^{-1}X(s)\},  \label{e.1.10}
\end{equation}
and $//_{s}(\sigma ):T_{o}M\rightarrow T_{\sigma (s)}M$ is parallel
translation along $\sigma $ relative to the Levi-Civita covariant derivative 
$\nabla .$ \ See \cite
{Eells58,Palais63,Eliasson67,Klingenberg:closed,Flaschel:Klingenberg} for
more details.

By polarization, equation \ref{e.1.9} defines a Riemannian metric on $%
\mathrm{H}(M).$ Similarly we may define a ``weak'' Riemannian metric $G^{0}$
on $\mathrm{H}(M)$ by setting 
\begin{equation}
G^{0}(X,X):=\int_{0}^{1}\left\langle X(s),X(s)\right\rangle ds
\label{e.1.11}
\end{equation}
for all $X\in T\mathrm{H}(M).$ Given these two metrics it is natural to
interpret $\mathcal{D}\sigma $ as either of the (non-existent) ``Riemannian
volume measures'' $\mathrm{Vol}_{G^{1}}$ or $\mathrm{Vol}_{G^{0}}$ with
respect to $G^{1}$ and $G^{0}$ respectively. Both interpretations of $%
\mathcal{D}\sigma $ are $\emph{formally}$ the same modulo an infinite
multiplicative constant, namely the ``determinant'' of $\frac{d}{ds}$ acting
on $\mathrm{H}(T_{o}M).$

As will be seen below in Theorem \ref{t.1.8}, the precise version of the
heuristic expressions (\ref{e.1.5}) and (\ref{e.1.7}) shows that \emph{%
depending on the choice of volume form on the path space}, we get a scalar
curvature correction term. %with (referring to the discussion above) 
%$\kappa =0$ in the case of $G^1$ and $\kappa = \frac{1}{6}$ 
%in the case of
%$G^0$. 

\subsection{Statement of the Main Results}

\label{s.1.3}

In order to state the main results, it is necessary to introduce finite
dimensional approximations to $\mathrm{H}(M),$ $G^{1},$ $G^{0},$ $\mathrm{Vol%
}_{G^{1}}$ and $\mathrm{Vol}_{G^{0}}.$

\begin{notation}
\label{n.1.4}$\mathrm{H}_{\mathcal{P}}(M)=\{\sigma \in \mathrm{H}(M)\cap
C^{2}(I\setminus \mathcal{P}):\nabla \sigma ^{\prime }(s)/ds=0\text{\textrm{%
\ for }}s\notin \mathcal{P}\}$ --- the piecewise geodesics paths in $\mathrm{%
H}(M)$ which change directions only at the partition points.\qed
\end{notation}

It is possible to check that $\mathrm{H}_{\mathcal{P}}(M)$ is a finite
dimensional submanifold of $\mathrm{H}(M).$ Moreover by Remark \ref{r.4.3}
below, $\mathrm{H}_{\mathcal{P}}(M)$ is diffeomorphic to $\left( \mathbb{R}%
^{d}\right) ^{n}.$ For $\sigma \in \mathrm{H}_{\mathcal{P}}(M),$ the tangent
space $T_{\sigma }\mathrm{H}_{\mathcal{P}}(M)$ can be identified with
elements $X\in T_{\sigma }\mathrm{H}_{\mathcal{P}}(M)$ satisfying the Jacobi
equations on $I\setminus \mathcal{P},$ see Proposition \ref{p.4.4} below for
more details. We will now introduce Riemann sum approximations to the
metrics $G^{1}$ and $G^{0}.$

\begin{definition}[The $\mathcal{P}$--Metrics]
\label{d.1.5} For each partition $\mathcal{P}=\{0=s_{0}<s_{1}<s_{2}<\cdots
<s_{n}=1\}$ of $[0,1],$ let $G_{\mathcal{P}}^{1}$ be the metric on $T\mathrm{%
H}_{\mathcal{P}}(M)$ given by 
\begin{equation}
G_{\mathcal{P}}^{1}(X,Y):=\sum_{i=1}^{n}\langle \frac{\nabla X(s_{i-1}+)}{ds}%
,\frac{\nabla Y(s_{i-1}+)}{ds}\rangle \Delta _{i}s  \label{e.1.12}
\end{equation}
for all $X,Y\in T_{\sigma }\mathrm{H}_{\mathcal{P}}(M)$ and $\sigma \in 
\mathrm{H}_{\mathcal{P}}(M).$ (We are writing $\frac{\nabla X(s_{i-1}+)}{ds}$
as a shorthand for $\lim_{s\searrow s_{i-1}}\frac{\nabla X(s)}{ds}.)$
Similarly, let $G_{\mathcal{P}}^{0}$ be the \emph{degenerate}\textbf{\ }%
metric on $\mathrm{H}_{\mathcal{P}}(M)$ given by 
\begin{equation}
G_{\mathcal{P}}^{0}(X,Y):=\sum_{i=1}^{n}\langle X(s_{i}),Y(s_{i})\rangle
\Delta _{i}s,  \label{e.1.13}
\end{equation}
for all $X,Y\in T_{\sigma }\mathrm{H}_{\mathcal{P}}(M)$ and $\sigma \in 
\mathrm{H}_{\mathcal{P}}(M).$\qed
\end{definition}

If $N^{p}$ is an oriented manifold equipped with a possibly degenerate
Riemannian metric $G,$ let $\mathrm{Vol}_{G}$ denote the $p$--form on $N$
determined by 
\begin{equation}
\mathrm{Vol}_{G}(v_{1},v_{2},\dots ,v_{p}):=\sqrt{\det \left( \left\{
G(v_{i},v_{j})\right\} _{i,j=1}^{p}\right) },  \label{e.1.14}
\end{equation}
where $\left\{ v_{1},v_{2},\dots ,v_{p}\right\} \subset T_{n}N$ is an
oriented basis and $n\in N.$ We will often identify a $p$--form on $N$ with
the Radon measure induced by the linear functional $f\in C_{c}(N)\rightarrow
\int_{N}f\mathrm{Vol}_{G}.$

\begin{definition}[$\mathcal{P}$--Volume Forms]
\label{d.1.6} Let $\mathrm{Vol}_{G_{\mathcal{P}}^{0}}$ and $\mathrm{Vol}_{G_{%
\mathcal{P}}^{1}}$ denote the volume forms on $\mathrm{H}_{\mathcal{P}}(M)$
determined by $G_{\mathcal{P}}^{0}$ and $G_{\mathcal{P}}^{1}$ in accordance
with equation (\ref{e.1.14}). \qed
\end{definition}

Given the above definitions, there are now two natural finite dimensional
``approximations'' to $\nu $ in equation (\ref{e.1.5}) given in the
following definition.

\begin{definition}[Approximates to Wiener Measure]
\label{d.1.7} For each partition $\mathcal{P}=\{0=s_{0}<s_{1}<s_{2}<\cdots
<s_{n}=1\}$ of $[0,1],$ let $\nu _{\mathcal{P}}^{0}$ and $\nu _{\mathcal{P}%
}^{1}$ denote measures on $\mathrm{H}_{\mathcal{P}}(M)$ defined by 
\begin{equation*}
\nu _{\mathcal{P}}^{0}:=\frac{1}{Z_{\mathcal{P}}^{0}}e^{-\frac{1}{2}E} 
\mathrm{Vol}_{G_{\mathcal{P}}^{0}}
\end{equation*}
and 
\begin{equation*}
\nu _{\mathcal{P}}^{1}=\frac{1}{Z_{\mathcal{P}}^{1}}e^{-\frac{1}{2}E} 
\mathrm{Vol}_{G_{\mathcal{P}}^{1}},
\end{equation*}
where $E:\mathrm{H}(M)\rightarrow \lbrack 0,\infty )$ is the energy
functional defined in equation (\ref{e.1.6}) and $Z_{\mathcal{P}}^{0}$ and $%
Z_{\mathcal{P}}^{1}$ are normalization constants given by 
\begin{equation}
Z_{\mathcal{P}}^{0}:=\prod_{i=1}^{n}(\sqrt{2\pi }\Delta _{i}s)^{d}\text{ and 
}Z_{\mathcal{P}}^{1}:=(2\pi )^{dn/2}.  \label{e.1.15}
\end{equation}
\qed
\end{definition}

We are now in a position to state the main results of this paper.

\begin{thm}
\label{t.1.8} Suppose that $f:W(M)\rightarrow \mathbb{R}$ is a bounded and
continuous, then 
\begin{equation}
\lim_{|\mathcal{P}|\rightarrow 0}\int_{\mathrm{H}_{\mathcal{P}}(M)}f(\sigma
)d\nu _{\mathcal{P}}^{1}(\sigma )=\int_{\mathrm{W}(M)}f(\sigma )d\nu (\sigma
)  \label{e.1.16}
\end{equation}
and 
\begin{equation}
\lim_{|\mathcal{P}|\rightarrow 0}\int_{\mathrm{H}_{\mathcal{P}}(M)}f(\sigma
)d\nu _{\mathcal{P}}^{0}(\sigma )=\int_{\mathrm{W}(M)}f(\sigma )e^{-\frac{1}{%
6}\int_{0}^{1}\mathrm{Scal}(\sigma (s))ds}d\nu (\sigma ),  \label{e.1.17}
\end{equation}
where $\mathrm{Scal}$ is the scalar curvature of $(M,g).$ \qed
\end{thm}

Equation (\ref{e.1.16}) is a special case of Theorem \ref{t.4.17} which is
proved in Section \ref{s.4.1} and equation (\ref{e.1.17}) is a special case
of Theorem \ref{t.6.1} which is proved in Section \ref{s.6} below. An easy
corollary of equation (\ref{e.1.17}) of this theorem is the following
``Euler approximation'' construction for the heat semi-group $e^{s\Delta /2}$
on $L^{2}(M,dx).$ The following corollary is a special case of Corollary \ref
{c.6.7}

\begin{cor}
\label{c.1.9}For $s>0$ let $Q_{s}$ be the symmetric integral operator on $%
L^{2}(M,dx)$ defined by the kernel 
\begin{equation*}
Q_{s}(x,y):=\left( 2\pi s\right) ^{-d/2}\exp \left( -\frac{1}{2s}d^{2}(x,y)+%
\frac{s}{12}\mathrm{Scal}(x)+\frac{s}{12}\mathrm{Scal}(y)\right) \text{ for
all }x,y\in M.
\end{equation*}
Then for all continuous functions $F:M\rightarrow \mathbb{R}$ and $x\in M,$ 
\begin{equation*}
(e^{\frac{s}{2}\Delta }F)(x)=\lim_{n\rightarrow \infty }(Q_{s/n}^{n}F)(x).
\end{equation*}
\qed
\end{cor}

In the literature one often finds ``verifications'' (or rather tests) of
path integral formulas like (\ref{e.1.7}) by studying the small $s$
asymptotics. This technique, known as ``loop expansion'' or ``WKB
approximation'', when applied in the manifold case leads to the insight that
the operator constructed from the Hamiltonian $\frac{1}{2}g^{ij}p_{i}p_{j}+V$
depends sensitively on choices made in the approximation scheme for the path
integral. Claims have been made that the correct form of the operator which
is the path integral quantization of the Hamiltonian ${\frac{1}{2}}%
g^{ij}p_{i}p_{j}+V$ is of the form $-\mathbf{\hbar }^{2}(\frac{1}{2}\Delta %
-\kappa \mathrm{Scal})+V$ where $\mathbf{\hbar }$ is Planck's constant, $%
\mathrm{Scal}$ is the scalar curvature of $(M,g)$ and $\kappa $ is a
constant whose value depends on the authors and their interpretation of the
path integrals. Values given in the literature include $\kappa =\frac{1}{12}%
,\kappa =\frac{1}{6}$ \cite{deWitt80}, $\kappa =\frac{1}{8}$, \cite[Eq.
(6.5.25)]{deWitt:super} all of which are computed by formal expansion
methods. The ambiguity in the path integral is analogous to the operator
ordering ambiguity appearing in pseudo-differential operator techniques for
quantization, see the paper by Fulling \cite{Fulling95} for a discussion of
this point. In \cite{Fulling95} it is claimed that depending on the choice
of covariant operator ordering, the correction term has $\kappa $ ranging
from $0$ (for Weyl quantization) to $\frac{1}{6}$. For a discussion in the
context of geometric quantization, see \cite[\S 9.7]{Woodhouse92}, where the
value $\kappa =\frac{1}{12}$ is given for the case of a real
polarization.  In addition to the above one also finds in the
literature claims, based on perturbation calculations, that noncovariant
correction terms are necessary in path integrals, see for example 
\cite{deboer95} and references therein. 

It should be stressed that in contrast to the informal calculations
mentioned above, the results presented in Theorem \ref{t.1.8} and Corollary 
\ref{c.1.9} involve only well defined quantities. Let us emphasize that the
scalar curvature term appearing in equation (\ref{e.1.17}) has the nature of
a Jacobian factor relating the two volume forms $\mathrm{Vol}_{G_{\mathcal{P}%
}^{0}}$ and $\mathrm{Vol}_{G_{\mathcal{P}}^{1}}$ on path space.

We conclude this discussion by mentioning the so called Onsager-Machlup
function of a diffusion process. The Onsager-Machlup function can be viewed
as an attempt to compute an ``ideal density'' for the probability measure on
path space induced by the diffusion process. In the paper \cite
{takahashi:watanabe}, the probability for a Brownian path to be found in a
small tubular $\epsilon $--neighborhood of a \emph{smooth} path $\sigma $
was computed to be asymptotic to 
\begin{equation*}
Ce^{-\lambda _{1}/\epsilon ^{2}}\cdot \exp \left( -{\frac{1}{2}}E(\sigma
(s))+\frac{1}{12}\int_{0}^{1}\mathrm{Scal}(\sigma (r))dr\right) ,
\end{equation*}
where $\lambda _{1}$ is the first eigenvalue for the Dirichlet problem on
the unit ball in $\mathbb{R}^{d}$ and $C$ is a constant. The expression ${%
\frac{1}{2}}E(\sigma )-\frac{1}{12}\int_{0}^{1}\mathrm{Scal}(\sigma (r))dr$
thus recovered from the Wiener measure on $\mathrm{W}(M)$ is in this context
viewed as the action corresponding to a Lagrangian for the Brownian motion.
It is intriguing to compare this formula with equations (\ref{e.1.16}) and (%
\ref{e.1.17}).

An important result in the analysis on path space, is the formula for
partial integration. Here we use the approximation result in Theorem \ref
{t.1.8} to give an alternative proof of this result.

\begin{thm}
\label{t.1.10}Let $k\in \mathrm{H}(\mathbb{R}^{d})\cap C^{1}([0,1];\mathbb{R}%
^{d})$, $\sigma \in W(M)$ and $X_{s}(\sigma )\in T_{\sigma (s)}M$ be the
solution to 
\begin{equation*}
\frac{\nabla }{ds}X_{s}(\sigma )+\frac{1}{2}\mathrm{Ric}X_{s}(\sigma )=%
\widetilde{//}_{s}(\sigma )k^{\prime }(s)\quad \text{with }X_{0}(\sigma )=0,
\end{equation*}
where $\widetilde{//}_{s}(\sigma )$ denotes stochastic parallel translation
along $\sigma ,$ see Definition \ref{d.4.15}. Then for all smooth cylinder
functions $f$ (see Definition \ref{d.7.15}) on $\mathrm{W}(M),$%
\begin{equation*}
\int_{\mathrm{W}(M)}Xf\,d\nu =\int_{\mathrm{W}(M)}f\left(
\int_{0}^{1}\langle k^{\prime },d\tilde{b}\rangle \right) d\nu .
\end{equation*}
Here. $\tilde{b}$ is the $\mathbb{R}^{d}$ -- valued Brownian motion which is
the anti-development of $\sigma $, see Definition \ref{d.4.15} and $Xf\,$is
the directional derivative of $f$ with respect to $X,$ see Definition \ref
{d.7.15}.\qed
\end{thm}

Section \ref{s.7} is devoted to the proof of this result whose precise
statement may be found in Theorem \ref{t.7.16}.

\begin{rem}
This theorem first appeared in Bismut \cite{Bismut84a} in the special case
where $f(\sigma )=F(\sigma (s))$ for some $F\in C^{\infty }(M)$ and $s\in
\lbrack 0,1]$ and then more generally in \cite{Driver92b}. Other proofs of
this theorem may be found in \cite
{Aida97,Airault96,Driver97b,Elworthy94a,Elworthy96a,Elworthy96c,Enchev:Stroock95,Fang93,Hsu95c,Hsu95d,Leandre93b,Leandre95e,Lyons96b,Norris95}%
.
\end{rem}

\textbf{Acknowledgments: }The authors thank David Elworthy, Aubrey Truman
and Daniel Stroock for helpful remarks. The second author is grateful to the
Mathematical Science Research Institute, l'Institut Henri Poincar\'{e} and
l'\'{E}cole Normal Sup\'{e}riure where some of this work was done.

\section{Basic Notations and Concepts\label{s.2}}

\subsection{Frame Bundle and Connections\label{s.2.1}}

Let $\pi :O(M)\rightarrow M$ denote the bundle of orthogonal frames on $M$.
An element $u\in O(M)$ is an isometry $u:\mathbb{R}^{d}\rightarrow T_{\pi
(u)}M$. We will make $O(M)$ into a pointed space by fixing $u_{0}\in \pi
^{-1}(o)$ once and for all. We will often use $u_{0}$ to identify the
tangent space $T_{o}M$ of $M$ at $o$ with $\mathbb{R}^{d}$.

Let $\theta $ denote the $\mathbb{R}^{d}$--valued form on $O(M)$ given by $%
\theta _{u}(\xi )=u^{-1}\pi _{\ast }\xi $ for all $u\in O(M)$, $\xi \in
T_{u}O(M)$ and let $\omega $ be the ${\mathfrak{so}}(d)$--valued connection
form on $O(M) $ defined by $\nabla $. Explicitly, if $s\rightarrow u(s)$ is
a smooth path in $O(M)$ then $\omega (u^{\prime }(0)):=u(0)^{-1}\nabla
u(s)/ds|_{t=0}, $ where $\nabla u(s)/ds$ is defined as in equation (\ref
{e.1.10}) with $X$ replaced by $u$. The forms $(\theta ,\omega )$ satisfy
the structure equations

\begin{subequations}
\label{e.2.1} 
\begin{align}
d\theta & =-\omega \wedge \theta ,  \label{e.2.1a} \\
d\omega & =-\omega \wedge \omega +\Omega  \label{e.2.1b}
\end{align}
where %$\Theta$ is the $\Re^d$--valued torsion 2--form and  
$\Omega $ is the ${\mathfrak{so}}(d)$--valued curvature 2--form on $O(M)$. 
%Denote by $(u,g)\mapsto ug$  the right
%action of $\OO(d)$ on $\OO(M)$ and for $A\in {\mathfrak{so}}(d)$, let 
%let $u \cdot A\in T_u \OO(M)$ denote the vertical tangent vector 
%defined by $u\cdot A :=\frac{d}{dr} |_{0}ue^{rA}$.
%With this notation, $\omega $ satisfies 
%\begin{equation}
%\omega_{ug}(\xi g)  =g^{-1}\omega_u(\xi)g  ,\qquad \omega_u (u \cdot A )= A
%\label{e.2.6}
%\end{equation}
%for $g\in \OO(d),A \in {\mathfrak{so}}(d)$ and $\xi \in T_u\OO(M)$. 
The horizontal lift $\mathcal{H}_{u}:T_{\pi (u)}M\rightarrow T_{u}O(M)$ is
uniquely defined by 
\begin{equation}
\theta \mathcal{H}_{u}u={\text{id}}_{\mathbb{R}^{d}},\quad \omega _{u}%
\mathcal{H}_{u}=0.  \label{e.2.2}
\end{equation}
%where $\omega _{u}:=\omega |_{T_{u}\OO(M)}$ is the restriction of $\omega $
%to $T_{u}\OO(M).$
\end{subequations}

\begin{definition}
\label{d.2.1}The \textbf{curvature tensor }$R$ of $\nabla $ is 
\begin{equation}
R(X,Y)Z=\nabla _{X}\nabla _{Y}Z-\nabla _{X}\nabla _{Y}Z-\nabla _{\lbrack
X,Y]}Z  \label{e.2.3}
\end{equation}
for all vector fields $X,Y$ and $Z$ on $M.$ The \textbf{Ricci tensor} of $%
(M,g)$ is $\mathrm{Ric}X:=\sum_{i=1}^{d}R(X,e_{i})e_{i}$ and the \textbf{%
scalar curvature }$\mathrm{Scal}$ is $\mathrm{Scal}=\sum_{i=1}^{d}\langle 
\mathrm{Ric}e_{i},e_{i}\rangle ,$ where $\{e_{i}\}$ is an orthonormal frame.
\end{definition}

The relationship between \ $\Omega $ and $R$ is: 
\begin{equation}
\Omega (\xi ,\eta )=u^{-1}R(\pi _{\ast }\xi ,\pi _{\ast }\eta )u=\Omega (%
\mathcal{H}_{u}\pi _{\ast }\xi ,\mathcal{H}_{u}\pi _{\ast }\eta )
\label{e.2.4}
\end{equation}
for all $u\in O(M)$ and $\xi ,\eta \in T_{u}O(M).$ The second equality in
equation (\ref{e.2.4}) follows from the fact that $\Omega $ is horizontal,
i.e. $\Omega (\xi ,\eta )$ depends only on the horizontal components of $\xi 
$ and $\eta .$

\subsection{Path spaces and the development map\label{s.2.2}}

Let $(M,o,\langle \cdot ,\cdot \rangle ,\nabla ),$ $(O(M),u_{0}),$ $\mathrm{W%
}(M),$ and $\mathrm{H}(M)$ be as above. \ We also let $\mathrm{H}(O(M))$ be
the set of finite energy paths $u:[0,1]\longrightarrow O(M)$ as defined in
equation (\ref{e.1.8}) with $M$ replaced by $O(M)$ and $o$ by $u_{0}.$

For $\sigma \in \mathrm{H}(M)$, let $s\mapsto u(s)$ be the horizontal lift
of $\sigma $ starting at $u_{0}$, i.e. $u$ is the solution of the ordinary
differential equation 
\begin{equation*}
u^{\prime }(s)=\mathcal{H}_{u(s)}\sigma ^{\prime }(s),\qquad u(0)=u_{0}.
\end{equation*}
Notice that this equation implies that $\omega (u^{\prime }(s))=0$ or
equivalently that $\nabla u(s)/ds=0$. Hence $u(s)=//_{s}(\sigma )u_{0}$,
where as before $//_{s}(\sigma )$ is the parallel translation operator along 
$\sigma $. Again since $u_{0}\in O(M)$ is fixed in this paper we will use $%
u_{0}$ to identify $T_{o}M$ with $\mathbb{R}^{d}$ and simply write $%
u(s)=//_{s}(\sigma )$. By smooth dependence of solutions of ordinary
differential equations on parameters, the map $\sigma \in \mathrm{H}%
(M)\mapsto //(\sigma )\in \mathrm{H}(O(M))$ is smooth. A proof of this fact
may be given using the material in Palais \cite{Palais63}, see also
Corollary 4.1 in \cite{Driver89}.

\begin{definition}[Cartan's Development Map]
\label{d.2.2}The \textbf{development map} $\phi :\mathrm{H}(\mathbb{R}%
^{d})\rightarrow \mathrm{H}(M)$ is defined, for $b\in \mathrm{H},$ by $\phi
(b)=\sigma \in \mathrm{H}(M)$ where $\sigma $ solves the functional
differential equation: 
\begin{equation}
\sigma ^{\prime }(s)=//_{s}(\sigma )b^{\prime }(s),\qquad \sigma (0)=o,
\label{e.2.5}
\end{equation}
see \cite{Cartan92,Kobayashi57,Eells:Elworthy:SI}.\qed
\end{definition}

It will be convenient to give another description of the development map $%
\phi $. Namely, if $b\in \mathrm{H}(\mathbb{R}^{d})$ and $\sigma =\phi
(b)\in \mathrm{H}(M)$ as defined in equation (\ref{e.2.5}) then $\sigma =\pi
(w)$ where $w(s)\in O(M)$ is the unique solution to the ordinary
differential equation 
\begin{equation}
w^{\prime }(s)=\mathcal{H}_{w(s)}w(s)b^{\prime }(s),\qquad w(0)=u_{0}.
\label{e.2.6}
\end{equation}
>From this description of $\phi $ and smooth dependence of solutions of
ordinary differential equations on parameters it can be seen that $\phi :%
\mathrm{H}(\mathbb{R}^{d})\rightarrow \mathrm{H}(M)$ is smooth. Furthermore, 
$\phi $ is injective by uniqueness of solutions to ordinary differential
equations.

The \textbf{anti--development map} $\phi ^{-1}:\mathrm{H}(M)\rightarrow 
\mathrm{H}(\mathbb{R}^{d})$ is given by $b=\phi ^{-1}(\sigma )$ where 
\begin{equation}
b(s)=\int_{0}^{s}//_{r}{}^{-1}(\sigma )\sigma ^{\prime }(r)dr.  \label{e.2.7}
\end{equation}
This inverse map $\phi ^{-1}$ is injective and smooth by the same arguments
as above. Hence $\phi :\mathrm{H}(\mathbb{R}^{d})\rightarrow \mathrm{H}(M)$
is a diffeomorphism of infinite dimensional Hilbert manifolds, see \cite
{Eells:Elworthy:SI}. % and Theorem \QTSN{ref}{t.9.4} of the Appendix.
However, as can be seen from equation (\ref{e.3.5}) \ below, $\phi $ is not
an isometry of the Riemannian manifolds $\mathrm{H}(M)$ and $\mathrm{H}(%
\mathbb{R}^{d})$ unless the curvature $\Omega $ of $M$ is zero. So the
geometry of $\mathrm{H}(\mathbb{R}^{d})$ and that of $(\mathrm{H}(M),G^{1})$
are not well related by $\phi $.

For each $h\in C^{\infty }(\mathrm{H}(M)\rightarrow H)$ and $\sigma \in 
\mathrm{H}(M)$, let $X^{h}(\sigma )\in T_{\sigma }\mathrm{H}(M)$ be given by 
\begin{equation}
X_{s}^{h}(\sigma ):=//_{s}(\sigma )h_{s}(\sigma )\text{ for all }s\in I,
\label{e.2.8}
\end{equation}
where for notational simplicity we have written $h_{s}(\sigma )$ for $%
h(\sigma )(s).$ The vector field 
%%%%As is verified in Corollary \ref{c.9.3} of the Appendix,%%
$X^{h}$ is a smooth vector field on $\mathrm{H}(M)$ for all $h\in \mathrm{H}$%
. The reader should also note that the map 
\begin{equation}
((\sigma ,h)\rightarrow X^{h}(\sigma )):\mathrm{H}(M)\times \mathrm{H}%
\rightarrow T\mathrm{H}(M)  \label{e.2.9}
\end{equation}
is an isometry of vector bundles.

\section{Differentials of the Development Map\label{s.3}}

For $u\in O(M)$ and $v,w\in T_{\pi (u)}M,$ let 
\begin{equation*}
R_{u}(v,w)=\Omega (\mathcal{H}_{u}v,\mathcal{H}_{u}w)=u^{-1}R(v,w)u
\end{equation*}
and for $a,b\in \mathbb{R}^{d}$ let 
\begin{equation*}
\Omega _{u}(a,b):=\Omega (\mathcal{H}_{u}ua,\mathcal{H}%
_{u}ub)=u^{-1}R(ua,ub)u.
\end{equation*}
For $\sigma \in \mathrm{H}(M)$ and $X\in T_{\sigma }\mathrm{H}(M)$, define $%
q_{s}(X)\in {\mathfrak{so}}(d)$ by 
\begin{equation}
q_{s}(X)=\int_{0}^{s}R_{u(r)}(\sigma ^{\prime }(r),X(r))dr.  \label{e.3.1}
\end{equation}
where $u=//(\sigma )$ is the horizontal lift of $\sigma $. 
%Then for $\sigma \in \HCM(M)$, $q: T_{\sigma} \HCM(M) \to \HCM({\mathfrak{so}}(d))$.

\begin{rem}
The one form $q_{s}$ in equation (\ref{e.3.1}) naturally appears as soon as
one starts to compute the differential of parallel translation operators,
see for example Theorem 2.2 in Gross \cite{Gross85} and Theorem 4.1 in \cite
{Driver89} and Theorem \ref{t.3.3} below.\qed
\end{rem}

\begin{notation}
\label{n.3.2}Given $A\in {\mathfrak{so}}(d)$ and $u\in O(M),$ let $u\cdot
A\in T_{u}O(M)$ denote the vertical tangent vector defined by $u\cdot A:=%
\frac{d}{dr}|_{0}ue^{rA}.$\qed
\end{notation}

\begin{thm}
\label{t.3.3} Let $\sigma \in \mathrm{H}(M)$, let $u=//(\sigma )$ be the
horizontal lift of $\sigma $ and let $b=\phi ^{-1}(\sigma )$. Then for $X\in
T_{\sigma }\mathrm{H}(M)$, 
\begin{align}
(//_{s}^{\ast }\omega )(X)& =q_{s}(X),  \label{e.3.2} \\
(//_{s}^{\ast }\theta )(X)& =u^{-1}(s)X(s),  \label{e.3.3} \\
(//_{\ast }X)(s)& =u(s)\cdot q_{s}(X)+\mathcal{H}_{u(s)}X(s),  \label{e.3.4}
\\
\intertext{and}
(\phi ^{\ast }X)(s)& =u^{-1}(s)X(s)-\int_{0}^{s}q_{r}(X)b^{\prime
}(r)dr,  \label{e.3.5}
\end{align}
where $\phi ^{\ast }X(b):=\phi _{\ast }^{-1}X(\phi (b)).$
\end{thm}

\begin{rem}
The results of this theorem may be found in one form or another in \cite
{Bismut84a,Driver89b,Driver92b,Cruzeiro96,Gross85,Leandre97c}. We will
nevertheless supply a proof to help fix our notation and keep the paper
reasonably self contained.
\end{rem}

\begin{proof}
Choose a one parameter family $t\mapsto \sigma _{t}$ of curves in $\mathrm{H}%
(M)$ such that $\sigma _{0}=\sigma $ and $\dot{\sigma}_{0}(s)=X(s)$ where $%
\dot{\sigma}_{t}(s)=\frac{d}{dt}\sigma _{t}(s)$. Let $u_{t}(s):=//_{s}(%
\sigma _{t})$, be the horizontal lift of $\sigma _{t},$ $u(s)=//_{s}(\sigma
) $, $u_{t}^{\prime }(s):=du_{t}(s)/ds,$ $\dot{u}_{t}(s):=du_{t}(s)/dt$ and $%
\dot{u}(s):=du_{t}(s)/dt|_{t=0}.$ (In general $t$--derivatives will be
denoted by a ``dot'' and $s$--derivatives will be denoted by a ``prime.'')
Notice, by definition, that 
\begin{equation*}
\dot{u}(s)=(//_{s})_{\ast }X=(//_{\ast }X)(s)
\end{equation*}
and $\omega (u_{t}^{\prime }(s))=0$ for all $(t,s).$\ The Cartan identity 
\begin{equation}
d\alpha (X,Y)=X\alpha (Y)-Y\alpha (X)-\alpha ([X,Y]),  \label{e.3.6}
\end{equation}
valid for any 1-form $\alpha $ and vector fields $X,Y,$ gives 
\begin{equation*}
0=\frac{d}{dt}\omega (u^{\prime })=d\omega (\dot{u},u^{\prime })+\frac{d}{ds}%
\omega (\dot{u})=\Omega (\dot{u},u^{\prime })+\frac{d}{ds}\omega (\dot{u}),
\end{equation*}
where we have used the structure equations (\ref{e.2.1b}) and $0=\omega
(u^{\prime })$ in the second equality. Setting $t=0$ and integrating the
previous equation relative to $s$ yields 
\begin{align*}
(//_{s}^{\ast }\omega )(X)& :=\omega ((//_{s})_{\ast }X)=\int_{0}^{s}\Omega
(u^{\prime }(0,r),\dot{u}(0,r))\,dr \\
& =\int_{0}^{s}R_{u_{0}(r)}(\pi _{\ast }u^{\prime }(0,r),\pi _{\ast }\dot{u}%
(0,r))\,dr \\
& =\int_{0}^{s}R_{u_{0}(r)}(\sigma ^{\prime }(r),X(r))\,dr,
\end{align*}
where we have made use of the fact that $\Omega $ is horizontal and the
relation $\sigma _{t}(s)=\pi (u_{t}(s)).$ This proves equation (\ref{e.3.2}%
). Equation (\ref{e.3.3}) is verified as follows: 
\begin{eqnarray*}
(//_{s}^{\ast }\theta )(X) &=&\theta ((//_{s})_{\ast }X)=\theta (\dot{u}%
(s))=u_{0}^{-1}(s)\frac{d}{dt}|_{t=0}\pi (u_{t}(s)) \\
&=&//_{s}^{-1}(\sigma )\frac{d}{dt}|_{t=0}\sigma _{t}(s)=//_{s}^{-1}(\sigma
)X(s).
\end{eqnarray*}

Recall that for $u\in O(M)$, $(\theta ,\omega ):T_{u}O(M)\rightarrow \mathbb{%
R}^{d}\times {\mathfrak{so}}(d)$ is an isomorphism. Therefore equations (\ref
{e.3.2}) and (\ref{e.3.3}) imply (\ref{e.3.4}), after taking into account
the definition of $\theta $ %and equation \eqref{e.2.6}.
and the identity, 
\begin{equation*}
\omega (u\cdot A):=u^{-1}\frac{\nabla }{dr}|_{r=0}ue^{rA}=A.
\end{equation*}

To prove equation (\ref{e.3.5}), let $b=\phi ^{-1}(\sigma )$ and $%
u(s)=//_{s}(\sigma ).$ Then 
\begin{equation*}
b(s)=\int_{0}^{s}u^{-1}(r)\sigma ^{\prime }(r)\,dr=\int_{0}^{s}\theta
(u^{\prime }(r))\,dr,
\end{equation*}
or equivalently, 
\begin{equation*}
b^{\prime }(s)=\theta (u^{\prime }(s)).
\end{equation*}
Therefore 
\begin{align*}
\frac{d}{ds}\phi _{\ast }^{-1}X(s)& =\frac{d}{dt}\theta (u_{t}^{\prime
}(s))|_{t=0} \\
& =\frac{d}{ds}\theta (\dot{u}(s))+d\theta (\dot{u}(s),u^{\prime }(s)) \\
& =\frac{d}{ds}(u^{-1}(s)X(s))-\omega \wedge \theta (\dot{u}(s),u^{\prime
}(s)) \\
& =\frac{d}{ds}(u^{-1}(s)X(s))-\omega (\dot{u}(s))\theta (u^{\prime }(s)) \\
& =\frac{d}{ds}(u^{-1}(s)X(s))-q_{s}(X)b^{\prime }(s),
\end{align*}
where we have used the equations (\ref{e.3.6}), (\ref{e.2.1a}), (\ref{e.3.2}%
) and the fact that $\omega (u^{\prime }(s))=0$. Integrating the last
equation relative to $s$ proves (\ref{e.3.5}).
\end{proof}

\subsection{Bracket Computation}

\label{s.3.1}

\begin{thm}[Lie Brackets]
\label{t.3.5} Let $h,k:\mathrm{H}(M)\rightarrow \mathrm{H}(\mathbb{R}^{d})$
be smooth functions. (We will write $h_{s}(\sigma )$ for $h(\sigma )(s).)$
Then $[X^{h},X^{k}]=X^{f(h,k)},$ where $f(h,k)$ is the smooth function $%
\mathrm{H}(M)\rightarrow \mathrm{H}(\mathbb{R}^{d})$ defined by 
\begin{equation}
f_{s}(h,k)(\sigma ):=X^{h}(\sigma )k_{s}-X^{k}(\sigma
)h_{s}+q_{s}(X^{k}(\sigma ))h_{s}-q_{s}(X^{h}(\sigma ))k_{s},  \label{e.3.7}
\end{equation}
where $q=//^{\ast }\omega $ as in equation (\ref{e.3.2}) and $X^{h}(\sigma
)k_{s}$ denotes derivative of $\sigma \rightarrow k_{s}(\sigma )$ by the
tangent vector $X^{h}(\sigma ).$
\end{thm}

\begin{rem}
This theorem also appears in equation (1.32) in Leandre \cite{Leandre94c},
equation (6.2.2) in Cruzeiro and Malliavin \cite{Cruzeiro96} and is Theorem
6.2 in \cite{Driver97a}. To some extent it is also contained in \cite
{Flaschel:Klingenberg}. Again for the readers convenience will supply a
short proof.\qed
\end{rem}

\begin{proof}
The vector fields $X^{h}$ and $X^{k}$ on $\mathrm{H}(M)$ are smooth, hence $%
[X^{h},X^{k}]$ is well defined. In order to simplify notation, we will
suppress the arguments $\sigma $ and $s$ from the proof of equation (\ref
{e.3.7}).

According to equation (\ref{e.3.3}), $h=(//^{\ast }\theta )(X^{h})$, $%
k=(//^{\ast }\theta )(X^{k})$, and $f(h,k)=(//^{\ast }\theta
)([X^{h},X^{k}]) $. Using equations (\ref{e.3.1}--\ref{e.3.6}) we find that 
\begin{eqnarray*}
f(h,k) &=&X^{h}\left[ (//^{\ast }\theta )(X^{k})\right] -X^{k}\left[
(//^{\ast }\theta )(X^{h})\right] -(d(//^{\ast }\theta ))(X^{h},X^{k}) \\
&=&X^{h}k-X^{k}h-(//^{\ast }d\theta )(X^{h},X^{k}) \\
&=&X^{h}k-X^{k}h+(//^{\ast }(\omega \wedge \theta ))(X^{h},X^{k}) \\
&=&X^{h}k-X^{k}h+(//^{\ast }\omega \wedge //^{\ast }\theta )(X^{h},X^{k}) \\
&=&X^{h}k-X^{k}h+q(X^{h})k-q(X^{k})h.
\end{eqnarray*}
\end{proof}

%\begin{rem}
%\label{r.3.6} Let $q=//^{*}\omega $ as in the above theorem, then 
%$$
%dq+q\wedge q=//^{*}(d\omega +\omega \wedge \omega )=//^{*}\Omega . 
%$$
%\end{rem}

\section{Finite Dimensional Approximations\label{s.4}}

\begin{definition}
\label{d.4.1}Let $\mathcal{P}=\{0=s_{0}<s_{1}<s_{2}<\cdots <s_{n}=1\}$ be a 
\textbf{partition} of $[0,1]$ and let $|\mathcal{P}|=\max_{i}|s_{i}-s_{i-1}|$
be the\textbf{\ norm} of the partition, $J_{i}:=(s_{i-1},s_{i}]$ for $%
i=1,2,\ldots ,n.$ For a function $k,$ let $\Delta _{i}k:=k(s_{i})-k(s_{i-1})$
and $\Delta _{i}s=s_{i}-s_{i-1}.$ For a piecewise continuous function on $%
[0,1]$, we will use the notation $f(s+)=\lim_{r\searrow s}f(r)$.\qed
\end{definition}

\begin{notation}
\label{n.4.2}$\mathrm{H}_{\mathcal{P}}=\{x\in \mathrm{H}\cap
C^{2}(I\setminus \mathcal{P}):x^{\prime \prime }(s)=0\text{\textrm{\ for }}%
s\notin \mathcal{P}\}$ --- the piecewise linear paths in $\mathrm{H}:=%
\mathrm{H}(\mathbb{R}^{d})$, which change directions only at the partition
points.\qed

\begin{rem}[Development]
\label{r.4.3} The development map $\phi :\mathrm{H}\rightarrow \mathrm{H}(M)$
has the property that $\phi (\mathrm{H}_{\mathcal{P}})=\mathrm{H}_{\mathcal{P%
}}(M),$ where $\mathrm{H}_{\mathcal{P}}(M)$ has been defined in Notation \ref
{n.1.4} above. Indeed, if $\sigma =\phi (b)$ with $b\in \mathrm{H}_{\mathcal{%
P}},$ then differentiating equation $(\ref{e.2.5})$ gives: 
\begin{equation*}
\frac{\nabla \sigma ^{\prime }(s)}{ds}=\frac{\nabla }{ds}(//_{s}(\sigma
)b^{\prime }(s))=//_{s}(\sigma )b^{\prime \prime }(s)=0\text{ for all }%
s\notin \mathcal{P}.
\end{equation*}
We will write $\phi _{\mathcal{P}}$ for $\phi |_{\mathrm{H}_{\mathcal{P}}}.$%
\qed
\end{rem}
\end{notation}

Because $\phi :\mathrm{H}\rightarrow \mathrm{H}(M)$ is a diffeomorphism and $%
\mathrm{H}_{\mathcal{P}}\subset \mathrm{H}$ is an embedded submanifold, it
follows that $\mathrm{H}_{\mathcal{P}}(M)$ is an embedded submanifold of $%
\mathrm{H}(M)$. Therefore for each $\sigma \in \mathrm{H}_{\mathcal{P}}(M)$, 
$T_{\sigma }\mathrm{H}_{\mathcal{P}}(M)$ may be viewed as a subspace of $%
T_{\sigma }\mathrm{H}(M)$. The next proposition explicitly identifies this
subspace.

\begin{prop}[Tangent Space]
\label{p.4.4} Let $\sigma \in \mathrm{H}_{\mathcal{P}}(M),$ then $X\in
T_{\sigma }\mathrm{H}(M)$ is in $T_{\sigma }\mathrm{H}_{\mathcal{P}}(M)$ if
and only if 
\begin{equation}
\frac{\nabla ^{2}}{ds^{2}}X(s)=R(\sigma ^{\prime }(s),X(s))\sigma ^{\prime
}(s)\text{\textrm{\ on }}I\setminus \mathcal{P}.  \label{e.4.1}
\end{equation}
Equivalently, letting $b=\phi ^{-1}(\sigma ),$ $u=//(\sigma )$ and $h\in 
\mathrm{H},$ then $X^{h}\in T_{\sigma }\mathrm{H}(M)$ defined in equation (%
\ref{e.2.8}) is in $T_{\sigma }\mathrm{H}_{\mathcal{P}}(M)$ if and only if 
\begin{equation}
h^{\prime \prime }(s)=\Omega _{u(s)}(b^{\prime }(s),h(s))b^{\prime }(s)\text{%
\textrm{\ on }}I\setminus \mathcal{P}.  \label{e.4.2}
\end{equation}
\end{prop}

\begin{proof}
Since $\mathrm{H}_{\mathcal{P}}(M)$ consists of piecewise geodesics, it
follows that for $\sigma \in \mathrm{H}_{\mathcal{P}}(M)$, any $X\in
T_{\sigma }\mathrm{H}_{\mathcal{P}}(M)$ must satisfy the Jacobi equation (%
\ref{e.4.1}) for $s\notin \mathcal{P}$. Equation \ref{e.4.2} is a
straightforward reformulation of this using the definitions.

It is instructive to give a direct proof of equation (\ref{e.4.2}). Since $%
\mathrm{H}_{\mathcal{P}}$ is a vector space, $T_{b}\mathrm{H}_{\mathcal{P}%
}\cong \mathrm{H}_{\mathcal{P}}$ for all $b\in \mathrm{H}_{\mathcal{P}}$.
Since $\phi _{\mathcal{P}}:\mathrm{H}_{\mathcal{P}}\rightarrow \mathrm{H}_{%
\mathcal{P}}(M)$ is a diffeomorphism, we must identify those vectors $X\in
T_{\sigma }\mathrm{H}(M)$ such that $\phi ^{\ast }X\in \mathrm{H}_{\mathcal{P%
}}$, i.e. those $X$ such that $(\phi ^{\ast }X)^{\prime \prime }:=0$ on $%
I\setminus \mathcal{P}$. Because $b\in \mathrm{H}_{\mathcal{P}}$ and hence $%
b^{\prime \prime }(s)=0$ on $I\setminus \mathcal{P},$ it follows from
equation (\ref{e.3.5}) that $(\phi ^{\ast }X)^{\prime \prime }=0$ 
%(\phi^* X)(s) = h(s) - \int_0^s  q_r(X)b'(r) dr .
%$$
%Demanding that $\frac{d^2}{ds^2}(\phi^* X) = 0$ 
on $I\setminus \mathcal{P}$ is equivalent to 
\begin{equation*}
0=h^{\prime \prime }(s)-\Omega _{u(s)}(b^{\prime }(s),h(s))b^{\prime }(s)%
\text{ on }I\setminus \mathcal{P}.
\end{equation*}
\end{proof}

\begin{rem}
\label{r.4.5}The metric $G_{\mathcal{P}}^{1}$ in Definition \ref{d.1.5}
above is easily seen to be non-degenerate because if $G_{\mathcal{P}%
}^{1}(X,X)=0$ then $\nabla X(s_{i}+)/ds=0$ for all $i.$ It then follows from
the continuity of $X$ and the fact that $X$ solves the Jacobi equation (\ref
{e.4.1}) that $X$ is zero. Also note that $G_{\mathcal{P}}^{1}$ is a
``belated'' Riemann sum approximation to the metric on $\mathrm{H}_{\mathcal{%
P}}(M)$ which is inherited from $G^{1}$ on $\mathrm{H}(M).$ Moreover, in the
case $M=\mathbb{R}^{d},$ the metric $G_{\mathcal{P}}^{1}$ is equal to $G^{1}$
on $T\mathrm{H}_{\mathcal{P}}(M).$\qed
\end{rem}

\begin{definition}
\label{d.4.6}Let $\mathrm{Vol}_{\mathcal{P}}$ be the Riemannian volume form
on $\mathrm{H}_{\mathcal{P}}$ equipped with the $H^{1}$ -- metric, $%
(h,k):=\int_{0}^{1}\langle h^{\prime }(s),k^{\prime }(s)\rangle ds.$\qed
\end{definition}

\begin{notation}
\label{n.4.7} Let $\mathcal{P}=\{0=s_{0}<s_{1}<s_{2}<\cdots <s_{n}=1\}$ be a
partition of $[0,1]$. For each $i=1,2,\ldots n,$ and $s\in (s_{i-1},s_{i}]$,
define 
\begin{align}
\hat{q}_{s}^{\mathcal{P}}(X)& =q_{s_{i-1}}(X)  \label{e.4.3} \\
\intertext{and}
\tilde{q}_{s}^{\mathcal{P}}(X)=q_{s}(X)-q_{s_{i-1}}(X)&
=\int_{s_{i-1}}^{s}\Omega _{u}(\sigma ^{\prime }(r),X(r))dr.  \label{e.4.4}
\end{align}
\end{notation}

Note that $q=\hat{q}^{\mathcal{P}}+\tilde{q}^{\mathcal{P}}$ and hence
equation (\ref{e.3.5}) becomes 
\begin{equation}
(\phi ^{\ast }X^{h})^{\prime }(s)=h^{\prime }(s)-q_{s}(X^{h})b^{\prime
}(s)=h^{\prime }(s)-\hat{q}_{s}(X^{h})b^{\prime }(s)-\tilde{q}%
_{s}(X^{h})b^{\prime }(s)  \label{e.4.5}
\end{equation}
for all $h\in \mathrm{H}(\mathbb{R}^{d}).$\qed

\begin{thm}
\label{t.4.8} $\phi _{\mathcal{P}}^{\ast }\mathrm{Vol}_{G_{\mathcal{P}}^{1}}=%
\mathrm{Vol}_{\mathcal{P}}$
\end{thm}

\begin{proof}
Let $\{h_{k}\}$ be an orthonormal basis for $\mathrm{H}_{\mathcal{P}}$, $%
b\in \mathrm{H}_{\mathcal{P}}$, $\sigma =\phi (b)$ and $u=//(\sigma )$.
Using the definitions of the volume form on a Riemannian manifold we must
show that 
\begin{equation*}
\det (G_{\mathcal{P}}^{1}(\phi _{\ast }h_{k},\phi _{\ast }h_{j}))=1,
\end{equation*}
where $\phi _{\ast }h_{k}:=\frac{d}{dt}|_{0}\phi (b+th_{k}).$

Let $H_{k}(s)=u^{-1}(s)(\phi _{\ast }(h_{k}))(s)$ and set 
\begin{equation*}
\langle H,K\rangle _{\mathcal{P}}:=\sum_{i=1}^{n}\langle H^{\prime
}(s_{i-1}+),K^{\prime }(s_{i-1}+)\rangle \Delta _{i}s.
\end{equation*}
Then $X^{H_{k}}=\phi _{\ast }(h_{k})$ and 
\begin{equation*}
\det (G_{\mathcal{P}}^{1}(\phi _{\ast }(h_{k}),\phi _{\ast }(h_{j})))=\det
(\langle H_{k},H_{j}\rangle _{\mathcal{P}}).
\end{equation*}
By equation (\ref{e.4.5}) 
\begin{equation*}
h_{k}^{\prime }=(\phi ^{\ast }(X^{H_{k}}))^{\prime }=H_{k}^{\prime
}-q(X^{H_{k}})b^{\prime }=H_{k}^{\prime }-\hat{q}(X^{H_{k}})b^{\prime }-%
\tilde{q}(X^{H_{k}})b^{\prime }
\end{equation*}
so that 
\begin{equation}
h_{k}^{\prime }+\hat{q}(X^{H_{k}})b^{\prime }=H_{k}^{\prime }-\tilde{q}%
(X^{H_{k}})b^{\prime }.  \label{e.4.6}
\end{equation}
Noting that $h_{k}^{\prime }$, $\hat{q}(X^{H_{k}})$, and $b^{\prime }$ are
all constant on $(s_{i-1},s_{i})$ and that $\tilde{q}%
_{s_{i-1}}(X^{H_{k}})=0, $ it follows that both sides of equation (\ref
{e.4.6}) are constant on $(s_{i-1},s_{i})$ and the constant value is $%
H_{k}^{\prime }(s_{i-1}+).$ Therefore 
\begin{align*}
\langle H_{k},H_{j}\rangle _{\mathcal{P}}& =\int_{0}^{1}\langle
H_{k}^{\prime }-\tilde{q}(X^{H_{k}})b^{\prime },H_{j}^{\prime }-\tilde{q}%
(X^{H_{j}})b^{\prime }\rangle \,ds \\
& =\int_{0}^{1}\langle h_{k}^{\prime }+\hat{q}(X^{H_{k}})b^{\prime
},h_{j}^{\prime }+\hat{q}(X^{H_{j}})b^{\prime }\rangle \,ds.
\end{align*}
Define the linear transformation, $T:\mathrm{H}_{\mathcal{P}}\rightarrow 
\mathrm{H}_{\mathcal{P}}$ by 
\begin{equation*}
(Th)(s)=\int_{0}^{s}\hat{q}_{r}(\phi _{\ast }h)b^{\prime }(r)dr.
\end{equation*}
We have just shown that 
\begin{align*}
\det (G_{\mathcal{P}}^{1}(\phi _{\ast }(h_{k}),\phi _{\ast }(h_{j})))& =\det
(\left\{ \langle (I+T)h_{k},(I+T)h_{j}\rangle _{\mathcal{P}}\right\} _{j,k})
\\
& =\det (\left\{ \langle h_{k},(I+T)^{\ast }(I+T)h_{j}\rangle _{\mathcal{P}%
}\right\} _{j,k}) \\
& =\det ((I+T)^{\ast }(I+T))=[\det (I+T)]^{2}.
\end{align*}
So to finish the proof it suffices to show that $\det (I+T)=1.$ This will be
done by showing that $T$ is nilpotent. For this we will make a judicious
choice of orthonormal basis for $\mathrm{H}_{\mathcal{P}}$. Let $%
\{e_{a}\}_{a=1}^{d}$ be an orthonormal basis for $T_{o}M\cong \mathbb{R}^{d}$
and define 
\begin{equation*}
h_{i,a}(s)=\left( \frac{1}{\sqrt{\Delta _{i}s}}\int_{0}^{s}1_{J_{i-1}}(r)dr%
\right) e_{a}
\end{equation*}
for $i=1,2,\ldots ,n$, $a=1,\ldots ,d$. Using the causality properties of $%
\phi $ and $\hat{q}$, it follows that $\phi _{\ast }h_{i,a}:=0$ on $%
[0,s_{i-1}]$ and hence $\hat{q}(\phi _{\ast }(h_{i,a})):=0$ on $[0,s_{i}).$
Thus for any $a,b$, $\langle Th_{i,a},h_{j,b}\rangle =0$ if $j\leq i$. This
shows that $T$ is nilpotent and hence finishes the proof.
\end{proof}

\begin{definition}
\label{d.4.9}Let $E_{\mathbb{R}^{d}}(b):=\int_{0}^{1}|b^{\prime }(s)|^{2}ds$
denote the energy of a path $b\in \mathrm{H.}$ For each partition $\mathcal{P%
}=\{0=s_{0}<s_{1}<s_{2}<\cdots <s_{n}=1\}$ of $[0,1],$ let $\mu _{\mathcal{P}%
}^{1}$ denote the volume form 
\begin{equation*}
\mu _{\mathcal{P}}^{1}=\frac{1}{Z_{\mathcal{P}}^{1}}e^{-\frac{1}{2}E_{%
\mathbb{R}^{d}}}\mathrm{Vol}_{\mathrm{H}_{\mathcal{P}}}
\end{equation*}
on $\mathrm{H}_{\mathcal{P}}$, where $Z_{\mathcal{P}}^{1}:=(2\pi )^{dn/2}.$
(By Lemma \ref{l.4.11} below, $\mu _{\mathcal{P}}^{1}$ is a probability
measure on $\mathrm{H}_{\mathcal{P}}.)$ \qed
\end{definition}

Let $b\in \mathrm{H}$ and $\sigma :=\phi (b)\in \mathrm{H}(M)$. Because
parallel translation is an isometry, it follows from equation (\ref{e.2.5})
that $E(b)=E(\sigma )$. As an immediate consequence of this identity and
Theorem \ref{t.4.8} is the following theorem.

\begin{thm}
\label{t.4.10} Let $\mu _{\mathcal{P}}^{1}$ (Definition \ref{d.4.9}) and $%
\nu _{\mathcal{P}}^{1}$ (Definition \ref{d.1.7}) be as above, then $\mu _{%
\mathcal{P}}^{1}$ is the pull back of $\nu _{\mathcal{P}}^{1}$ by $\phi _{%
\mathcal{P}},$ i.e. $\mu _{\mathcal{P}}^{1}=\phi _{\mathcal{P}}^{\ast }\nu _{%
\mathcal{P}}^{1}$.\qed
\end{thm}

Before exploring the consequences of this last theorem, we will make a few
remarks about the measure $\mu _{\mathcal{P}}^{1}$. Let $\pi _{\mathcal{P}}:%
\mathrm{W}(\mathbb{R}^{d})\rightarrow (\mathbb{R}^{d})^{n}$ be given by $\pi
_{\mathcal{P}}(x):=(x(s_{1}),x(s_{2}),\ldots ,x(s_{n}))$. Note that $\pi _{%
\mathcal{P}}:\mathrm{H}_{\mathcal{P}}\rightarrow (\mathbb{R}^{d})^{n}$ is a
linear isomorphism of finite dimensional vector spaces. We will denote the
inverse of $\pi _{\mathcal{P}}|_{\mathrm{H}_{\mathcal{P}}}$ by $i_{\mathcal{P%
}}$.

\begin{lemma}
\label{l.4.11} Let $dy_{1}dy_{2}\cdots dy_{n}$ denote the standard volume
form on $(\mathbb{R}^{d})^{n}$ and $y_{0}:=0$ by convention. Then 
\begin{equation}
i_{\mathcal{P}}^{\ast }\mu _{\mathcal{P}}^{1}=\frac{1}{Z_{\mathcal{P}}^{1}}%
\left( \prod_{i=1}^{n}(\Delta _{i}s)^{-d/2}\exp \{-\frac{1}{2\Delta _{i}s}%
|y_{i}-y_{i-1}|^{2}\}\right) dy_{1}dy_{2}\cdots dy_{n}  \label{e.4.7}
\end{equation}
where $Z_{\mathcal{P}}^{1}$ is defined in equation (\ref{e.1.15}). Using the
explicit value on $Z_{\mathcal{P}}^{1}$, this equation may also be written
as 
\begin{equation}
i_{\mathcal{P}}^{\ast }\mu _{\mathcal{P}}^{1}=\left(
\prod_{i=1}^{n}p_{\Delta _{i}s}(y_{i-1},y_{i})\right) dy_{1}dy_{2}\cdots
dy_{n},  \label{e.4.8}
\end{equation}
where $p_{s}(x,y):=(2\pi s)^{-d/2}\exp \{-|x-y|^{2}/2s\}$ is the heat kernel
on $\mathbb{R}^{d}.$ In particular $i_{\mathcal{P}}^{\ast }\mu _{\mathcal{P}%
}^{1}$ and hence $\mu _{\mathcal{P}}^{1}$ are probability measures.
\end{lemma}

\begin{proof}
Let $x\in \mathrm{H}_{\mathcal{P}}$, then 
\begin{equation*}
E(x)=\int_{0}^{1}|x^{\prime }(s)|^{2}ds=\sum_{i=1}^{n}|\frac{\Delta _{i}x}{%
\Delta _{i}s}|^{2}\Delta _{i}s=\sum_{i=1}^{n}\frac{1}{\Delta _{i}s}|\Delta
_{i}x|^{2}.
\end{equation*}
Hence if $x=i_{\mathcal{P}}(y),$ then 
\begin{equation}
\int_{0}^{1}|x^{\prime }(s)|^{2}ds=\sum_{i=1}^{n}\frac{1}{\Delta _{i}s}%
|y_{i}-y_{i-1}|^{2}=\sum_{i=1}^{n}|\xi _{i}|^{2}.  \label{e.4.9}
\end{equation}
where $\xi _{i}:=(\Delta _{i}s)^{-1/2}(y_{i}-y_{i-1}).$ This last equation
shows that the linear transformation 
\begin{equation*}
x\in \mathrm{H}_{\mathcal{P}}\rightarrow \{(\Delta
_{i}s)^{-1/2}(x(s_{i})-x(s_{i-1})\}_{i=1}^{n}\in (\mathbb{R}^{d})^{n}
\end{equation*}
is an isometry of vector spaces and therefore 
\begin{equation}
i_{\mathcal{P}}^{\ast }\mathrm{Vol}_{\mathcal{P}}=d\xi _{1}d\xi _{2}\cdots
d\xi _{n}.  \label{e.4.10}
\end{equation}
Now an easy computation shows that 
\begin{equation}
d\xi _{1}d\xi _{2}\cdots d\xi _{n}=\left( \prod_{i=1}^{n}(\Delta
_{i}s)^{-d/2}\right) dy_{1}dy_{2}\cdots dy_{n}.  \label{e.4.11}
\end{equation}
>From equations ((\ref{e.4.9}) -- (\ref{e.4.11})), we see that equation \ref
{e.4.7} is valid.
\end{proof}

\begin{notation}
\label{n.4.12}Let $\{B(s)\}_{s\in \lbrack 0,1]}$ be the standard $\mathbb{R}%
^{d}$---valued Brownian motion on $(\mathrm{W}(\mathbb{R}^{d}),\mu )$ as in
Notation \ref{n.1.2}. Given a partition $\mathcal{P}$ of $[0,1]$ as above,
set $B_{\mathcal{P}}:=i_{\mathcal{P}}\circ \pi _{\mathcal{P}}(B).$ The
explicit formula for $B_{\mathcal{P}}$ is: 
\begin{equation*}
B_{\mathcal{P}}(s)=B(s_{i-1})+(s-s_{i-1})\frac{\Delta _{i}B}{\Delta _{i}s}%
\text{\textrm{\ if }}s\in (s_{i-1},s_{i}],
\end{equation*}
where $\Delta _{i}B:=B(s_{i})-B(s_{i-1})$. We will also denote the
expectation relative to $\mu $ by $\mathbb{E},$ so that $\mathbb{E}\left[ f%
\right] =\int_{\mathrm{W}(\mathbb{R}^{d})}fd\mu .$ \qed
\end{notation}

Note that $B_{\mathcal{P}}$ is the unique element in $\mathrm{H}_{\mathcal{P}%
}$ such that $B_{\mathcal{P}}=B$ on $\mathcal{P}$. We now have the following
easy corollary of Lemma \ref{l.4.11} and the fact that the right side of
equation (\ref{e.4.8}) is the distribution of $(B(s_{1}),B(s_{2}),\ldots
,B(s_{n})).$

\begin{cor}
\label{c.4.13} The law of $B_{\mathcal{P}}$ and the law of $\phi (B_{%
\mathcal{P}})$ (with respect to $\mu )$ is $\mu _{\mathcal{P}}^{1}$ and $\nu
_{\mathcal{P}}^{1}$ respectively.\qed
\end{cor}

\subsection{Limits of the finite dimensional approximations\label{s.4.1}}

Let us recall the following Wong and Zakai type approximation theorem for
solutions to Stratonovich stochastic differential equations.

\begin{thm}
\label{t.4.14} Let $f:\mathbb{R}^{d}\times \mathbb{R}^{n}\rightarrow $End$(%
\mathbb{R}^{d},\mathbb{R}^{n})$ and $f_{0}:\mathbb{R}^{d}\times \mathbb{R}%
^{n}\rightarrow \mathbb{R}^{n}$ be twice differentiable with bounded
continuous derivatives. Let $\xi _{0}\in \mathbb{R}^{n}$ and $\mathcal{P}$
be a partition of $[0,1].$ Further let $B$ and $B_{\mathcal{P}}$ be as in
Notation \ref{n.4.12} and $\xi _{\mathcal{P}}(s)$ denote the solution to the
ordinary differential equation: 
\begin{equation}
\xi _{\mathcal{P}}^{\prime }(s)=f(\xi _{\mathcal{P}}(s))B_{\mathcal{P}%
}^{\prime }(s)+f_{0}(\xi _{\mathcal{P}}(s)),\qquad \xi _{\mathcal{P}}(0)=\xi
_{0}  \label{e.4.12}
\end{equation}
and $\xi $ denote the solution to the Stratonovich stochastic differential
equation, 
\begin{equation}
d\xi (s)=f(\xi (s))\delta B(s)+f_{0}(\xi (s))ds,\qquad \xi (0)=\xi _{0}.
\label{e.4.13}
\end{equation}
Then, for any $\gamma \in (0,\frac{1}{2}),$ $p\in \lbrack 1,\infty )$, there
is a constant $C(p,\gamma )<\infty $ depending only on $f$ and $M$, so that 
\begin{equation}
\lim_{|\mathcal{P}|\rightarrow 0}\mathbb{E}\left[ \sup_{s\leq 1}|\xi _{%
\mathcal{P}}(s)-\xi (s)|^{p}\right] \leq C(p,\gamma )|\mathcal{P}|^{\gamma
p}.  \label{e.4.14}
\end{equation}
\end{thm}

This theorem is a special case of Theorem 5.7.3 and Example 5.7.4 in Kunita 
\cite{Kunita90}. Theorems of this type have a long history starting with
Wong and Zakai \cite{Wong:Zakai65b,Wong:Zakai67}. The reader may also find
this and related results in the following \emph{partial}\textbf{\ }list of
references: \cite
{Amit91,Bally89,Bally89b,Bismut81,Blum84,Doss79,Elworthy78,Guo-Shu82,Ikeda81,Jorgensen75,Kaneko:Nakao,Kurtz:Protter:1991a,Kurtz:Protter:1991b,Lyons96,Malliavin78b,Malliavin97,McShane72,McShane74,Moulinier88,Nakao:Yamato,Pinsky78,Stroock:Taniguchi1994,StVar69b,Stroock:Varadhan72,Sussmann91}%
.  The theorem as stated here may be found in \cite{DHu}.

\begin{definition}
\label{d.4.15}

\begin{enumerate}
\item  Let $u$ be the solution to the {Stratonovich stochastic differential
equation} 
\begin{equation*}
\delta u=\mathcal{H}_{u}u\delta B,\qquad u(0)=u_{0}.
\end{equation*}
Notice that $u$ may be viewed as $\mu $ -- a.e. defined function from $%
\mathrm{W}(\mathbb{R}^{d})\rightarrow \mathrm{W}(O(M)).$

\item  Let $\tilde{\phi}:=\pi \circ u:\mathrm{W}(\mathbb{R}^{d})\rightarrow 
\mathrm{W}(M).$ This map is will be called the stochastic development map.

\item  Let $/\tilde{/}_{\cdot }(\sigma )$ denote stochastic parallel
translation relative to the probability space $(\mathrm{W}(M),\nu ).$ That
is $/\tilde{/}_{\cdot }(\sigma )$ is a stochastic extension of $//_{\cdot
}(\sigma ).$

\item  Let $\tilde{b}(s)=\int_{0}^{s}/\tilde{/}_{r}{}^{-1}(\sigma )\delta
\sigma (r),$ where $\delta \sigma (r)$ denotes the Stratonovich differential.
\end{enumerate}

\qed
\end{definition}

\begin{rem}
\label{r.4.16}Using Theorem \ref{t.4.14}, one may show that $\tilde{\phi}$
is a ``stochastic extension'' of $\phi ,$ i.e. $\tilde{\phi}=\lim_{|\mathcal{%
P}|\rightarrow 0}\phi (B_{\mathcal{P}}).$ Moreover, the law of $\tilde{\phi}$
(i.e. $\mu \tilde{\phi}^{-1})$ is the Wiener measure $\nu $ on $W(M).$ It is
also well known that $\tilde{b}$ is a standard $\mathbb{R}^{d}$ -- valued
Brownian motion on $(\mathrm{W}(M),\nu )$ and that the law of $u$ under $\mu 
$ on $\mathrm{W}(\mathbb{R}^{d})$ and the law of $/\tilde{/}$ under $\nu $
are equal. \qed
\end{rem}

The fact that $\tilde{\phi}$ has a ``stochastic extension'' seems to have
first been observed by Eells and Elworthy \cite{Eells:Elworthy:SI} who used
ideas of Gangolli \cite{Gangolli1964}.The relationship of the stochastic
development map to stochastic differential equations on the orthogonal frame
bundle $O(M)$ of $M$ is pointed out in Elworthy \cite
{Elworthy74,Elworthy75,Elworthy78}. The frame bundle point of view has also
been developed by Malliavin, see for example \cite
{Malliavin78,Malliavin78b,Malliavin79}. For a more detailed history of the
stochastic development map, see pp. 156--157 in Elworthy \cite{Elworthy78}.
The results in the previous remark are all standard and may be found in the
previous references and also in \cite{Emery,Ikeda81,Kunita90,Malliavin97}.
For a fairly self contained short exposition of these results the reader may
wish to consult Section 3 in \cite{Driver92b}. Using Theorem \ref{t.4.14}
and Corollary \ref{c.4.13} above, we get the following limit theorem for $%
\nu _{\mathcal{P}}^{1}$.

\begin{thm}
\label{t.4.17} Suppose that $F:W(O(M))\rightarrow \mathbb{R}$ is a
continuous and bounded function and for $\sigma \in \mathrm{H}(M)$ we let $%
f(\sigma ):=F(//_{\cdot }(\sigma )).$ Then 
\begin{equation}
\lim_{|\mathcal{P}|\rightarrow 0}\int_{\mathrm{H}_{\mathcal{P}}(M)}f(\sigma
)d\nu _{\mathcal{P}}^{1}(\sigma )=\int_{\mathrm{W}(M)}\tilde{f}(\sigma )d\nu
(\sigma ),  \label{e.4.15}
\end{equation}
where $\tilde{f}(\sigma ):=F(/\tilde{/}_{\cdot }(\sigma )).$
\end{thm}

\begin{proof}
By Remark \ref{r.4.16} 
\begin{equation}
\int_{\mathrm{W}(M)}\tilde{f}(\sigma )d\nu (\sigma )=\mathbb{E}[\tilde{f}%
(u)].  \label{e.4.16}
\end{equation}

By embedding $O(M)$ into $\mathbb{R}^{D}$ for some $D\in \mathbb{N}$ and
extending the map $v\mapsto \mathcal{H}_{u}uv$ to a compact neighborhood of $%
O(M)\subset \mathbb{R}^{D}$, we may apply Theorem \ref{t.4.14} to conclude
that 
\begin{equation}
\lim_{|\mathcal{P}|\rightarrow 0}\mathbb{E}\left[ \sup_{0\leq s\leq 1}\left|
u_{\mathcal{P}}(s)-u(s)\right| _{^{\mathbb{R}^{D}}}^{p}\right] =0,
\label{e.4.17}
\end{equation}
where $u_{\mathcal{P}}$ solves equation (\ref{e.2.6}) with $b$ replaced by $%
B_{\mathcal{P}}.$ But the law of $u_{\mathcal{P}}$ is equal to the law of $%
//(\cdot )$ under $\nu _{\mathcal{P}}^{1}$, see Corollary \ref{c.4.13}.
Therefore, 
\begin{equation}
\int_{\mathrm{H}_{\mathcal{P}}(M)}f(\sigma )d\nu _{\mathcal{P}}^{1}(\sigma )=%
\mathbb{E}[f(u_{\mathcal{P}})].  \label{e.4.18}
\end{equation}
The limit in equation (\ref{e.4.15}) now easily follows from (\ref{e.4.16}--%
\ref{e.4.18}) and the dominated convergence theorem.
\end{proof}

\section{The $L^{2}$ metric\label{s.5}}

In section \ref{s.4} we considered the metric $G_{\mathcal{P}}^{1}$ (see
Definition \ref{d.1.5}) on $\mathrm{H}_{\mathcal{P}}(M)$ and the associated
finite dimensional approximations of the Wiener measure $\nu $ on $\mathrm{W}%
(M)$. It was found that under the development map $\phi _{\mathcal{P}}$, the
volume form with respect to. $G_{\mathcal{P}}^{1}$ pulls back to the volume
form of a flat metric on $\mathrm{H}_{\mathcal{P}}(\mathbb{R}^{d})$, see
Theorem \ref{t.4.8}. As a consequence, we found that under the development
map $\phi _{\mathcal{P}}$, the volume form $\nu _{\mathcal{P}}^{1}$ on $%
\mathrm{H}_{\mathcal{P}}(M)$ pulls back to the Gaussian density $\mu _{%
\mathcal{P}}^{1}$ on $\mathrm{H}_{\mathcal{P}}(\mathbb{R}^{d})$.

\begin{definition}
\label{d.5.1}Let $M^{\mathcal{P}}:=M^{n}$ and $\pi _{\mathcal{P}}:\mathrm{W}%
(M)\rightarrow M^{\mathcal{P}}$ denote the projection 
\begin{equation}
\pi _{\mathcal{P}}(\sigma ):=(\sigma (s_{1}),\ldots ,\sigma (s_{n})).
\label{e.5.1}
\end{equation}
We will also use the same notation for the restriction of $\pi _{\mathcal{P}%
} $ to $\mathrm{H}(M)$ and $\mathrm{H}_{\mathcal{P}}(M).$\qed
\end{definition}

In this section we will consider two further models for the geometry on path
space, namely the degenerate $L^{2}$-``metric'' $G_{\mathcal{P}}^{0}$
defined in Definition \ref{d.1.5} on $\mathrm{H}_{\mathcal{P}}(M)$ and the
product manifold $M^{\mathcal{P}}$ with its ``natural'' metric.

\begin{rem}
\label{r.5.2}The form $G_{\mathcal{P}}^{0}$ is non-negative but fails to be
definite precisely at $\sigma \in \mathrm{H}_{\mathcal{P}}(M)$ for which $%
\sigma (s_{i})$ is conjugate to $\sigma (s_{i-1})$ along $\sigma
([s_{i-1},s_{i}])$ for some $i.$ In this case there exists a nonzero $X\in T%
\mathrm{H}_{\mathcal{P}}(M)$ for which $G_{\mathcal{P}}^{0}(X,X)=0$. Hence, $%
\mathrm{Vol}_{G_{\mathcal{P}}^{0}}$ will also be zero for such $\sigma \in 
\mathrm{H}_{\mathcal{P}}(M)$. \qed
\end{rem}

\begin{definition}
\label{d.5.3} Let $M^{\mathcal{P}}$ be as in Definition \ref{d.5.1}. For $%
\mathbf{x}=(x_{1},x_{2},\ldots ,x_{n})\in M^{\mathcal{P}}$, let 
\begin{equation}
E_{\mathcal{P}}(\mathbf{x}):=\sum_{i=1}^{n}\frac{d^{2}(x_{i-1},x_{i})}{%
\Delta _{i}s}  \label{e.5.2}
\end{equation}
where $d$ is the geodesic distance on $M$. Let $g_{\mathcal{P}}$ be the
Riemannian metric on $M^{\mathcal{P}}$ given by 
\begin{equation}
g_{\mathcal{P}}=\left( \Delta _{1}s\right) g\times \left( \Delta
_{2}s\right) g\times \cdots \times \left( \Delta _{n}s\right) g,
\label{e.5.3}
\end{equation}
i.e. if $\mathbf{v}=(v_{1},v_{2},\dots ,v_{n})\in TM^{n}=(TM)^{n}$ then 
\begin{equation*}
g_{\mathcal{P}}(\mathbf{v},\mathbf{v}):=\sum_{i=1}^{n}g(v_{i},v_{i})\Delta
_{i}s.
\end{equation*}
Let the normalizing constant $Z_{\mathcal{P}}^{0}$ be given by equation \ref
{e.1.15} and let $\gamma _{\mathcal{P}}$ denote the measure on $M^{\mathcal{P%
}}$ defined by 
\begin{equation}
\gamma _{\mathcal{P}}(d\mathbf{x}):=\frac{1}{Z_{\mathcal{P}}^{0}}\exp \left(
-\frac{1}{2}E_{\mathcal{P}}(\mathbf{x})\right) \mathrm{Vol}_{g_{\mathcal{P}%
}}(d\mathbf{x})  \label{e.5.4}
\end{equation}
where $\mathrm{Vol}_{g_{\mathcal{P}}}$ denotes volume form on $M^{\mathcal{P}%
}$ defined with respect to. $g_{\mathcal{P}}$. \qed
\end{definition}

\begin{rem}
\label{r.5.4}An easy computation shows that 
\begin{equation}
\mathrm{Vol}_{g_{\mathcal{P}}}=\left( \prod_{i=1}^{n}(\Delta
_{i}s)^{d/2}\right) \times \mathrm{Vol}_{g}^{n},  \label{e.5.5}
\end{equation}
where $\mathrm{Vol}_{g}$ is the volume measure on $(M,g)$ and $\mathrm{Vol}%
_{g}^{n}$ denotes the $n$--fold product of $\mathrm{Vol}_{g}$with itself.%
\qed
\end{rem}

The next proposition shows the relationship between $\nu _{\mathcal{P}}^{0}$
(defined in Definition \ref{d.1.7} above) and $\gamma _{\mathcal{P}}$. For
the statement we need to define a subset of paths $\sigma $ in $\mathrm{H}_{%
\mathcal{P}}(M)$ such that each geodesic piece $\sigma ([s_{i-1},s_{i}])$ is
short. The formal definition is as follows.

\begin{definition}
\label{d.5.5}\ 

\begin{enumerate}
\item  For any $\epsilon >0$, let 
\begin{equation*}
\mathrm{H}_{\mathcal{P}}^{\epsilon }(M):=\{\sigma \in \mathrm{H}_{\mathcal{P}%
}(M):\int_{s_{i-1}}^{s_{i}}|\sigma ^{\prime }(s)|ds<\epsilon \quad \text{%
\textrm{for }}i=1,2,\ldots ,n\}.
\end{equation*}

\item  For any $\epsilon >0$, let 
\begin{equation*}
M_{\epsilon }^{\mathcal{P}}=\{\mathbf{x}\in M^{\mathcal{P}%
}:d(x_{i-1},x_{i})<\epsilon \quad \text{\textrm{for }}i=1,2,\ldots ,n\}
\end{equation*}
where $d$ is the geodesic distance on $(M,g)$ and $x_{0}:=o.$\qed
\end{enumerate}
\end{definition}

\begin{prop}
\label{p.5.6} For $\epsilon >0$ less than the injectivity radius of $M$, we
have

\begin{enumerate}
\item  \label{point.1} $G_{\mathcal{P}}^{0}$ is a Riemannian metric on $%
\mathrm{H}_{\mathcal{P}}^{\epsilon }(M)$.

\item  \label{point.2} The image of $\mathrm{H}_{\mathcal{P}}^{\epsilon }(M)$
under $\pi _{\mathcal{P}}$ is $M_{\epsilon }^{\mathcal{P}}$ and the map 
\begin{equation*}
\pi _{\mathcal{P}}:(\mathrm{H}_{\mathcal{P}}^{\epsilon }(M),G_{\mathcal{P}%
}^{0})\rightarrow (M_{\epsilon }^{\mathcal{P}},g_{\mathcal{P}})
\end{equation*}
is an isometry, where $g_{\mathcal{P}}$ is the metric on $M^{\mathcal{P}}$
in equation (\ref{e.5.3}).

\item  \label{point.3} $\pi _{\mathcal{P}}^{\ast }\gamma _{\mathcal{P}}=\nu
_{\mathcal{P}}^{0}$ on $\mathrm{H}_{\mathcal{P}}^{\epsilon }(M)$.
\end{enumerate}
\end{prop}

\begin{proof}
Because $\epsilon $ is less than the injectivity radius of $M$, it follows
that any $X\in T_{\sigma }\mathrm{H}_{\mathcal{P}}(M)$ is determined by its
values on the partition points $\mathcal{P}$. Therefore, if $G_{\mathcal{P}%
}^{0}(X,X)=0$ for $X\in T_{\sigma }\mathrm{H}_{\mathcal{P}}^{\epsilon }(M),$
then $X:=0$. This proves the first item. The second item is a triviality.
The last item is proved by noting that for $\sigma \in \mathrm{H}_{\mathcal{P%
}}^{\epsilon }(M)$, $\sigma |_{[s_{i-1},s_{i}]}$ is a minimal length
geodesic joining $\sigma (s_{i-1})$ to $\sigma (s_{i}),$ and therefore 
\begin{equation}
\int_{s_{i-1}}^{s_{i}}|\sigma ^{\prime }(s)|^{2}ds=\left( \frac{d(\sigma
(s_{i-1}),\sigma (s_{i}))}{\Delta _{i}s}\right) ^{2}\Delta _{i}s=\frac{%
d^{2}(\sigma (s_{i-1}),\sigma (s_{i}))}{\Delta _{i}s}.  \label{e.5.6}
\end{equation}
Summing this last equation on $i$ shows, 
\begin{equation}
E(\sigma )=\int_{0}^{1}|\sigma ^{\prime }(s)|^{2}ds=\sum_{i=1}^{n}\frac{%
d^{2}(\sigma (s_{i-1}),\sigma (s_{i}))}{\Delta _{i}s}=E_{\mathcal{P}}(\pi _{%
\mathcal{P}}(\sigma )).  \label{e.5.7}
\end{equation}
Hence by the definition of $\gamma _{\mathcal{P}}$, the fact that $\pi _{%
\mathcal{P}}$ is an isometry on $\mathrm{H}_{\mathcal{P}}^{\epsilon }(M)$
(point \ref{point.2} above), and (\ref{e.5.7}) above, we find that on $%
\mathrm{H}_{\mathcal{P}}^{\epsilon }(M)$, 
\begin{equation*}
\pi _{\mathcal{P}}^{\ast }\gamma _{\mathcal{P}}=\frac{1}{Z_{\mathcal{P}}^{0}}%
e^{-E/2}\mathrm{Vol}_{G_{\mathcal{P}}^{0}}=\nu _{\mathcal{P}}^{0}.
\end{equation*}
\end{proof}

Note that in general, for $\mathbf{x}\in M^{\mathcal{P}}$, $\pi _{\mathcal{P}%
}^{-1}(\mathbf{x})$ has more than one element, and may even fail to be a
discrete subset. Therefore using the product manifold $M^{\mathcal{P}}$ as a
model for $\mathrm{H}_{\mathcal{P}}(M)$ requires some care. The important
aspect of the isometric subsets $M_{\epsilon }^{\mathcal{P}}$ and $\mathrm{H}%
_{\mathcal{P}}^{\epsilon }(M)$ is that in a precise sense they have nearly
full measure with respect to. $\gamma _{\mathcal{P}},\nu _{\mathcal{P}}^{1}$
and $\nu _{\mathcal{P}}^{0}$. This will be proved in section \ref{s.5.1}
below.

Before carrying out these estimates we will finish this section by comparing 
$\nu _{\mathcal{P}}^{0}$ to $\nu _{\mathcal{P}}^{1}$.

\begin{notation}
\label{n.5.7} Let $\mathbb{R}^{d\mathcal{P}}$ denote the Euclidean space $(%
\mathbb{R}^{d})^{n}$ equipped with the product inner product defined in the
same way as $g_{\mathcal{P}}$ in equation (\ref{e.5.3}) with $\mathbb{R}^{d} 
$ replacing $TM.$\qed
\end{notation}

To simplify notation throughout this section, let 
\begin{equation}
\sigma \in \mathrm{H}_{\mathcal{P}}(M),\quad b:=\phi ^{-1}(\sigma ),\quad
u:=//(\sigma ),\quad \text{and }A(s):=\Omega _{u(s)}(b^{\prime }(s),\cdot
)b^{\prime }(s).  \label{e.5.8}
\end{equation}
Note that since $b\in \mathrm{H}_{\mathcal{P}}(\mathbb{R}^{d})$, 
\begin{equation}
b^{\prime }(s)=\Delta _{i}b/\Delta _{i}s\quad \text{\textrm{and}}\quad
A(s)=\Omega _{u(s)}(\frac{\Delta _{i}b}{\Delta _{i}s},\cdot )\frac{\Delta
_{i}b}{\Delta _{i}s}  \label{e.5.9}
\end{equation}
for $s\in (s_{i-1},s_{i}]$. Let us also identify $X\in T_{\sigma }\mathrm{H}%
_{\mathcal{P}}(M)$ with $h:=u^{-1}X$. Recall from Proposition \ref{p.4.4}
that $h:[0,1]\rightarrow \mathbb{R}^{d}$ is a piecewise smooth function such
that $h(0)=0$ and Equation (\ref{e.4.2}) holds, i.e. 
\begin{equation}
h^{\prime \prime }=Ah\text{ on }I\setminus \mathcal{P}\text{ and }h(0)=0\in 
\mathbb{R}^{d}.  \label{e.5.10}
\end{equation}
In order to compare $\mathrm{Vol}_{G_{\mathcal{P}}^{0}}$ and $\mathrm{Vol}%
_{G_{\mathcal{P}}^{1}}$ it is useful to define two linear maps 
\begin{eqnarray*}
J_{0} &:&(T_{\sigma }\mathrm{H}_{\mathcal{P}}(M),G_{\mathcal{P}%
}^{0})\rightarrow \mathbb{R}^{d\mathcal{P}} \\
J_{1} &:&(T_{\sigma }\mathrm{H}_{\mathcal{P}}(M),G_{\mathcal{P}%
}^{1})\rightarrow \mathbb{R}^{d\mathcal{P}}
\end{eqnarray*}
by 
\begin{equation*}
J_{0}(X)=(h(s_{1}),h(s_{2}),\ldots ,h(s_{n}))
\end{equation*}
and 
\begin{equation*}
J_{1}(X)=(h^{\prime }(s_{0}+),h^{\prime }(s_{1}+),\ldots ,h^{\prime
}(s_{n-1}+))
\end{equation*}
where $h:=u^{-1}X$ as above.

It follows from the definition of $G_{\mathcal{P}}^{0}$ and the metric on $%
\mathbb{R}^{d\mathcal{P}}$ that if $\sigma $ is such that $J_{0}$ is
injective, then $J_{0}$ is an isometry. By point \ref{point.2} of
Proposition \ref{p.5.6} this holds on $\mathrm{H}_{\mathcal{P}}^{\epsilon
}(M)$. However, by Remark \ref{r.5.2}{\ there is in general a nonempty
subset of $\mathrm{H}_{\mathcal{P}}(M)$ where $J_{0}$ fails to be injective.
Clearly, $J_{0}$ fails to be injective precisely where $G_{\mathcal{P}}^{0}$
fails to be positive definite. Similarly, it is immediate from the
definitions and the fact that $G_{\mathcal{P}}^{1}$ is a nondegenerate
Riemann metric, see Remark \ref{r.4.5}, that $J_{1}$ is an isometry at all $%
\sigma \in \mathrm{H}_{\mathcal{P}}(M)$. }

To simplify notation, let $V$ denote the vector space $(\mathbb{R}^{d})^{n}$
and let $T=T_{\mathcal{P}}(\sigma )$ be defined by $T:=J_{0}\circ
J_{1}^{-1}. $ Thus $T:V\rightarrow V$ is the unique linear map such that 
\begin{equation}
T(h^{\prime }(s_{0}+),h^{\prime }(s_{1}+),\ldots ,h^{\prime
}(s_{n-1}+))=(h(s_{1}),h(s_{2}),\ldots ,h(s_{n}))  \label{e.5.11}
\end{equation}
for all $h=u^{-1}X$ with $X\in T_{\sigma }\mathrm{H}_{\mathcal{P}}(M)$. With
this notation it follows that 
\begin{eqnarray}
\mathrm{Vol}_{G_{\mathcal{P}}^{0}} &=&J_{0}^{\ast }\mathrm{Vol}_{\mathbb{R}%
^{d\mathcal{P}}} =(T\circ J_{1})^{\ast }\mathrm{Vol}_{\mathbb{R}^{d\mathcal{P%
}}}  \notag \\
&=&J_{1}^{\ast }T^{\ast }\mathrm{Vol}_{\mathbb{R}^{d\mathcal{P}}}=\det
(T)J_{1}^{\ast }\mathrm{Vol}_{\mathbb{R}^{d\mathcal{P}}}  \notag \\
&=&\det (T)\mathrm{Vol}_{G_{\mathcal{P}}^{1}}.  \label{e.5.12}
\end{eqnarray}
Note that in this computation $\sigma \in \mathrm{H}_{\mathcal{P}}(M)$ is
fixed and we treat $\mathrm{Vol}_{G_{\mathcal{P}}^{0}}, \mathrm{Vol}_{G_{%
\mathcal{P}}^{1}}$ as elements of the exterior algebra $\wedge ^{dn}\left(
T_{\sigma }^{\ast }\mathrm{H}_{\mathcal{P}}(M)\right) $ at some fixed $%
\sigma $ and $\mathrm{Vol}_{\mathbb{R}^{d\mathcal{P}}}$ as an element of $%
\wedge ^{dn}\left( (\mathbb{R}^{d\mathcal{P}})^{\ast }\right) $.

Our next task is to compute $\det(T)$.

\begin{lemma}
\label{l.5.8} Let $Z_{i-1}(s)$ denote the $d\times d$ matrix--valued
solution to 
\begin{equation}
Z_{i-1}^{\prime \prime }(s)=A(s)Z_{i-1}(s)\quad \text{\textrm{with }}\quad
Z_{i-1}(s_{i-1})=0\quad \text{\textrm{and }}\quad Z_{i-1}^{\prime
}(s_{i-1})=I.  \label{e.5.13}
\end{equation}
Then 
\begin{equation*}
\det (T_{\mathcal{P}}(\sigma ))=\prod_{i=1}^{n}\det (Z_{i-1}(s_{i})).
\end{equation*}
\end{lemma}

\begin{proof}
We start by noting that for $\sigma \in \mathrm{H}_{\mathcal{P}}(M)$ such
that $G_{\mathcal{P}}^{0}$ is nondegenerate, then $\det (Z_{i-1})\neq 0$ for 
$i=1,2,\dots ,n$. To see this assume that $\det (Z_{i-1})=0$ for some $i$.
In view of the fact that $Z$ solves the Jacobi equation (\ref{e.5.13}), this
is equivalent to the existence of a vector field $X_{i-1}$ along $\sigma
([s_{i-1},s_{i}])$ which solves (\ref{e.4.1}) for $s\in \lbrack
s_{i-1},s_{i}]$ and which satisfies 
\begin{equation*}
X_{i-1}(s_{i-1})=0,\qquad X_{i-1}(s_{i})=0.
\end{equation*}
Define $X$ by 
\begin{equation*}
X(s)=\left\{ 
\begin{array}{ll}
X_{i-1}(s), & s\in \lbrack s_{i-1},s_{i}] \\ 
0 & s\in \lbrack 0,1]\setminus \lbrack s_{i-1},s_{i}]
\end{array}
\right.
\end{equation*}

Then $X\in T_{\sigma }\mathrm{H}_{\mathcal{P}}(M)$ and it is clear from the
construction that $G_{\mathcal{P}}^{0}(X,X)=0$. Thus for such $\sigma $, $%
\mathrm{Vol}_{G_{\mathcal{P}}^{0}}|_{\sigma }=0$. Hence we may without loss
of generality restrict our considerations to the case when $\det
(Z_{i-1})\neq 0 $ for all $i$.

Let $C_{i-1}(s)$ be the $d\times d$ matrix--valued solutions to 
\begin{equation*}
C_{i-1}^{\prime \prime }(s)=A(s)C_{i-1}(s)\quad \text{\textrm{with }}\quad
C_{i-1}(s_{i-1})=I\quad \text{\textrm{and }}\quad C_{i-1}^{\prime
}(s_{i-1})=0.
\end{equation*}
For $i\in \{1,2,\ldots ,n\}$ and $h=u^{-1}X$ with $X\in T_{\sigma }\mathrm{H}%
_{\mathcal{P}}(M)$ let 
\begin{equation*}
k(s):=C_{i-1}(s)h(s_{i-1})+Z_{i-1}(s)h^{\prime }(s_{i-1}+).
\end{equation*}
Then $k^{\prime \prime }=Ak$ on $(s_{i-1},s_{i})$, $k(s_{i-1})=h(s_{i-1})$
and $k^{\prime }(s_{i-1})=h^{\prime }(s_{i-1}+)$. Since $h$ satisfies the
same linear differential equation with initial conditions at $s_{i-1}$, it
follows that $h=k$ on $[s_{i-1},s_{i}]$ and in particular that 
\begin{equation*}
h(s_{i})=C_{i-1}(s_{i})h(s_{i-1})+Z_{i-1}(s_{i})h^{\prime }(s_{i-1}+).
\end{equation*}
Solving this equation for $h^{\prime }(s_{i-1}+)$ gives 
\begin{equation*}
h^{\prime }(s_{i-1}+)=Z_{i-1}(s_{i})^{-1}(h(s_{i})-C_{i-1}(s_{i})h(s_{i-1}))
\end{equation*}
from which it follows that $T^{-1}(\xi _{1},\xi _{2},\ldots ,\xi _{n})=(\eta
_{1},\eta _{2},\ldots ,\eta _{n})$ where 
\begin{equation*}
\eta _{i}=\alpha _{i}\xi _{i}-\beta _{i}\xi _{i-1}\quad \text{\textrm{for }}%
\quad i=1,2,\ldots ,n,
\end{equation*}
$\alpha _{i}:=Z_{i-1}(s_{i})^{-1}$ and $\beta
_{i}:=Z_{i-1}(s_{i})^{-1}C_{i-1}(s_{i})$. (In the previous displayed
equation $\xi _{0}$ should be interpreted as $0$.) Thus the linear
transformation $T^{-1}:V\rightarrow V$ may be written in block lower
triangular form as 
\begin{equation*}
T^{_{-1}}=\left[ 
\begin{array}{ccccc}
\alpha _{1} & 0 & 0 & \cdots & 0 \\ 
\beta _{2} & \alpha _{2} & 0 & \cdots & 0 \\ 
0 & \beta _{3} & \alpha _{3} & \ddots & 0 \\ 
\vdots & \vdots & \ddots & \ddots & 0 \\ 
0 & 0 & \cdots & \beta _{n} & \alpha _{n}
\end{array}
\right]
\end{equation*}
and hence for $\sigma \in \mathrm{H}_{\mathcal{P}}(M)$ so that $G_{\mathcal{P%
}}^{0}$ is nondegenerate, 
\begin{equation*}
\det (T^{-1})=\prod_{i=1}^{n}\det (\alpha _{i})=\prod_{i=1}^{n}\det
(Z_{i-1}(s_{i})^{-1}).
\end{equation*}
It follows by the above arguments, that for all $\sigma \in \mathrm{H}_{%
\mathcal{P}}(M)$ 
\begin{equation*}
\det (T)=\prod_{i=1}^{n}\det (Z_{i-1}(s_{i}))
\end{equation*}
\end{proof}

As a consequence, we have the key theorem relating $\nu _{\mathcal{P}}^{0}$
to $\nu _{\mathcal{P}}^{1}$.

\begin{thm}
\label{t.5.9}Let 
\begin{equation}
\rho _{\mathcal{P}}(\sigma ):=\prod_{i=1}^{n}\det (\frac{Z_{i-1}(s_{i})}{%
\Delta _{i-1}s}),  \label{e.5.14}
\end{equation}
then $\nu _{\mathcal{P}}^{0}=\rho _{\mathcal{P}}\nu _{\mathcal{P}}^{1}$.
\end{thm}

\begin{proof}
>From the Definition \ref{d.1.6} for $\nu _{\mathcal{P}}^{0}$, Equation (\ref
{e.5.12}) and Lemma \ref{e.5.10} we find that 
\begin{eqnarray*}
\nu _{\mathcal{P}}^{0} &=&\frac{1}{Z_{\mathcal{P}}^{0}}e^{-\frac{1}{2}E}%
\mathrm{Vol}_{G_{\mathcal{P}}^{0}} \\
&=&\frac{1}{Z_{\mathcal{P}}^{0}}e^{-\frac{1}{2}E}\prod_{i=1}^{n}\det
(Z_{i-1}(s_{i}))\mathrm{Vol}_{G_{\mathcal{P}}^{1}} \\
&=&\frac{1}{Z_{\mathcal{P}}^{0}}e^{-\frac{1}{2}E}\prod_{i=1}^{n}(\Delta
_{i-1}s)^{d}\cdot \prod_{i=1}^{n}\det (\frac{1}{\Delta _{i-1}s}%
Z_{i-1}(s_{i}))\mathrm{Vol}_{G_{\mathcal{P}}^{1}}.
\end{eqnarray*}
Equation (\ref{e.5.14}) now follows from Definition \ref{d.1.7} (for $\nu _{%
\mathcal{P}}^{1}$) and the expressions for $Z_{\mathcal{P}}^{1}$ and $Z_{%
\mathcal{P}}^{0}$ in equation (\ref{e.1.15}).

Using this result and Bishop's Comparison Theorem we have the following
estimate on $\rho _{\mathcal{P}}(\sigma ).$
\end{proof}

\begin{cor}
\label{c.5.10}Let $K>0$ be such that $\mathrm{Ric}\geq -(d-1)KI$ (for
example take $K $ to be a bound on $\Omega )$ then 
\begin{equation}
\rho _{\mathcal{P}}(\sigma )\leq \prod_{i=1}^{n}\left( \frac{\sinh (\sqrt{K}%
|\Delta _{i}b|)}{\sqrt{K}|\Delta _{i}b|}\right) ^{d-1}.  \label{e.5.15}
\end{equation}
\end{cor}

\begin{proof}
The proof amounts to applying Theorem 3.8 on p. 120 \cite{Chavel93} to each
of the $Z_{i-1}(s_{i})$'s above. In order to use this theorem one must keep
in mind that $\frac{\Delta _{i}b}{\Delta _{i}s}$ is not a unit vector and
the estimate given in \cite{Chavel93} corresponds to the determinant of $%
Z_{i-1}(s_{i})$ restricted $\left\{ \xi :=\frac{\Delta _{i}b}{\Delta _{i}s}%
\right\} ^{\perp }.$ Noting that $Z_{i-1}(s_{i})\xi =\Delta _{i}s\cdot \xi $
and accounting for the aforementioned discrepancies, Theorem 3.8 in \cite
{Chavel93} gives the estimate 
\begin{equation*}
\det \left( Z_{i-1}(s_{i})\right) \leq \left( \frac{\sinh (\sqrt{K}|\Delta
_{i}b|)}{\sqrt{K}|\Delta _{i}b|/\Delta _{i}s}\right) ^{d-1}\Delta _{i}s
\end{equation*}
or equivalently that 
\begin{equation*}
\det (\frac{1}{\Delta _{i-1}s}Z_{i-1}(s_{i}))\leq \left( \frac{\sinh (\sqrt{K%
}|\Delta _{i}b|)}{\sqrt{K}|\Delta _{i}b|}\right) ^{d-1}.
\end{equation*}
This clearly implies the estimate in equation (\ref{e.5.15}).
\end{proof}

\subsection{Estimates of the measure of $\mathrm{H}_{\mathcal{P}}^{\protect%
\epsilon }(M)$ and $M_{\protect\epsilon }^{\mathcal{P}}$\label{s.5.1}}

We will need the following Lemma, which is again a consequence of Bishop's
comparison theorem.

\begin{lemma}
\label{l.5.11}Let $\omega _{d-1}$ denote the surface area of the unit sphere
in $\mathbb{R}^{d},$ $R$ be the diameter of $M$ and let $K\geq 0$ such that $%
\mathrm{Ric}\geq -(d-1)KI.$ Then for all $F:[0,R]\rightarrow \lbrack
0,\infty ]$, 
\begin{equation*}
\int_{M}F(d(o,\cdot ))dvol\leq \omega _{d-1}\int_{0}^{R}r^{d-1}F(r)\left( 
\frac{\sinh (\sqrt{K}r)}{\sqrt{K}r}\right) ^{d-1}dr
\end{equation*}
\end{lemma}

\begin{proof}
See Equations (2.48) on p. 72 (3.15) on p. 113, and Theorem 3.8 on p. 120 in
Chavel \cite{Chavel93}.
\end{proof}

We are now ready to estimate the measures of $M_{\epsilon }^{\mathcal{P}}$
and $\mathrm{H}_{\mathcal{P}}^{\epsilon }(M)$. We start by considering $%
\gamma _{\mathcal{P}}(M^{\mathcal{P}}\setminus M_{\epsilon }^{\mathcal{P}})$.

\begin{prop}
\label{p.5.12} Fix $\epsilon >0$ and let $M_{\epsilon }^{\mathcal{P}}$ be as
in Definition \ref{d.5.5} and let $\gamma _{\mathcal{P}}$ be the measure on $%
M^{\mathcal{P}}$ defined by (\ref{e.5.4}). Then there is a constants $%
C<\infty $ such that 
\begin{equation*}
\gamma _{\mathcal{P}}(M^{\mathcal{P}}\setminus M_{\epsilon }^{\mathcal{P}%
})\leq C\exp (-\frac{\epsilon ^{2}}{4|\mathcal{P}|}).
\end{equation*}
\end{prop}

\begin{proof}
Let $f:[0,\infty )^{n}\rightarrow \lbrack 0,\infty )$ be a measurable
function. Let $d\mathbf{x}=\prod_{i=1}^{n}\mathrm{Vol}_{g}(dx_{i})$ and note
that 
\begin{equation}
d\mathbf{x}=d\mathrm{Vol}_{g_{\mathcal{P}}}(\mathbf{x})\prod_{i=1}^{n}(%
\Delta _{i}s)^{-d/2}.  \label{e.5.16}
\end{equation}
An application of Lemma \ref{l.5.11} and Fubini's theorem proves 
\begin{multline*}
\int_{M^{\mathcal{P}}}f(d(o,x_{1}),d(x_{1},x_{2}),\ldots
,d(x_{n-1},x_{n}))\gamma _{\mathcal{P}}(d\mathbf{x}) \\
\leq \int_{\lbrack 0,\infty )^{n}}f(r_{1},r_{2},\ldots ,r_{n})\exp \left(
-\sum_{i=1}^{n}\frac{r_{i}^{2}}{2\Delta _{i}s}\right) \prod_{i=1}^{n}\left( 
\frac{\sinh (\sqrt{K}r_{i})}{\sqrt{K}r_{i}}\right) ^{d-1}\frac{\omega
_{d-1}r_{i}^{d-1}dr_{i}}{(2\pi \Delta _{i-1}s)^{d/2}}.
\end{multline*}
As usual let $\{B(s)\}_{s\in \lbrack 0,1]}$ be a standard $\mathbb{R}^{d}$%
--valued Brownian motion in Notation \ref{n.4.12} and $\Delta
_{i}B=B(s_{i})-B(s_{i-1})$. Noting that 
\begin{equation*}
\exp \left( -\sum_{i=1}^{n}\frac{r_{i}^{2}}{2\Delta _{i}s}\right)
\prod_{i=1}^{n}\frac{\omega _{d-1}r_{i}^{d-1}dr_{i}}{(2\pi \Delta
_{i-1}s)^{d/2}}
\end{equation*}
is the distribution of $(|\Delta _{1}B|,|\Delta _{2}B|,\ldots ,|\Delta
_{n}B|),$ the above inequality may be written as: 
\begin{align}
& \int_{M^{\mathcal{P}}}f(d(o,x_{1}),d(x_{1},x_{2}),\ldots
,d(x_{n-1},x_{n}))\gamma _{\mathcal{P}}(d\mathbf{x})  \label{e.5.17} \\
& \leq \mathbb{E}\left[ f(|\Delta _{1}B|,|\Delta _{2}B|,\ldots ,|\Delta
_{n}B|)\prod_{i=1}^{n}\left( \frac{\sinh (\sqrt{K}|\Delta _{i}B|)}{\sqrt{K}%
|\Delta _{i}B|}\right) ^{d-1}\right] .  \notag
\end{align}
For $i\in \{1,2,\ldots ,n\},$ let $\mathcal{A}_{i}:=\{\mathbf{x}\in M^{%
\mathcal{P}}:d(x_{i-1},x_{i})\geq \epsilon \}$ so that $M^{\mathcal{P}%
}\setminus M_{\epsilon }^{\mathcal{P}}=\cup _{i=1}^{n}\mathcal{A}_{i}$ and 
\begin{equation}
\gamma _{\mathcal{P}}(M^{\mathcal{P}}\setminus M_{\epsilon }^{\mathcal{P}%
})\leq \sum_{i=1}^{n}\gamma _{\mathcal{P}}(\mathcal{A}_{i}).  \label{e.5.18}
\end{equation}
Since $1_{\mathcal{A}_{i}}(\mathbf{x})=\chi _{\epsilon }(d(x_{i-1},x_{i})),$
where $\chi _{\epsilon }(r)=1_{r\geq \epsilon }$, we find from equation (\ref
{e.5.17}) that 
\begin{align}
\gamma _{\mathcal{P}}(\mathcal{A}_{i})& \leq \mathbb{E}\left[ \chi
_{\epsilon }(|\Delta _{i}B|)\prod_{j=1}^{n}\left( \frac{\sinh (\sqrt{K}%
|\Delta _{j}B|)}{\sqrt{K}|\Delta _{j}B|}\right) ^{d-1}\right]  \notag \\
& =\mathbb{E}\left[ \chi _{\epsilon }(|\Delta _{i}B|)\left( \frac{\sinh (%
\sqrt{K}|\Delta _{i}B|)}{\sqrt{K}|\Delta _{i}B|}\right) ^{d-1}\right]
\prod_{j\neq i}\psi (\sqrt{\Delta _{j}s}),  \label{e.5.19}
\end{align}
where $\psi $ is defined in equation (\ref{e.8.19}) of the Appendix. An
application of Lemma \ref{l.8.7} now completes the proof in view of (\ref
{e.5.18}) and (\ref{e.5.19}).
\end{proof}

We also have the following analogue of Proposition \ref{p.5.12}.

\begin{prop}
\label{p.5.13} For any $\epsilon >0$ there is a constant $C<\infty $ such
that 
\begin{equation*}
\nu _{\mathcal{P}}^{1}(\mathrm{H}_{\mathcal{P}}(M)\setminus \mathrm{H}_{%
\mathcal{P}}^{\epsilon }(M))\leq C\exp (-\frac{\epsilon ^{2}}{4|\mathcal{P}|}%
).
\end{equation*}
\end{prop}

\begin{proof}
Let us recall that $\phi (\mathrm{H}_{\mathcal{P}}(\mathbb{R}^{d}))=\mathrm{H%
}_{\mathcal{P}}(M)$ and let us note that $\phi (\mathrm{H}_{\mathcal{P}%
}^{\epsilon }(\mathbb{R}^{d}))=\mathrm{H}_{\mathcal{P}}^{\epsilon }(M)$. By
Theorem \ref{t.4.10} and Corollary \ref{c.4.13} this implies that 
\begin{eqnarray*}
\nu _{\mathcal{P}}^{1}(\mathrm{H}_{\mathcal{P}}(M)\setminus \mathrm{H}_{%
\mathcal{P}}^{\epsilon }(M)) &=&\mu _{\mathcal{P}}^{1}(\mathrm{H}_{\mathcal{P%
}}(\mathbb{R}^{d})\setminus \mathrm{H}_{\mathcal{P}}^{\epsilon }(\mathbb{R}%
^{d})) \\
&=&\mu (\{\max \{|\Delta _{i-1}B|:i=1,2,\ldots ,n\}\geq \epsilon \}) \\
&\leq &\sum_{i=1}^{n}\mu (|\Delta _{i-1}B|\geq \epsilon ) \\
&=&\sum_{i=1}^{n}\mathbb{E}\left[ \chi _{\epsilon }(|\Delta _{i}B|)\right] \\
&\leq &Ce^{-\frac{\epsilon ^{2}}{4|\mathcal{P}|}}.
\end{eqnarray*}
where as above $\chi _{\epsilon }(r)=1_{r\geq \epsilon }.$ The last
inequality follows from Lemma \ref{l.8.7} with $K=0$.
\end{proof}

Finally we consider $\nu _{\mathcal{P}}^{0}(\mathrm{H}_{\mathcal{P}%
}(M)\setminus \mathrm{H}_{\mathcal{P}}^{\epsilon }(M))$.

\begin{prop}
\label{p.5.14} Let $\epsilon >0$. Then there is a constant $C<\infty $ such
that 
\begin{equation}
\nu _{\mathcal{P}}^{0}(\mathrm{H}_{\mathcal{P}}(M)\setminus \mathrm{H}_{%
\mathcal{P}}^{\epsilon }(M))\leq C\exp (-\frac{\epsilon ^{2}}{4|\mathcal{P}|}%
).  \label{e.5.20}
\end{equation}
\end{prop}

\begin{proof}
Let $B$ be the standard $\mathbb{R}^{d}$ valued Brownian motion. For $%
i=1,2,\dots ,n$, let $\mathcal{A}_{i}=\{|\Delta _{i}B|>\epsilon \}$ and set $%
\mathcal{A}=\cup _{i=1}^{n}\mathcal{A}_{i}$. Then $\mathrm{H}_{\mathcal{P}%
}(M)\setminus \mathrm{H}_{\mathcal{P}}^{\epsilon }(M)=\phi _{\mathcal{P}}(%
\mathcal{A})$ where $\phi _{\mathcal{P}}:\mathrm{H}_{\mathcal{P}}(\mathbb{R}%
^{d})\rightarrow \mathrm{H}_{\mathcal{P}}(M)$ denotes the development map.

By Theorem \ref{t.5.9}, $\nu _{\mathcal{P}}^{0}=\rho _{\mathcal{P}}\nu _{%
\mathcal{P}}^{1}$, where $\rho _{\mathcal{P}}$ is given by (\ref{e.5.14}).
By Theorem \ref{t.4.10} and Corollary \ref{c.5.10} above, 
\begin{equation*}
\nu _{\mathcal{P}}^{0}(\mathrm{H}_{\mathcal{P}}(M)\setminus \mathrm{H}_{%
\mathcal{P}}^{\epsilon }(M))=\int_{\mathcal{A}}\rho _{\mathcal{P}}(\phi (B_{%
\mathcal{P}}))d\mu _{\mathcal{P}}^{1}\leq \int_{\mathcal{A}%
}\prod_{i=1}^{n}\left( \frac{\sinh (\sqrt{K}|\Delta _{i}B|)}{\sqrt{K}|\Delta
_{i}B|}\right) ^{d-1}d\mu ,
\end{equation*}
wherein we have used the fact that the distribution of $\left\{ \Delta
_{i}B_{\mathcal{P}}\right\} _{i}$ under $\mu _{\mathcal{P}}^{1}$ is the same
as the distribution of $\left\{ \Delta _{i}B\right\} _{i}$ under $\mu .$
Thus arguing as in the proof of Proposition \ref{p.5.12} we have with $\chi
_{\epsilon }=1_{r\geq \epsilon }$, 
\begin{align*}
\nu _{\mathcal{P}}^{0}(\phi _{\mathcal{P}}(\mathcal{A}))& \leq
\sum_{i=1}^{n}\nu _{\mathcal{P}}^{0}(\phi _{\mathcal{P}}(\mathcal{A}_{i})) \\
& \leq \sum_{i=1}^{n}\mathbb{E}\left[ \chi _{\epsilon }(|\Delta
_{i}B|)\prod_{j=1}^{n}\left( \frac{\sinh (\sqrt{K}|\Delta _{j}B|)}{\sqrt{K}%
|\Delta _{j}B|}\right) ^{d-1}\right] \\
& =\sum_{i=1}^{n}\mathbb{E}\left[ \chi _{\epsilon }(|\Delta _{i}B|)\left( 
\frac{\sinh (\sqrt{K}|\Delta _{i}B|)}{\sqrt{K}|\Delta _{i}B|}\right) ^{d-1}%
\right] \prod_{j\neq i}\psi (\sqrt{\Delta _{j}s}).
\end{align*}
where $\psi $ is defined in equation (\ref{e.8.19}) of the Appendix. An
application of Lemma \ref{l.8.7} in the Appendix completes the proof.
\end{proof}

\section{Convergence of $\protect\nu _{\mathcal{P}}^{0}$ to Wiener Measure%
\label{s.6}}

This section is devoted to the proof of the following Theorem.

\begin{thm}
\label{t.6.1}Let $F:W(O(M))\rightarrow \mathbb{R}$ be a continuous and
bounded function and set $f(\sigma ):=F(//_{\cdot }(\sigma ))$ for $\sigma
\in \mathrm{H}(M).$ Then 
\begin{equation*}
\lim_{|\mathcal{P}|\rightarrow 0}\int_{\mathrm{H}_{\mathcal{P}}(M)}f(\sigma
)d\nu _{\mathcal{P}}^{0}(\sigma )=\int_{W(M)}\tilde{f}(\sigma )e^{-\frac{1}{6%
}\int_{0}^{1}\mathrm{Scal}(\sigma (s))ds}d\nu (\sigma ),
\end{equation*}
where $\tilde{f}(\sigma ):=F(/\tilde{/}_{\cdot }(\sigma ))$ and $/\tilde{/}%
_{r}(\sigma )$ is stochastic parallel translation, see Definition \ref
{d.4.15}. \qed
\end{thm}

Because of Theorem \ref{t.4.17}, in order to prove this theorem it will
suffice to compare $\nu _{\mathcal{P}}^{1}$ with $\nu _{\mathcal{P}}^{0}$.
Of course the main issue is to compare $\mathrm{Vol}_{G_{\mathcal{P}}^{0}}$
with $\mathrm{Vol}_{G_{\mathcal{P}}^{1}}$. In view of Proposition \ref
{p.5.14} and the boundedness of $f$ and $\mathrm{Scal}$, 
\begin{equation*}
\left| \int_{\mathrm{H}_{\mathcal{P}}(M)\setminus \mathrm{H}_{\mathcal{P}%
}^{\epsilon }(M)}f(\sigma )d\nu _{\mathcal{P}}^{0}(\sigma )\right| \leq
C||f||_{\infty }e^{-\frac{\epsilon ^{2}}{4|\mathcal{P}|}}
\end{equation*}
which tends to zero faster than any power of $|\mathcal{P}|$. Therefore, it
suffices to prove that for any $\epsilon >0$ smaller than the injectivity
radius of $M$, 
\begin{equation}
\lim_{|\mathcal{P}|\rightarrow 0}\int_{\mathrm{H}_{\mathcal{P}}^{\epsilon
}(M)}f(\sigma )d\nu _{\mathcal{P}}^{0}(\sigma )=\int_{W(M)}\tilde{f}(\sigma
)e^{-\frac{1}{6}\int_{0}^{1}\mathrm{Scal}(\sigma (s))ds}d\nu (\sigma ).
\label{e.6.1}
\end{equation}

\subsection{Estimating the Radon Nikodym Derivative\label{s.6.1}}

In this section we will continue to use the notation set out in equation (%
\ref{e.5.8}).

\begin{prop}
\label{p.6.2} Suppose that $A$ is given by equation (\ref{e.5.9})\ and that $%
Z_{i-1}$ is defined as in Lemma \ref{l.5.8}. Let $\Lambda $ be an upper
bound for both the norms of the curvature tensor $R$ (or equivalently $%
\Omega )$ and its covariant derivative $\nabla R.$ Then 
\begin{equation}
Z_{i-1}(s_{i})=\Delta _{i}s(I+\frac{1}{6}\Omega _{u(s_{i-1})}(\Delta
_{i}b,\cdot )\Delta _{i}b+\mathcal{E}_{i-1}),  \label{e.6.2}
\end{equation}
where 
\begin{equation}
|\mathcal{E}_{i-1}|\leq \frac{1}{6}(2\Lambda |\Delta _{i}b|^{3}+\frac{1}{2}%
\Lambda ^{2}|\Delta _{i}b|^{4})\cosh (\sqrt{\Lambda }|\Delta _{i}b|).
\label{e.6.3}
\end{equation}
In particular, if $\epsilon >0$ is given and it is assumed that $|\Delta
_{i}b|\leq \epsilon $ for all $i$, then 
\begin{equation}
|\mathcal{E}_{i-1}|\leq C|\Delta _{i}b|^{3},  \label{e.6.4}
\end{equation}
where $C=C(\epsilon ,R,\nabla R)=\frac{1}{6}(2\Lambda +\frac{1}{2}\Lambda
^{2}\epsilon )\cosh (\sqrt{\Lambda }\epsilon )$.
\end{prop}

\begin{proof}
By Lemma \ref{l.8.3} of the Appendix, 
\begin{equation}
Z_{i-1}(s_{i})=\Delta _{i}sI+\frac{\Delta _{i}s^{3}}{6}\Omega _{u(s_{i-1})}(%
\frac{\Delta _{i}b}{\Delta _{i}s},\cdot )\frac{\Delta _{i}b}{\Delta _{i}s}%
+\Delta _{i}s\mathcal{E}_{i-1},  \label{e.6.5}
\end{equation}
with $\mathcal{E}_{i-1}$ satisfying the estimate, 
\begin{equation}
|\mathcal{E}_{i-1}|=\frac{1}{6}(2K_{1}(\Delta _{i}s)^{3}+\frac{1}{2}%
K^{2}(\Delta _{i}s)^{4})\cosh (\sqrt{K}\Delta _{i}s)  \label{e.6.6}
\end{equation}
where $K:=\sup_{s\in (s_{i-1},s_{i})}|A(s)|$ and $K_{1}:=\sup_{s\in
(s_{i-1},s_{i})}|A^{\prime }(s)|$.

By (\ref{e.5.9}), for $s\in \lbrack s_{i-1},s_{i}]$, 
\begin{equation*}
|A(s)|\leq \Lambda |\Delta _{i}b|^{2}(\Delta _{i}s)^{-2}
\end{equation*}
and hence $K(\Delta _{i}s)^{2}\leq \Lambda |\Delta _{i}b|^{2}$.

Since $u^{\prime }(s)=\mathcal{H}_{u(s)}u(s)b^{\prime }(s),$ we see for $%
s_{i-1}<s\leq s_{i}$ that 
\begin{eqnarray*}
A^{\prime }(s) &=&(D\Omega )_{u(s)}(b^{\prime }(s),b^{\prime }(s),\cdot
)b^{\prime }(s) \\
&=&(\Delta _{i}s)^{-3}(D\Omega )_{u(s)}(\Delta _{i}b,\Delta _{i}b,\cdot
)\Delta _{i}b,
\end{eqnarray*}

where $(D\Omega )_{u(s)}(b^{\prime }(s),\cdot ,\cdot ):=\frac{d}{ds}\Omega
_{u(s)}.$ Therefore $|A^{\prime }(s)|\leq \Lambda (\Delta _{i}s)^{-3}|\Delta
_{i}b|^{3}$ which combined with equation (\ref{e.6.6}) proves equation (\ref
{e.6.3}).
\end{proof}

\begin{prop}
\label{p.6.3} Let $\Psi (U)$ be given as in Lemma \ref{l.8.1} of the
Appendix and define 
\begin{equation}
U_{i-1}:=\frac{1}{6}\Omega _{u(s_{i-1})}(\Delta _{i}b,\cdot )\Delta _{i}b+%
\mathcal{E}_{i-1},  \label{e.6.7}
\end{equation}
where $\mathcal{E}_{i-1}$ is defined in Proposition \ref{p.6.2}. Then 
\begin{equation}
\rho _{\mathcal{P}}(\sigma )=\exp (W_{\mathcal{P}}(\sigma ))\exp (-\frac{1}{6%
}\mathcal{R}_{\mathcal{P}}(\sigma )),  \label{e.6.8}
\end{equation}
where 
\begin{equation*}
\mathcal{R}_{\mathcal{P}}(\sigma ):=\sum_{i=1}^{n}\langle \mathrm{Ric}%
_{u(s_{i-1})}\Delta _{i}b,\Delta _{i}b\rangle
\end{equation*}
and 
\begin{equation}
W_{\mathcal{P}}(\sigma ):=\sum_{i=1}^{n}({\text{tr}}{\mathcal{E}}_{i-1}+\Psi
(-U_{i-1})).  \label{e.6.9}
\end{equation}
Moreover there exists $\epsilon _{0}>0$ and $C_{1}<\infty $ such that for
all $\epsilon \in (0,\epsilon _{0}],$ 
\begin{equation}
|W_{\mathcal{P}}(\sigma )|\leq C_{1}\sum_{i=1}^{n}|\Delta _{i}b|^{3}\text{%
\quad for all \quad }\sigma \in H_{\mathcal{P}}^{\epsilon }(M).
\label{e.6.10}
\end{equation}
\end{prop}

\begin{proof}
Recall that by definition, the trace of the linear map $v\mapsto \Omega
_{u(s_{i-1})}(\Delta _{i}b,v)\Delta _{i}b$ equals $-\langle \mathrm{Ric}%
_{u(s_{i-1})}\Delta _{i}b,\Delta _{i}b\rangle $ and hence 
\begin{equation*}
{\text{tr}} U_{i-1}=-\frac{1}{6}\langle \mathrm{Ric}_{u(s_{i-1})}\Delta
_{i}b,\Delta _{i}b\rangle +{\text{tr}}\mathcal{E}_{i-1}.
\end{equation*}
>From the definitions of $\mathcal{R}_{\mathcal{P}}$ and $W_{\mathcal{P}}$,
we get using Lemma \ref{l.5.8} and Lemma \ref{l.8.1}, 
\begin{eqnarray*}
\rho _{\mathcal{P}}(\sigma ) &=&\prod_{i=1}^{n}\exp \left( -\frac{1}{6}%
\langle \mathrm{Ric}_{u(s_{i-1})}\Delta _{i}b,\Delta _{i}b\rangle +{\text{tr}%
}\mathcal{E}_{i-1}+\Psi (-U_{i-1})\right) \\
&=&\exp \left( -\frac{1}{6}\sum_{i=1}^{n}\langle \mathrm{Ric}%
_{u(s_{i-1})}\Delta _{i}b,\Delta _{i}b\rangle \right) \exp \left(
\sum_{i=1}^{n}{\text{tr}}\mathcal{E}_{i-1}+\sum_{i=1}^{n}\Psi
(-U_{i-1})\right)
\end{eqnarray*}
which proves equation (\ref{e.6.8}).

Letting $\Lambda $ be a bound on the curvature tensor $\Omega $, it follows
using equation (\ref{e.6.4}) that 
\begin{eqnarray*}
|U_{i-1}| &\leq &\frac{1}{6}|\Omega _{u(s_{i-1})}(\Delta _{i}b,\cdot )\Delta
_{i}b|+|\mathcal{E}_{i-1}| \\
&\leq &\frac{\Lambda }{6}|\Delta _{i}b|^{2}+C|\Delta _{i}b|^{3} \\
&\leq &(C\epsilon +\frac{\Lambda }{6})|\Delta _{i}b|^{2}\leq (C\epsilon +%
\frac{\Lambda }{6})\epsilon ^{2}\leq \frac{1}{2}
\end{eqnarray*}
for $\epsilon $ sufficiently small. So, using Lemma \ref{l.8.1} of the
Appendix, $W_{\mathcal{P}}$ satisfies the estimate, 
\begin{eqnarray*}
|W_{\mathcal{P}}(\sigma )| &\leq &\sum_{i=1}^{n}(|{\text{tr}}{\mathcal{E}}%
_{i-1}|+|\Psi (-U_{i-1})|) \\
&\leq &d\sum_{i=1}^{n}(|\mathcal{E}_{i-1}|+|U_{i-1}|^{2}(1-|U_{i-1}|)^{-1})
\\
&\leq &d\sum_{i=1}^{n}\left[ C|\Delta _{i}b|^{3}+2\left( (C\epsilon +\frac{%
\Lambda }{6})|\Delta _{i}b|^{2}\right) ^{2}\right] \\
&\leq &C_{1}\sum_{i=1}^{n}|\Delta _{i}b|^{3}.
\end{eqnarray*}
\end{proof}

Let $\mathcal{S}_{\mathcal{P}}:\mathrm{H}_{\mathcal{P}}(M)\rightarrow 
\mathbb{R}$ be given as 
\begin{equation}
\mathcal{S}_{\mathcal{P}}(\sigma ):=\sum_{i=1}^{n}\mathrm{Scal}(\sigma
(s_{i-1}))\Delta _{i}s,  \label{e.6.11}
\end{equation}
where $\mathrm{Scal}$ is the scalar curvature of $(M,\langle \cdot ,\cdot
\rangle ).$

\begin{prop}
\label{p.6.4}Let $p\in \mathbb{R}$ and $\epsilon >0.$ Then there exists $%
C=C(p,\epsilon ,M)<\infty $ such that 
\begin{equation}
1-Ce^{-\frac{\epsilon ^{2}}{4|\mathcal{P}|}}\leq \int_{\mathrm{H}_{\mathcal{P%
}}^{\epsilon }(M)}e^{p(\mathcal{R}_{\mathcal{P}}(\sigma )-\mathcal{S}_{%
\mathcal{P}}(\sigma ))}d\nu _{\mathcal{P}}^{1}(\sigma )\leq e^{CK^{2}|%
\mathcal{P}|}-Ce^{-\frac{\epsilon ^{2}}{4|\mathcal{P}|}},  \label{e.6.12}
\end{equation}
and hence 
\begin{equation}
\left| \int_{\mathrm{H}_{\mathcal{P}}^{\epsilon }(M)}e^{p(\mathcal{R}_{%
\mathcal{P}}(\sigma )-\mathcal{S}_{\mathcal{P}}(\sigma ))}d\nu _{\mathcal{P}%
}^{1}(\sigma )\right| \leq e^{CK^{2}|\mathcal{P}|}-Ce^{-\frac{\epsilon ^{2}}{%
4|\mathcal{P}|}}  \label{e.6.13}
\end{equation}
for all partitions $\mathcal{P}$ with $|\mathcal{P}|$ sufficiently small.
\end{prop}

\begin{proof}
Let $u_{\mathcal{P}}$ be the solution to equation (\ref{e.2.6}) with $b$
replaced by $B_{\mathcal{P}},$ $\mathbf{R}_{i}:=\mathrm{Ric}_{u_{\mathcal{P}%
}(s_{i-1})},$ and 
\begin{equation*}
Y:=e^{p\sum_{i=1}^{n}(\langle \mathbf{R}_{i}\Delta _{i}B,\Delta _{i}B\rangle
-{\text{tr}}(\mathbf{R}_{i})\Delta _{i}s)}.
\end{equation*}
By Theorem \ref{t.4.10}, the distribution of $e^{p(\mathcal{R}_{\mathcal{P}}-%
\mathcal{S}_{\mathcal{P}})}$ under $\nu _{\mathcal{P}}^{1}$ is the same as
the distribution of $Y$ under $\mu .$ Therefore, 
\begin{equation*}
\int_{\mathrm{H}_{\mathcal{P}}^{\epsilon }(M)}e^{p(\mathcal{R}_{\mathcal{P}%
}(\sigma )-\mathcal{S}_{\mathcal{P}}(\sigma ))}d\nu _{\mathcal{P}%
}^{1}(\sigma )=\int_{\mathcal{A}^{c}}Yd\mu ,
\end{equation*}
where $\mathcal{A}:=\cup _{i=1}^{n}\mathcal{A}_{i}$ and $\mathcal{A}%
_{i}:=\{|\Delta _{i}B|\geq \epsilon \}$ as in the proof of Proposition \ref
{p.5.14}. By Proposition \ref{p.8.8} of the Appendix 
\begin{equation*}
1\leq \int_{\mathrm{W}(\mathbb{R}^{d})}Yd\mu =\int_{\mathcal{A}^{c}}Yd\mu
+\int_{\mathcal{A}}Yd\mu \leq e^{dp^{2}K^{2}|\mathcal{P}|},
\end{equation*}
where $K$ is a bound on $\mathrm{Ric}.$ Therefore, 
\begin{equation*}
1-\int_{\mathcal{A}}Yd\mu \leq \int_{\mathrm{H}_{\mathcal{P}}^{\epsilon
}(M)}e^{p(\mathcal{R}_{\mathcal{P}}(\sigma )-\mathcal{S}_{\mathcal{P}%
}(\sigma ))}d\nu _{\mathcal{P}}^{1}(\sigma )\leq e^{dp^{2}K^{2}|\mathcal{P}%
|}-\int_{\mathcal{A}}Yd\mu .
\end{equation*}
So to finish the proof it suffices to show that $\int_{\mathcal{A}}Yd\mu
\leq C\exp \left( -\frac{\epsilon ^{2}}{4|\mathcal{P}|}\right) $.

Noting that 
\begin{equation*}
\left| \sum_{i=1}^{n}(\langle \mathbf{R}_{i}\Delta _{i}B,\Delta _{i}B\rangle
-{\text{tr}}(\mathbf{R}_{i})\Delta _{i}s)\right| \leq K\left(
\sum_{i=1}^{n}|\Delta _{i}B|^{2}+d\right) ,
\end{equation*}
\begin{align}
\int_{\mathcal{A}}Yd\mu & \leq \int_{\mathcal{A}}\exp \left(
K|p|(\sum_{j=1}^{n}|\Delta _{i}B|^{2}+d)\right) d\mu  \notag \\
& \leq \sum_{i}\int_{\mathcal{A}_{i}}\exp \left( K|p|(\sum_{j=1}^{n}|\Delta
_{i}B|^{2}+d)\right) d\mu  \notag \\
& \leq \sum_{i}\mathbb{E}\left[ \exp \left( K|p|(\sum_{j=1}^{n}|\Delta
_{i}B|^{2}+d)\right) \right] \mathbb{E}\left[ \chi _{\epsilon }(|\Delta
_{i}B|)e^{K|p||\left( \Delta _{i}B|^{2}+d\right) }\right] ,  \label{e.6.14}
\end{align}
were $\chi _{\epsilon }(r)=1_{r\geq \epsilon }.$ The first factor of each
term in the sum is bounded by Lemma \ref{l.8.5}. Using the same type of
argument as in the proof of Lemma \ref{l.8.6} one shows for $|\mathcal{P}|$
sufficiently small that there is a constant $C<\infty $ such that 
\begin{equation*}
\mathbb{E}\left[ \chi _{\epsilon }(|\Delta _{i}B|)e^{K|p|\left( |\Delta
_{i}B|^{2}+d\right) }\right] =\mathbb{E}\left[ \chi _{\epsilon }(\sqrt{%
\Delta _{i}s}|B(1)|)e^{K|p|\left( \Delta _{i}s|B(1)|^{2}+d\right) }\right]
\leq C(\Delta _{i}s)e^{-\frac{\epsilon ^{2}}{4|\mathcal{P}|}}.
\end{equation*}
Hence the sum in equation \ref{e.6.14} may be estimated to give $\int_{%
\mathcal{A}}Yd\mu \leq C\exp \left( -\frac{\epsilon ^{2}}{4|\mathcal{P}|}%
\right) .$
\end{proof}

\begin{cor}
\label{c.6.5}Let $\mathcal{S}_{\mathcal{P}}:\mathrm{H}_{\mathcal{P}%
}(M)\rightarrow \mathbb{R}$ be given as in equation (\ref{e.6.11}). Then for
all $\epsilon >0 $ sufficiently small there is a constant $C=C(\epsilon )$
such that 
\begin{equation}
\int_{\mathrm{H}_{\mathcal{P}}^{\epsilon }(M)}\left| \rho _{\mathcal{P}}-e^{-%
\frac{1}{6}\mathcal{S}_{\mathcal{P}}}\right| d\nu _{\mathcal{P}}^{1}\leq C%
\sqrt{|\mathcal{P}|}  \label{e.6.15}
\end{equation}
for all partitions $\mathcal{P}$ with $|\mathcal{P}|$ sufficiently small.
\end{cor}

\begin{proof}
Let $C$ be a generic constant depending on the geometry and the dimension of 
$M.$ Let $J$ denote the left side of equation (\ref{e.6.15}) and let $K$ be
a constant so that $|\mathrm{Scal}|\leq K$. Then 
\begin{align*}
J& =\int_{\mathrm{H}_{\mathcal{P}}^{\epsilon }(M)}\left| \rho _{\mathcal{P}%
}-e^{-\frac{1}{6}\mathcal{S}_{\mathcal{P}}}\right| d\nu _{\mathcal{P}}^{1} \\
& =\int_{\mathrm{H}_{\mathcal{P}}^{\epsilon }(M)}\left| e^{-\frac{1}{6}%
\mathcal{R}_{\mathcal{P}}}e^{W_{\mathcal{P}}}-e^{-\frac{1}{6}\mathcal{S}_{%
\mathcal{P}}}\right| d\nu _{\mathcal{P}}^{1} \\
& \leq e^{K}\int_{\mathrm{H}_{\mathcal{P}}^{\epsilon }(M)}\left| e^{-\frac{1%
}{6}(\mathcal{R}_{\mathcal{P}}-\mathcal{S}_{\mathcal{P}})}e^{W_{\mathcal{P}%
}}-1\right| d\nu _{\mathcal{P}}^{1}\leq I+II,
\end{align*}
where 
\begin{equation*}
I:=e^{K}\int_{\mathrm{H}_{\mathcal{P}}^{\epsilon }(M)}\left| e^{-\frac{1}{6}(%
\mathcal{R}_{\mathcal{P}}-\mathcal{S}_{\mathcal{P}})}-1\right| e^{W_{%
\mathcal{P}}}d\nu
\end{equation*}
and 
\begin{equation*}
II:=e^{K}\int_{\mathrm{H}_{\mathcal{P}}^{\epsilon }(M)}\left| e^{W_{\mathcal{%
P}}}-1\right| d\nu .
\end{equation*}
Since $|e^{a}-1|\leq e^{|a|}-1\leq |a|e^{|a|}$ for all $a\in \mathbb{R},$ 
\begin{equation}
\int_{\mathrm{H}_{\mathcal{P}}^{\epsilon }(M)}\left| e^{W_{\mathcal{P}%
}}-1\right| d\nu _{\mathcal{P}}^{1}\leq \int_{\mathrm{H}_{\mathcal{P}%
}^{\epsilon }(M)}|W_{\mathcal{P}}|e^{|W_{\mathcal{P}}|}d\nu _{\mathcal{P}%
}^{1}  \label{e.6.16}
\end{equation}
By Proposition \ref{p.6.3} there exist $\epsilon _{0}>0$ such that $|W_{%
\mathcal{P}}(\sigma )|\leq C\sum_{i=1}^{n}|\Delta _{i}b|^{3}$ on $\mathrm{H}%
_{\mathcal{P}}^{\epsilon }(M))$ for $\epsilon <\epsilon _{0}$. Therefore,
with the aid of Theorem \ref{t.4.10}, 
\begin{eqnarray*}
\int_{\mathrm{H}_{\mathcal{P}}^{\epsilon }(M)}|W_{\mathcal{P}}|e^{|W_{%
\mathcal{P}}|}d\nu _{\mathcal{P}}^{1} &\leq &C\sum_{i=1}^{n}\int_{\mathrm{H}%
_{\mathcal{P}}^{\epsilon }(M)}|\Delta _{i}b|^{3}\exp \left( C\epsilon
_{0}\sum_{i=1}^{n}|\Delta _{i}b|^{2}\right) d\nu _{\mathcal{P}}^{1} \\
&\leq &C\sum_{i=1}^{n}\int_{\mathrm{H}_{\mathcal{P}}(M)}|\Delta
_{i}b|^{3}\exp \left( C\epsilon _{0}\sum_{i=1}^{n}|\Delta _{i}b|^{2}\right)
d\nu _{\mathcal{P}}^{1} \\
&=&C\sum_{i=1}^{n}\int_{\mathrm{W}(\mathbb{R}^{d})}|\Delta _{i}B|^{3}\exp
\left( C\epsilon _{0}\sum_{i=1}^{n}|\Delta _{i}B|^{2}\right) d\mu \\
&=&C\sum_{i=1}^{n}\mathbb{E}\left[ |\Delta _{i}B|^{3}\exp \left( C\epsilon
_{0}|\Delta _{i}B|^{2}\right) \right] \mathbb{E}\exp \left( C\epsilon
_{0}\sum_{i:i\neq j}^{n}|\Delta _{i}B|^{2}\right) .
\end{eqnarray*}
By Lemma \ref{l.8.5}, $\lim \sup_{|\mathcal{P}|\rightarrow 0}\mathbb{E}\left[
e^{C\epsilon _{0}\sum_{i:i\neq j}^{n}|\Delta _{i}B|^{2}}\right]
=e^{dC\epsilon _{0}}<\infty $ and hence 
\begin{eqnarray}
II &\leq &2e^{K}Ce^{dC\epsilon _{0}}\sum_{i=1}^{n}\mathbb{E}\left[ |\Delta
_{i}B|^{3}\exp \left( C\epsilon _{0}|\Delta _{i}B|^{2}\right) \right]  \notag
\\
&=&2e^{K}Ce^{dC\epsilon _{0}}\sum_{i=1}^{n}(\Delta _{i}s)^{3/2}\mathbb{E}%
\left[ |B(1)|^{3}\exp \left( C\epsilon _{0}\Delta _{i}s|B(1)|^{2}\right) %
\right]  \notag \\
&\leq &2e^{K}Ce^{dC\epsilon _{0}}\mathbb{E}\left[ |B(1)|^{3}\exp \left(
C\epsilon _{0}|\mathcal{P}||B(1)|^{2}\right) \right] \sqrt{|\mathcal{P}|} 
\notag \\
&\leq &C\sqrt{|\mathcal{P}|}.  \label{e.6.17}
\end{eqnarray}
for all partitions $\mathcal{P}$ with $|\mathcal{P}|$ sufficiently small.

To estimate $I$, apply Holder's inequality to get 
\begin{equation*}
I^{2}\leq e^{2K}\left( \int_{\mathrm{H}_{\mathcal{P}}^{\epsilon }(M)}\left|
e^{-\frac{1}{6}(\mathcal{R}_{\mathcal{P}}-\mathcal{S}_{\mathcal{P}%
})}-1\right| ^{2}d\mu \right) \left( \int_{\mathrm{H}_{\mathcal{P}%
}^{\epsilon }(M)}e^{2|W_{\mathcal{P}}|}d\mu \right) .
\end{equation*}
The second term is bounded by the above arguments. Expanding the square
gives 
\begin{eqnarray*}
\left| e^{-\frac{1}{6}(\mathcal{R}_{\mathcal{P}}-\mathcal{S}_{\mathcal{P}%
})}-1\right| ^{2} &=&(e^{-\frac{1}{3}(\mathcal{R}_{\mathcal{P}}-\mathcal{S}_{%
\mathcal{P}})}-1)-2(e^{-\frac{1}{6}(\mathcal{R}_{\mathcal{P}}-\mathcal{S}_{%
\mathcal{P}})}-1) \\
&\leq &\left| e^{-\frac{1}{3}(\mathcal{R}_{\mathcal{P}}-\mathcal{S}_{%
\mathcal{P}})}-1\right| +2\left| e^{-\frac{1}{6}(\mathcal{R}_{\mathcal{P}}-%
\mathcal{S}_{\mathcal{P}})}-1\right| .
\end{eqnarray*}
By equation (\ref{e.6.13}) of Proposition \ref{p.6.4} to each term above,
there is a constant $C=C(\epsilon ,M)<\infty ,$ such that 
\begin{equation*}
I^{2}\leq Ce^{2K}\left( e^{CK^{2}|\mathcal{P}|}-Ce^{-\frac{\epsilon ^{2}}{4|%
\mathcal{P}|}}\right) \leq C|\mathcal{P}|
\end{equation*}
for all partitions $\mathcal{P}$ with $|\mathcal{P}|$ sufficiently small.
>From this we see that 
\begin{equation*}
I\leq C|\mathcal{P}|^{1/2}
\end{equation*}
which together with (\ref{e.6.17}) proves the Corollary.
\end{proof}

\subsection{Proof of Theorem \ref{t.6.1}}

\label{s.6.2}

To simplify notation, let $\rho :W(M)\rightarrow (0,\infty )$ be given by 
\begin{equation}
\rho (\sigma ):=\exp \left( -\frac{1}{6}\int_{0}^{1}\mathrm{Scal}(\sigma
(s))\,ds\right) ,  \label{e.6.18}
\end{equation}
where $\mathrm{Scal}$ is the scalar curvature of $(M,g).$ Recall, by the
remark following Theorem \ref{t.6.1}, to prove Theorem\ref{t.6.1} it
suffices to prove equation (\ref{e.6.1}) for some $\epsilon >0.$ Let $F:%
\mathrm{W}(O(M))\rightarrow \mathbb{R},$ $f:\mathrm{H}(M)\rightarrow \mathbb{%
R},$ and $\tilde{f}:\mathrm{W}(M)\rightarrow \mathbb{R}$ be as in the
statement of Theorem \ref{t.6.1}. Then by Corollary \ref{c.6.5}\ and
Proposition \ref{p.5.13}, for $\epsilon >0$ sufficiently small and for
partitions $\mathcal{P}$ with $|\mathcal{P}|$ sufficiently small, 
\begin{eqnarray*}
\int_{\mathrm{H}_{\mathcal{P}}^{\epsilon }(M)}fd\nu _{\mathcal{P}}^{0}
&=&\int_{\mathrm{H}_{\mathcal{P}}^{\epsilon }(M)}f\rho _{\mathcal{P}}d\nu _{%
\mathcal{P}}^{1} \\
&=&\int_{\mathrm{H}_{\mathcal{P}}^{\epsilon }(M)}fe^{-\frac{1}{6}\mathcal{S}%
_{\mathcal{P}}}d\nu _{\mathcal{P}}^{1}+\tilde{\epsilon}_{\mathcal{P}} \\
&=&\int_{\mathrm{H}_{\mathcal{P}}(M)}fe^{-\frac{1}{6}\mathcal{S}_{\mathcal{P}%
}}d\nu _{\mathcal{P}}^{1}+\epsilon _{\mathcal{P}},
\end{eqnarray*}
and $|\epsilon _{\mathcal{P}}|\leq C\Vert f\Vert _{\infty }|\mathcal{P}%
|^{1/2}$ where $C$ is a constant independent of $\mathcal{P}.$ Because of
Theorem \ref{t.4.17}, to finish the proof, it suffices to show that 
\begin{equation*}
\lim_{|\mathcal{P}|\rightarrow 0}\int_{\mathrm{H}_{\mathcal{P}}(M)}f(e^{-%
\frac{1}{6}\mathcal{S}_{\mathcal{P}}}-\rho )d\nu _{\mathcal{P}}^{1}=0.
\end{equation*}
As above, let $B$ be the $\mathbb{R}^{d}$ --Brownian motion in Notation \ref
{n.1.2}, $B_{\mathcal{P}}$ be its piecewise linear approximation, $\sigma _{%
\mathcal{P}}=\phi (B_{\mathcal{P}})$ and $u_{\mathcal{P}}:=//(\sigma _{%
\mathcal{P}})$. If $\Lambda $ is a constant such that $|\mathrm{Scal}| \leq
\Lambda $ and $|\nabla \mathrm{Scal}|\leq \Lambda ,$ then 
\begin{eqnarray}
&&\left| \int_{\mathrm{H}_{\mathcal{P}}(M)}f(e^{-\frac{1}{6}\mathcal{S}_{%
\mathcal{P}}}-\rho )d\nu _{\mathcal{P}}^{1}\right|  \notag \\
&&\qquad \leq \mathbb{E}\left[ \left| f(u_{\mathcal{P}})\left( e^{-\frac{1}{6%
}\int_{0}^{1}\mathrm{Scal}(\sigma _{\mathcal{P}}(\underline{s}))ds}-e^{-%
\frac{1}{6}\int_{0}^{1}\mathrm{Scal}(\sigma _{\mathcal{P}}(s))ds}\right)
\right| \right]  \notag \\
&&\qquad \leq \Vert f\Vert _{\infty }e^{\Lambda /6}\mathbb{E}\left[
\int_{0}^{1}|\mathrm{Scal}(\sigma _{\mathcal{P}}(s))-\mathrm{Scal}(\sigma _{%
\mathcal{P}}(\underline{s}))|ds\right]  \label{e.6.19}
\end{eqnarray}
wherein the last step we used the inequality $|e^{a}-e^{b}|\leq e^{\max
(a,b)}|a-b|$. For $s\in \lbrack s_{i-1},s_{i})$, we have 
\begin{equation*}
|\mathrm{Scal}(\sigma _{\mathcal{P}}(s))-\mathrm{Scal}(\sigma _{\mathcal{P}%
}(s_{i-1}))|\leq \Lambda |\Delta _{i}B|
\end{equation*}
and hence 
\begin{eqnarray*}
\left| \int_{\mathrm{H}_{\mathcal{P}}(M)}f(e^{-\frac{1}{6}\mathcal{S}_{%
\mathcal{P}}}-\rho )d\nu _{\mathcal{P}}^{1}\right| &\leq &\Vert f\Vert
_{\infty }e^{\Lambda /6}\Lambda \sum_{i=1}^{n}\mathbb{E}|\Delta _{i}B|\Delta
_{i}s \\
&=&\Vert f\Vert _{\infty }e^{\Lambda /6}\Lambda \mathbb{E}%
|B(1)|\sum_{i=1}^{n}\left( \Delta _{i}s\right) ^{3/2} \\
&\leq &C\Vert f\Vert _{\infty }|\mathcal{P}|^{1/2}.
\end{eqnarray*}
This finishes the proof of Theorem \ref{t.6.1}.\qed

\begin{definition}
\label{d.6.6}Let $\mathcal{P}$ be a partition of $[0,1].$ To every point $%
\mathbf{x}\in M^{\mathcal{P}}$ we will associate a path $\sigma _{\mathbf{x}%
}\in \mathrm{H}_{\mathcal{P}}(M)$ as follows. If for each $i,$ there is a
unique minimal geodesic joining $x_{i-1}$ to $x_{i},$ let $\sigma _{\mathbf{x%
}}$ be the uniqe path in $\mathrm{H}_{\mathcal{P}}(M)$ such that $\sigma _{%
\mathbf{x}}(s_{i})=x_{i}$ and $\int_{s_{i-1}}^{s_{i}}|\sigma ^{\prime
}(s)|ds=d(x_{i-1},x_{i})$ for $i=1,2,\ldots ,n$. Otherwise set $\sigma _{%
\mathbf{x}}(s):=o$ for all $s.$
\end{definition}

\begin{cor}
\label{c.6.7}Let $\alpha \in \lbrack 0,1],$ $F:W(O(M))\rightarrow \mathbb{R}$
be a continuous and bounded function and set $f(\sigma ):=F(//_{\cdot
}(\sigma ))$ for $\sigma \in \mathrm{H}(M).$ Then for $\alpha \in \lbrack
0,1],$ 
\begin{equation*}
\lim_{|\mathcal{P}|\rightarrow 0}\int_{M^{\mathcal{P}}}f(\sigma _{\mathbf{x}%
})e^{\frac{1}{6}\sum_{i=1}^{n}\left( \alpha \mathrm{Scal}(x_{i-1})+(1-\alpha
)\mathrm{Scal}(x_{i})\right) \Delta _{i}s}d\gamma _{\mathcal{P}}(\mathbf{x}%
)=\int_{W(M)}\tilde{f}(\sigma )d\nu (\sigma ),
\end{equation*}
where $\tilde{f}(\sigma ):=F(/\tilde{/}_{\cdot }(\sigma ))$ and $/\tilde{/}%
_{r}(\sigma )$ is stochastic parallel translation, see Definition \ref
{d.4.15}.
\end{cor}

\begin{rem}
\label{r.6.7b}It would be possible to prove this corollary using standard
weak convergence arguments appropriatly generalized to the manifold setting,
see for example Section 10 in Stroock and Varadhan \cite{StVar69b} and \cite
{StVar79} and Ethier and Kurtz \cite{EtKur86}.
\end{rem}

\begin{proof}
For $\sigma \in \mathrm{H}(M)$, let $\chi _{\mathcal{P},\alpha }(\sigma )=e^{%
\frac{1}{6}\sum_{i=1}^{n}\left( \alpha \mathrm{Scal}(\sigma
(s_{i-1}))+(1-\alpha )\mathrm{Scal}(\sigma (s_{i}))\right) \Delta _{i}s}$.
Let $\Lambda $ be a constant such that $|\mathrm{Scal}|\leq \Lambda $ and $%
|\nabla \mathrm{Scal}|\leq \Lambda $. Then $\chi _{\mathcal{P},\alpha
}(\sigma )\leq e^{\Lambda /6}$ so by Proposition \ref{p.5.12} 
\begin{equation*}
\int_{M^{\mathcal{P}}\setminus M_{\epsilon }^{\mathcal{P}}}f(\sigma _{%
\mathbf{x}})\chi _{\mathcal{P},\alpha }(\sigma _{\mathbf{x}})d\gamma _{%
\mathcal{P}}(\mathbf{x})=\epsilon _{\mathcal{P}}
\end{equation*}
where $\epsilon _{\mathcal{P}}\leq C||f||_{\infty }|\mathcal{P}|^{1/2}$.
Therefore it is sufficient to consider $\int_{M_{\epsilon }^{\mathcal{P}%
}}f(\sigma _{\mathbf{x}})\chi _{\mathcal{P},\alpha }(\sigma _{\mathbf{x}%
})d\gamma _{\mathcal{P}}(\mathbf{x})$. By Propositions \ref{p.5.6} we have 
\begin{equation*}
\int_{M_{\epsilon }^{\mathcal{P}}}f(\sigma _{\mathbf{x}})\chi _{\mathcal{P}%
,\alpha }(\sigma _{\mathbf{x}})d\gamma _{\mathcal{P}}(\mathbf{x})=\int_{%
\mathrm{H}_{\mathcal{P}}^{\epsilon }(M)}f(\sigma )\chi _{\mathcal{P},\alpha
}(\sigma )d\nu _{\mathcal{P}}^{0}(\sigma )
\end{equation*}
Let $\rho (\sigma )$ be given by (\ref{e.6.18}). Arguing as in the proof of
Theorem \ref{t.6.1}, the Corollary will follow if 
\begin{equation*}
\lim_{|\mathcal{P}|\rightarrow 0}\int_{\mathrm{H}_{\mathcal{P}}(M)}f(\sigma
)(\chi _{\mathcal{P},\alpha }(\sigma )\rho (\sigma )-1)d\nu _{\mathcal{P}%
}^{1}(\sigma )=0
\end{equation*}
Let $\sigma _{\mathcal{P}},B_{\mathcal{P}}$ be as in the proof of Theorem 
\ref{t.6.1}. We estimate as in the proof of Theorem \ref{t.6.1}, 
\begin{multline*}
\left| \lim_{|\mathcal{P}|\rightarrow 0}\int_{\mathrm{H}_{\mathcal{P}%
}(M)}f(\sigma )(\chi _{\mathcal{P},\alpha }(\sigma )\rho (\sigma )-1)d\nu _{%
\mathcal{P}}^{1}(\sigma )\right| \\
\leq ||f||_{\infty }e^{\Lambda /6}\mathbb{E}\left[ \left|
\sum_{i=1}^{n}\left( \alpha \mathrm{Scal}(\sigma _{\mathcal{P}%
}(s_{i-1}))+(1-\alpha )\mathrm{Scal}(\sigma _{\mathcal{P}}(s_{i}))\right)
\Delta _{i}s-\int_{0}^{1}\mathrm{Scal}(\sigma _{\mathcal{P}}(s))ds\right| %
\right] \\
\leq C||f||_{\infty }|\mathcal{P}|^{1/2}
\end{multline*}
which completes the proof of Corollary \ref{c.6.7}.
\end{proof}

\section{Partial Integration Formulas\label{s.7}}

As an application of Theorem \ref{t.4.17}, we will derive the known
integration by parts formula for the measure $\nu .$ This will be
accomplished by taking limits of the finite dimensional integration by parts
formulas for the measure $\nu _{\mathcal{P}}^{1}.$ The main result appears
at the end of this section in Theorem \ref{t.7.16}. A similar method for
proving integration by parts formula for laws of solutions to stochastic
differential equations has been used by Bell \cite{Bell87,Bell95}.

\subsection{Integration by parts for the approximate measures\label{s.7.1}}

The two ingredients for computing the integration by parts formula for the
form $\nu _{\mathcal{P}}^{1}$ is the differential of $E$ and the Lie
derivative of $\mathrm{Vol}_{G_{\mathcal{P}}^{1}}.$ The following lemma may
be found in any book on Riemannian geometry. We will supply the short proof
for the readers convenience.

\begin{lemma}
\label{l.7.1} Let $Y\in T_{\sigma }\mathrm{H}(M),$ then 
\begin{equation}
YE=dE(Y)=2\int_{0}^{1}\langle \sigma ^{\prime }(s),\frac{\nabla Y(s)}{ds}%
\rangle ds.  \label{e.7.1}
\end{equation}
\end{lemma}

\begin{proof}
Choose a one parameter family of paths $\sigma _{t}\in \mathrm{H}(M)$ such
that $\sigma _{0}=\sigma $ and $\frac{d}{dt}|_{t=0}\sigma _{t}=Y.$ Then 
\begin{equation*}
YE=\frac{d}{dt}|_{t=0}\int_{0}^{1}|\sigma _{t}^{\prime
}(s)|^{2}ds=2\int_{0}^{1}\langle \frac{\nabla }{dt}\sigma _{t}^{\prime
}(s)|_{t=0},\sigma ^{\prime }(s)\rangle ds.
\end{equation*}
Since $\nabla $ has zero torsion, 
\begin{equation*}
\frac{\nabla }{dt}\sigma _{t}^{\prime }(s)|_{t=0}=\frac{\nabla }{ds}\frac{d}{%
dt}|_{t=0}\sigma _{t}(s)=\frac{\nabla }{ds}Y(s).
\end{equation*}
The last two equations clearly imply equation (\ref{e.7.1}).
\end{proof}

To compute the Lie derivative of $\mathrm{Vol}_{G_{\mathcal{P}}^{1}} $ is
will be useful to have an orthonormal frame on $\mathrm{H}_{\mathcal{P}}(M)$
relative to $G_{\mathcal{P}}^{1}.$ We will construct such a frame in the
next lemma.

\begin{notation}
\label{n.7.2} Given $\sigma \in \mathrm{H}_{\mathcal{P}}(M)$, let $\mathrm{H}%
_{\mathcal{P},\sigma }$be the subspace of $\mathrm{H}$ given by 
\begin{equation}
\mathrm{H}_{\mathcal{P},\sigma }:=\{v\in \mathrm{H}:v^{\prime \prime
}(s)=\Omega _{u(s)}(b^{\prime }(s),v(s))b^{\prime }(s),\quad \forall s\notin 
\mathcal{P}\},  \label{e.7.2}
\end{equation}
where $u=//(\sigma )$ and $b=\phi ^{-1}(\sigma )$.\qed
\end{notation}

Because of equation (\ref{e.4.2}) of Proposition \ref{p.4.4}, $v\in \mathrm{H%
}_{\mathcal{P},\sigma }$ if and only if $X^{v}(\sigma ):=//(\sigma )v\in
T_{\sigma }\mathrm{H}_{\mathcal{P}}(M)$.

\begin{lemma}[$G_{\mathcal{P}}$--orthonormal frame]
\label{l.7.3} Let $\mathcal{P}$ be a partition of $[0,1]$ and $G_{\mathcal{P}%
}^{1}$ be as in equation (\ref{e.1.12}) above. Also let $\{e_{a}\}_{a=1}^{d}$
be an orthonormal frame for $T_{o}M\cong \mathbb{R}^{d}$. For $\sigma \in 
\mathrm{H}_{\mathcal{P}}(M),$ $i=1,2,\ldots ,n$ and $a=1,\ldots ,d$ let $%
h_{i,a}(s,\sigma ):=v(s)$ be determined (uniquely) by:

\begin{enumerate}
\item  $v\in \mathrm{H}_{\mathcal{P},\sigma }$.

\item  $v^{\prime }(s_{j}+)=0$ if $j\neq i-1$.

\item  $v^{\prime }(s_{i-1}+)=\frac{1}{\sqrt{\Delta _{i}s}}e_{a}.$
\end{enumerate}

Then $\{X^{h_{a,i}},i=1,\ldots ,n,\quad a=1,\ldots ,d\}$ is a globally
defined orthonormal frame for $(\mathrm{H}_{\mathcal{P}}(M),G_{\mathcal{P}%
}^{1})$.
\end{lemma}

\begin{proof}
This lemma is easily verified using the definition of $G_{\mathcal{P}}^{1}$
in equation (\ref{e.1.12}), the identity 
\begin{equation*}
\frac{\nabla X^{v}(\sigma )(s+)}{ds}=//_{s}(\sigma )v^{\prime }(s+),
\end{equation*}
and the fact that $//_{s}(\sigma )$ is an isometry.
\end{proof}

\begin{definition}
\label{d.7.4} Let $PC^{1}$ denote the set of $k\in \mathrm{H}$ which are
piecewise $C^{1}.$ Given $k\in PC^{1},$ define $k_{\mathcal{P}}:\mathrm{H}_{%
\mathcal{P}}(M)\rightarrow \mathrm{H}$ by requiring $k_{\mathcal{P}}(\sigma
)\in \mathrm{H}_{\mathcal{P},\sigma }$ for all $\sigma \in \mathrm{H}_{%
\mathcal{P}}(M)$ and $k_{\mathcal{P}}^{\prime }(\sigma ,s+)=k^{\prime }(s+)$
for all $s\in \mathcal{P}\setminus \{1\}.$ Note that with this definition of 
$k_{\mathcal{P}},$ $X^{k_{\mathcal{P}}}$ is the unique tangent vector field
on $\mathrm{H}_{\mathcal{P}}(M)$ such that 
\begin{equation*}
\frac{\nabla X^{k_{\mathcal{P}}}(s+)}{ds}=\frac{\nabla X^{k}(s+)}{ds}\quad 
\text{ for all \thinspace\ }s\in \mathcal{P}\setminus \{1\}.
\end{equation*}
\end{definition}

\begin{lemma}
\label{l.7.5} If $k\in PC^{1},$ then $L_{X^{k_{\mathcal{P}}}}\mathrm{Vol}%
_{G_{\mathcal{P}}^{1}}=0.$
\end{lemma}

\begin{proof}
Recall that on a general Riemannian manifold 
\begin{equation*}
L_{X}\mathrm{Vol}=-\sum_{i}\langle L_{X}e_{i},e_{i}\rangle \mathrm{Vol}%
=\sum_{i}\langle \lbrack e_{i},X],e_{i}\rangle \mathrm{Vol},
\end{equation*}
where $\{e_{i}\}$ is an orthonormal frame. Therefore we must show that 
\begin{equation}
\sum_{i=1}^{n}\sum_{a=1}^{d}G_{\mathcal{P}}^{1}([X^{h_{a,i}},X^{k_{\mathcal{P%
}}}],X^{h_{a,i}})=0.  \label{e.7.3}
\end{equation}
Suppressing $\sigma \in \mathrm{H}_{\mathcal{P}}(M)$ from the notation and
using Theorem \ref{t.3.5} to expand the Lie bracket, we find 
\begin{align*}
G_{\mathcal{P}}^{1}([X^{h_{a,i}},X^{k_{\mathcal{P}}}],X^{h_{a,i}})&
=\sum_{j=1}^{n}\langle (X^{h_{a,i}}k_{\mathcal{P}}-X^{k_{\mathcal{P}%
}}h_{a,i})^{\prime },h_{a,i}^{\prime }\rangle |_{(s_{j-1}+)}\Delta _{j}s \\
& \quad +\sum_{j=1}^{n}\langle (q(X^{k_{\mathcal{P}%
}})h_{a,i}-q(X^{h_{a,i}})k_{\mathcal{P}})^{\prime },h_{a,i}^{\prime }\rangle
|_{(s_{j-1}+)}\Delta _{j}s.
\end{align*}
For $s\in \mathcal{P}\setminus \{1\}$, $(X^{h_{a,i}}k_{\mathcal{P}})^{\prime
}(s+)=X^{h_{a,i}}k_{\mathcal{P}}^{\prime }(s+)=0$, since $k^{\prime }(s+)$
is independent of $\sigma .$ For the same reason, $(X^{k_{\mathcal{P}%
}}h_{a,i})^{\prime }(s+)=0$ for $s\in \mathcal{P}\setminus \{1\}$. Moreover
for $s\in \mathcal{P}\setminus \{1\}$, 
\begin{equation*}
\langle (q(X^{k_{\mathcal{P}}})h_{a,i})^{\prime },h_{a,i}^{\prime }\rangle
|_{s+}=\langle q(X^{k_{\mathcal{P}}})h_{a,i}^{\prime }+R_{u}(\sigma ^{\prime
},X^{k_{\mathcal{P}}})h_{a,i},h_{a,i}^{\prime }\rangle |_{s+}=0.
\end{equation*}
because $q(X^{k_{\mathcal{P}}})$ is skew symmetric and because either $%
h_{a,i}(s+)$ or $h_{a,i}^{\prime }(s+)$ are equal to zero for all $s\in 
\mathcal{P}\setminus \{1\}$. Similarly, 
\begin{equation*}
\langle (q(X^{h_{a,i}})k_{\mathcal{P}})^{\prime },h_{a,i}^{\prime }\rangle
|_{s+}=\langle q(X^{h_{a,i}})k_{\mathcal{P}}^{\prime }+R_{u}(\sigma ^{\prime
},X^{h_{a,i}})k_{\mathcal{P}},h_{a,i}^{\prime }\rangle |_{s+}=0
\end{equation*}
because for all $s\in \mathcal{P}\setminus \{1\}$, either $%
q_{s+}(X^{h_{a,i}})=0$ or $h_{a,i}^{\prime }(s+)=0$ and either $%
h_{a,i}(s+)=0 $ or $h_{a,i}^{\prime }(s+)=0$. Thus every term in the sum in
equation (\ref{e.7.3}) is zero.
\end{proof}

\begin{thm}
\label{t.7.6} Suppose that $k\in PC^{1}$, $\mathcal{P}$ is a partition of $%
[0,1],$ $b\in \mathrm{H}_{\mathcal{P}}$ and $\sigma =\phi (b)\in \mathrm{H}_{%
\mathcal{P}}(M).$ Then 
\begin{equation}
(L_{X^{k_{\mathcal{P}}}}\nu _{\mathcal{P}}^{1})_{\sigma }=-\left(
\sum_{i=1}^{n}\langle k^{\prime }(s_{i-1}+),\Delta _{i}b\rangle \right) (\nu
_{\mathcal{P}}^{1})_{\sigma },  \label{e.7.4}
\end{equation}
i.e. the divergence of $X^{k_{\mathcal{P}}}$ relative to the volume form $%
\nu _{\mathcal{P}}^{1}$ is 
\begin{equation}
(\mathrm{div}_{\nu _{\mathcal{P}}^{1}}X^{k_{\mathcal{P}}})(\sigma
)=-\sum_{i=1}^{n}\langle k^{\prime }(s_{i-1}+),\Delta _{i}b\rangle .
\label{e.7.5}
\end{equation}
\qed
\end{thm}

\begin{proof}
By Lemma \ref{l.7.5}, 
\begin{equation*}
(L_{X^{k_{\mathcal{P}}}}\nu _{\mathcal{P}}^{1})_{\sigma }=\left[ -\frac{1}{2}%
(X^{k_{\mathcal{P}}}E)(\sigma )\right] \cdot (\nu _{\mathcal{P}%
}^{1})_{\sigma }
\end{equation*}
and by Lemma \ref{l.7.1}, 
\begin{eqnarray*}
(X^{k_{\mathcal{P}}}E)(\sigma ) &=&2\int_{0}^{1}\langle \sigma ^{\prime }(s),%
\frac{\nabla X^{k_{\mathcal{P}}}(\sigma )(s)}{ds}\rangle \,ds \\
&=&2\int_{0}^{1}\langle //_{s}(\sigma )b^{\prime }(s),//_{s}(\sigma )k_{%
\mathcal{P}}^{\prime }(\sigma ,s)\rangle \,ds \\
&=&2\sum_{i=1}^{n}\int_{J_{i}}\langle b^{\prime }(s),k_{\mathcal{P}}^{\prime
}(\sigma ,s)\rangle \,ds.
\end{eqnarray*}
Now for $s\in J_{i}:=(s_{i-1},s_{i}],$%
\begin{eqnarray*}
\langle b^{\prime }(s),k_{\mathcal{P}}^{\prime }(\sigma ,s)\rangle
&=&\langle b^{\prime }(s_{i-1}+),k_{\mathcal{P}}^{\prime }(\sigma
,s_{i-1}+)\rangle +\int_{s_{i-1}}^{s}b^{\prime }(r)\cdot k_{\mathcal{P}%
}^{\prime \prime }(\sigma ,r)dr \\
&=&\langle b^{\prime }(s_{i-1}+),k_{\mathcal{P}}^{\prime }(s_{i-1}+)\rangle
\\
&&+\int_{s_{i-1}}^{s}\langle b^{\prime }(s),\Omega _{u(r)}(b^{\prime }(r),k_{%
\mathcal{P}}(\sigma ,r))b^{\prime }(r)\rangle \,dr \\
&=&\langle b^{\prime }(s_{i-1}+),k_{\mathcal{P}}^{\prime }(s_{i-1}+)\rangle ,
\end{eqnarray*}
wherein the last equality we used the skew adjointness of $\Omega
_{u(r)}(b^{\prime }(r),k_{\mathcal{P}}(\sigma ,r))$ and the fact that $%
b^{\prime }(s)=b^{\prime }(r)=\Delta _{i}b/\Delta _{i}s$ for all $s,r\in
J_{i}.$ Combining the previous three displayed equations proves equation (%
\ref{e.7.4}).
\end{proof}

\begin{cor}
\label{c.7.7} Let $k\in PC^{1},$ $\mathcal{P}$ be a partition of $[0,1]$ as
above, and let $f:\mathrm{H}_{\mathcal{P}}(M)\rightarrow \mathbb{R}$ be a $%
C^{1}$ function for which $f$ and its differential is bounded, then 
\begin{equation}
\int_{\mathrm{H}_{\mathcal{P}}(M)}\left( X^{k_{\mathcal{P}}}f\right) \nu _{%
\mathcal{P}}^{1}=\int_{\mathrm{H}_{\mathcal{P}}(M)}f\left(
\sum_{i=1}^{n}\langle k^{\prime }(s_{i-1}+),\Delta _{i}b\rangle \right) \nu
_{\mathcal{P}}^{1}.  \label{e.7.6}
\end{equation}
where in this formula $\Delta _{i}b$ is the to be understood as the function
on $\mathrm{H}(M)$ defined by 
\begin{equation}
\Delta _{i}b(\sigma ):=\phi ^{-1}(\sigma )(s_{i})-\phi ^{-1}(\sigma
)(s_{i-1}).  \label{e.7.7}
\end{equation}
\qed
\end{cor}

\begin{proof}
First assume that $f$ has compact support. Then by Stoke's theorem 
\begin{eqnarray*}
0 &=&\int_{\mathrm{H}_{\mathcal{P}}(M)}d\left[ i_{X^{k_{\mathcal{P}}}}\left(
f\nu _{\mathcal{P}}^{1}\right) \right] =\int_{\mathrm{H}_{\mathcal{P}%
}(M)}L_{X^{k_{\mathcal{P}}}}\left( f\nu _{\mathcal{P}}^{1}\right) \\
&=&\int_{\mathrm{H}_{\mathcal{P}}(M)}\left[ (X^{k_{\mathcal{P}}}f)\nu _{%
\mathcal{P}}^{1}+fL_{X^{k_{\mathcal{P}}}}\nu _{\mathcal{P}}^{1}\right]
\end{eqnarray*}
which combined with equation (\ref{e.7.4}) proves equation (\ref{e.7.6}).
For the general case choose $\chi \in C_{c}^{\infty }(\mathbb{R})$ such that 
$\chi $ is one in a neighborhood of $0.$ Define $\chi _{n}:=\chi (\frac{1}{n}%
E(\cdot ))\in C_{c}^{\infty }(\mathrm{H}_{\mathcal{P}}(M))$ and $f_{n}:=\chi
_{n}f\in C_{c}^{\infty }(\mathrm{H}_{\mathcal{P}}(M)).$ Observe that 
\begin{eqnarray*}
\left( X^{k_{\mathcal{P}}}f_{n}\right) &=&\chi _{n}\cdot X^{k_{\mathcal{P}%
}}f+\frac{1}{n}f\cdot \chi ^{\prime }(\frac{1}{n}E(\cdot ))X^{k_{\mathcal{P}%
}}E \\
&=&\chi _{n}\cdot X^{k_{\mathcal{P}}}f+\frac{1}{n}f\cdot \chi ^{\prime }(%
\frac{1}{n}E(\cdot ))\left( \sum_{i=1}^{n}\langle k^{\prime
}(s_{i-1}+),\Delta _{i}b\rangle \right) ,
\end{eqnarray*}
wherein the last equality we have used the formula for $X^{k_{\mathcal{P}}}E$
computed in the proof of Lemma \ref{t.7.6}. Because of Theorem \ref{t.4.10}, 
$\sum_{i=1}^{n}\langle k^{\prime }(s_{i-1}+),\Delta _{i}b\rangle $ is a
Gaussian random variable on $(\mathrm{H}_{\mathcal{P}}(M),\nu _{\mathcal{P}%
}^{1})$ and hence is in $L^{p}$ for all $p\in \lbrack 1,\infty ).$ Also 
\begin{eqnarray*}
|X^{k_{\mathcal{P}}}f| &\leq &C\sqrt{G_{\mathcal{P}}^{1}(X^{k_{\mathcal{P}%
}},X^{k_{\mathcal{P}}})} \\
&=&C\sqrt{\sum_{i=1}^{n}\langle k^{\prime }(s_{i-1}+),k^{\prime
}(s_{i-1}+)\rangle \Delta _{i}s}\leq C\Vert k^{\prime }\Vert _{\infty ,}
\end{eqnarray*}
where $C$ is bound on the differential of $f.$ Using these remarks and the
dominated convergence theorem, we may replace $f$ by $f_{n}$ in equation \ref
{e.7.6} and pass to the limit to conclude that equation (\ref{e.7.6}) holds
for bounded $f$ with bounded derivatives.
\end{proof}

\begin{rem}
\label{r.7.8}Obviously Corollary \ref{c.7.7} holds for more general
functions $f.$ For example the above proof works if $f$ and $df$ are in $%
L^{p}(\mathrm{H}_{\mathcal{P}}(M),\nu _{\mathcal{P}}^{1})$ for some $p>1.$%
\qed
\end{rem}

We would like to pass to the limit as $|\mathcal{P}|\rightarrow 0$ in
equation (\ref{e.7.6}) of Corollary \ref{c.7.7}. The right side of this
equation is easily dealt with using Theorem \ref{t.4.17}. In order to pass
to the limit on the left side of equation (\ref{e.7.6}) it will be necessary
to understand the limiting behavior of $k_{\mathcal{P}}$ as $|\mathcal{P}%
|\rightarrow 0.$ This is the subject of the next subsection.

\subsection{The limit of $k_{\mathcal{P}}$}

\label{s.7.2}

\begin{notation}
Let $\mathcal{P}=\{0=s_{0}<s_{1}<s_{2}<\cdots <s_{n}=1\}$ be a partition of $%
[0,1]$ and for $r\in (s_{j-1},s_{j}]$, let $\underline{r}:=s_{j-1}$. For $%
k\in PC^{1}$, define $||k^{\prime }||_{1,\mathcal{P}}$ and $|||k^{\prime
}|||_{\mathcal{P}}$ by 
\begin{align}
||k^{\prime }||_{1,\mathcal{P}}& =\sum_{i=1}^{n}|k^{\prime
}(s_{i-1}+)|\Delta _{i}s,  \label{e.7.8} \\
\intertext{and}
|||k^{\prime }|||_{\mathcal{P}}& =\int_{0}^{1}|k^{\prime
}(r)-k^{\prime }(\underline{r})|dr.  \label{e.7.9}
\end{align}
Note that $|||k^{\prime }|||_{\mathcal{P}}=0$ if $k\in \mathrm{H}_{\mathcal{P%
}}.$
\end{notation}

\begin{lemma}
\label{l.7.10} Let $\mathcal{P}$ be a partition of $[0,1]$, $\sigma \in 
\mathrm{H}_{\mathcal{P}}(M)$, $b=\phi ^{-1}(\sigma )$, $u=//(\sigma )$, $%
k\in PC^{1}$ and $k_{\mathcal{P}}(\sigma ,\cdot )$ be as in definition \ref
{d.7.4}. Then with $\Delta _{i}b$ given by (\ref{e.7.7}) and $||k^{\prime
}||_{1,\mathcal{P}}$ given by (\ref{e.7.8}), 
\begin{align}
|k_{\mathcal{P}}(\sigma ,s)|& \leq ||k^{\prime }||_{1,\mathcal{P}}e^{\frac{1%
}{2}\Lambda \sum_{j=1}^{n}|\Delta _{j}b|^{2}}\quad \forall s\in \lbrack 0,1]
\label{e.7.10} \\
\intertext{and}
|k_{\mathcal{P}}(\sigma ,s)-k_{\mathcal{P}}(\sigma ,s_{i-1})|&
\leq \left( |k^{\prime }(s_{i-1}+)|\Delta _{i}s+\frac{1}{2}|k_{\mathcal{P}%
}(\sigma ,s_{i-1})|\Lambda \left| \Delta _{i}b\right| ^{2}\right) \cosh 
\sqrt{\Lambda }\left| \Delta _{i}b\right| ,  \label{e.7.11} \\
\intertext{and}
|k_{\mathcal{P}}(\sigma ,s)-k_{\mathcal{P}}(\sigma ,s_{i-1})|&
\leq |k^{\prime }(s_{i-1}+)|\Delta _{i}s+\frac{1}{2}\Lambda |\Delta
_{i}b|^{2}\Vert k_{\mathcal{P}}(\sigma ,\cdot )\Vert _{\infty }\quad \forall
s\in (s_{i-1},s_{i}],  \label{e.7.12}
\end{align}
where $\Lambda $ is a bound on the curvature tensor.
\end{lemma}

\begin{proof}
Let $\kappa (\cdot ):=k_{\mathcal{P}}(\sigma ,\cdot )\in \mathrm{H}_{%
\mathcal{P},\sigma }$ and $A(s):=\Omega _{u(s)}(b^{\prime }(s),\cdot
)b^{\prime }(s).$ By Definition \ref{d.7.4} of $k_{\mathcal{P}}$, $\kappa $
satisfies 
\begin{align}
\kappa ^{\prime \prime }(s)& =A(s)\kappa (s)\text{ for all }s\notin \mathcal{%
P}  \label{e.7.13} \\
\intertext{and}
\kappa ^{\prime }(s+)& =k^{\prime }(s+)\qquad \forall s\in 
\mathcal{P}\setminus \{1\}.  \label{e.7.14}
\end{align}
Noting that $\left| \Omega _{u(s)}(b^{\prime }(s),\cdot )b^{\prime
}(s)\right| \leq \Lambda \left| b^{\prime }(s)\right| ^{2}=\Lambda \frac{%
\left| \Delta _{i}b\right| ^{2}}{\Delta _{i}s^{2}}$for $s\in
(s_{i-1},s_{i}], $ Lemma \ref{l.8.2} of the Appendix and equation (\ref
{e.7.12}) implies that 
\begin{eqnarray*}
|\kappa (s)-\kappa (s_{i-1})| &\leq &|\kappa (s_{i-1})|\left( \cosh \sqrt{%
\Lambda }\left| \Delta _{i}b\right| -1\right) +|k^{\prime }(s_{i-1}+)|\Delta
_{i}s\frac{\sinh \sqrt{\Lambda }\left| \Delta _{i}b\right| }{\sqrt{\Lambda }%
\left| \Delta _{i}b\right| } \\
&\leq &\left( |k^{\prime }(s_{i-1}+)|\Delta _{i}s+\frac{1}{2}|\kappa
(s_{i-1})|\Lambda \left| \Delta _{i}b\right| ^{2}\right) \cosh \sqrt{\Lambda 
}\left| \Delta _{i}b\right|
\end{eqnarray*}
where we have made use of the elementary inequalities 
\begin{equation}
\cosh (a)-1\leq \frac{1}{2}a^{2}\cosh (a),\quad \text{ and }\frac{\sinh (a)}{%
a}\leq \cosh (a)\quad \forall a\in \mathbb{R}.  \label{e.7.15}
\end{equation}
In particular, equation (\ref{e.7.11}) is valid and 
\begin{align}
|\kappa (s)|& \leq |\kappa (s_{i-1})|\cosh \sqrt{\Lambda }\left| \Delta
_{i}b\right| +|k^{\prime }(s_{i-1}+)|\Delta _{i}s\frac{\sinh \sqrt{\Lambda }%
\left| \Delta _{i}b\right| }{\sqrt{\Lambda }\left| \Delta _{i}b\right| } 
\notag \\
& \leq \left( |\kappa (s_{i-1})|+|k^{\prime }(s_{i-1}+)|\Delta _{i}s\right)
\exp \{\frac{1}{2}\Lambda \left| \Delta _{i}b\right| ^{2}\},  \label{e.7.16}
\end{align}
since $\cosh (a)\leq e^{\alpha ^{2}/2}$ for all $a.$ Using the fact that $%
\kappa (s_{0})=\kappa (0)=0$ and an inductive argument, equation (\ref
{e.7.16}) with $s=s_{i}$ implies that 
\begin{equation*}
|\kappa (s_{i})|\leq \left( \sum_{j=1}^{i}|k^{\prime }(s_{j-1}+)|\Delta
_{j}s\right) e^{\frac{1}{2}\Lambda \sum_{j=1}^{i}|\Delta _{j}b|^{2}}.
\end{equation*}
Combining this last equation with equation (\ref{e.7.16}) proves the bound
in equation (\ref{e.7.10}).
\end{proof}

In the rest of this section, unless otherwise stated, $C$ will be a generic
constant depending only on the geometry of $M$ and $C(\gamma ,p)$ will be a
generic constant depending only on $\gamma ,p$ and the geometry of $M$.

\begin{thm}
\label{t.7.11}Let $k\in PC^{1}$ and $B$ and $B_{\mathcal{P}}$ be the $%
\mathbb{R}^{d}$--valued processes defined in Notation \ref{n.1.2} and
Notation \ref{n.4.12} respectively. Also let $u$ be the $O(M)$--valued
process which solves the Stratonovich stochastic differential equation 
\begin{equation}
\delta u=\mathcal{H}_{u}u\delta B,\qquad u(0)=u_{0},  \label{e.7.17}
\end{equation}
$u_{\mathcal{P}}=//(\phi (B_{\mathcal{P}}))$ and $z_{\mathcal{P}}=k_{%
\mathcal{P}}(\phi (B_{\mathcal{P}}),\cdot ).$ (Note by Theorem \ref{t.4.14}
that $u=\lim_{|\mathcal{P}|\rightarrow \infty }\phi (B_{\mathcal{P}})$ is a
stochastic extension of $\phi .)$ Let $z$ denote the solution to the
(random) ordinary differential equation 
\begin{equation}
z^{\prime }(s)+\frac{1}{2}\mathrm{Ric}_{u(s)}z(s)=k^{\prime }(s),\qquad
z(0)=0.  \label{e.7.18}
\end{equation}
Then for $\gamma \in (0,\frac{1}{2})$, $p\in \lbrack 1,\infty )$, 
\begin{equation*}
\mathbb{E}\left[ \sup_{s\in \lbrack 0,1]}|z_{\mathcal{P}}(s)-z(s)|^{p}\right]
\leq C(\gamma ,p)\left( ||k^{\prime }||_{1,\mathcal{P}}^{p}|\mathcal{P}%
|^{\gamma p}+|||k^{\prime }|||_{\mathcal{P}}^{p}\right) .
\end{equation*}
\end{thm}

We will prove this theorem after the next two lemmas. Before doing this let
us note that $z_{\mathcal{P}}$ in Theorem \ref{t.7.11} above is determined
by 
\begin{equation}
z_{\mathcal{P}}^{\prime \prime }(s)=A(s)z_{\mathcal{P}}(s)\text{ for }%
s\notin \mathcal{P},\text{ }z_{\mathcal{P}}(0)=0,\text{ and }z_{\mathcal{P}%
}^{\prime }(s+)=k^{\prime }(s+)\text{ }\forall s\in \mathcal{P}\setminus
\{1\},  \label{e.7.19}
\end{equation}
where 
\begin{equation}
A(s):=\Omega _{u_{\mathcal{P}}(s)}(\frac{\Delta _{i}B}{\Delta _{i}s},\cdot )%
\frac{\Delta _{i}B}{\Delta _{i}s}\text{ when }s\in (s_{i-1},s_{i}].
\label{e.7.20}
\end{equation}

\begin{lemma}
\label{l.7.12} Let $\delta _{i}$ be defined by 
\begin{equation*}
\delta _{i}:=z_{\mathcal{P}}(s_{i})+\int_{0}^{s_{i}}\left( \frac{1}{2} 
\mathrm{Ric}_{u_{\mathcal{P}}(\underline{r})}z_{\mathcal{P}}(r)-k^{\prime }(%
\underline{r}+)\right) dr.
\end{equation*}
Then there for all $p\in \lbrack 1,\infty )$ and $\gamma \in (0,1/2)$ there
is a constant $C=C(p,\gamma ,\Lambda )<\infty $ such that 
\begin{equation*}
\mathbb{E}\left[ \max_{i}|\delta _{i}|^{p}\right] \leq C||k^{\prime }||_{1,%
\mathcal{P}}^{p}|\mathcal{P}|^{\gamma p},
\end{equation*}
where $\Lambda $ is a bound on $\Omega $ and its horizontal derivative.
\end{lemma}

\begin{proof}
With out loss of generality, we can assume that $p\geq 2.$ Throughout the
proof $C$ will denote generic constant depending only on $p,$ $\gamma ,$ $%
\Lambda ,$ and possibly the dimension of $M.$ By Taylor's theorem with
integral remainder and equation (\ref{e.7.19}) and equation (\ref{e.7.20})
we have 
\begin{align}
z_{\mathcal{P}}(s_{i})& =z_{\mathcal{P}}(s_{i-1})+z_{\mathcal{P}}^{\prime
}(s_{i-1}+)\Delta _{i}s+\int_{s_{i-1}}^{s_{i}}(s_{i}-r)z_{\mathcal{P}%
}^{\prime \prime }(r))dr  \notag \\
& =z_{\mathcal{P}}(s_{i-1})+k^{\prime }(s_{i-1}+)\Delta _{i}s  \notag \\
& +\int_{s_{i-1}}^{s_{i}}(s_{i}-r)\Omega _{u_{\mathcal{P}}(r)}(B_{\mathcal{P}%
}^{\prime }(r),z_{\mathcal{P}}(r))B_{\mathcal{P}}^{\prime }(r)dr  \notag \\
& =z_{\mathcal{P}}(s_{i-1})+k^{\prime }(s_{i-1}+)\Delta _{i}s+\frac{1}{2}%
\Omega _{u_{\mathcal{P}}(s_{i-1})}(\Delta _{i}B,z_{\mathcal{P}%
}(s_{i-1}))\Delta _{i}B+\beta _{i}  \label{e.7.21}
\end{align}
where 
\begin{equation}
\beta _{i}=\int_{s_{i-1}}^{s_{i}}(s_{i}-r)(\Omega _{u_{\mathcal{P}}(r)}(B_{%
\mathcal{P}}^{\prime }(r),z_{\mathcal{P}}(r))-\Omega _{u_{\mathcal{P}%
}(s_{i-1})}(B_{\mathcal{P}}^{\prime }(r),z_{\mathcal{P}}(s_{i-1})))B_{%
\mathcal{P}}^{\prime }(r)dr.  \label{e.7.22}
\end{equation}
By It\^{o}'s lemma, 
\begin{eqnarray*}
\Omega _{u_{\mathcal{P}}(s_{j-1})}(\Delta _{j}B,z_{\mathcal{P}%
}(s_{j-1}))\Delta _{j}B &=&\int_{s_{j-1}}^{s_{j}}\Omega _{u_{\mathcal{P}%
}(s_{j-1})}(B(r)-B(s_{j-1}),z_{\mathcal{P}}(s_{j-1}))dB(r) \\
&&+\int_{s_{j-1}}^{s_{j}}\Omega _{u_{\mathcal{P}}(s_{j-1})}(dB(r),z_{%
\mathcal{P}}(s_{j-1}))B(r)-B(s_{j-1}) \\
&&-\mathrm{Ric}_{u_{\mathcal{P}}(\underline{r})}z_{\mathcal{P}%
}(s_{j-1})\Delta _{j}s.
\end{eqnarray*}
Using this equation and the fact that $z_{\mathcal{P}}(0)=0,$ we may sum
equation (\ref{e.7.21}) on $i$ to find 
\begin{equation}
z_{\mathcal{P}}(s_{i})=\int_{0}^{s_{i}}\left( k^{\prime }(\underline{r}+)-%
\frac{1}{2}\mathrm{Ric}_{u_{\mathcal{P}}(\underline{r})}z_{\mathcal{P}%
}(r)\right) dr+M_{s_{i}}^{\mathcal{P}}+\sum_{j=1}^{i}\beta _{j},
\label{e.7.23}
\end{equation}
where $M^{\mathcal{P}}$ is the $\mathbb{R}^{d}$ -- valued Martingale, 
\begin{eqnarray*}
M_{s}^{\mathcal{P}} &:&=\int_{0}^{s}\Omega _{u_{\mathcal{P}}(\underline{r}%
)}(B(r)-B(\underline{r}),z_{\mathcal{P}}(\underline{r}))dB(r) \\
&&\qquad +\int_{0}^{s}\Omega _{u_{\mathcal{P}}(\underline{r})}(dB(r),z_{%
\mathcal{P}}(\underline{r}))\left( B(r)-B(\underline{r})\right) .
\end{eqnarray*}

Therefore $\delta _{i}=M_{s_{i}}^{\mathcal{P}}+\sum_{j=1}^{i}\beta _{j}.$

By the martingale moment inequality \cite[Prop. 3.26]{karatzas:shreve}, 
\begin{equation}
\mathbb{E}\left[ \sup_{s}|M_{s}^{\mathcal{P}}|^{p}\right] \leq C_{p}\mathbb{E%
}\left[ \langle M^{\mathcal{P}}\rangle _{1}^{p/2}\right]  \label{e.7.24}
\end{equation}
where $C_{p}$ is a constant and $\langle M^{\mathcal{P}}\rangle $ is the
quadratic variation of $M^{\mathcal{P}}.$ It is easy to estimate $\langle M^{%
\mathcal{P}}\rangle _{1}$ by 
\begin{equation*}
\langle M^{\mathcal{P}}\rangle _{1}\leq 2d\Lambda ^{2}\int_{0}^{1}\left|
B(r)-B(\underline{r})\right| ^{2}\left| z_{\mathcal{P}}(\underline{r}%
)\right| ^{2}dr
\end{equation*}
and hence by Jensen's inequality 
\begin{equation*}
\langle M^{\mathcal{P}}\rangle _{1}^{p/2}\leq \left( 2d\right) ^{p/2}\Lambda
^{p}\int_{0}^{1}\left| B(r)-B(\underline{r})\right| ^{p}\left| z_{\mathcal{P}%
}(\underline{r})\right| ^{p}dr.
\end{equation*}
Because \{$z_{\mathcal{P}}(\underline{r})\}_{r\in \lbrack 0,1]}$ is adapted
to the filtration generated by $B$ we may use the independence of the
increments of $B$ along with scaling to find 
\begin{eqnarray*}
\mathbb{E}\langle M^{\mathcal{P}}\rangle _{1}^{p/2} &\leq &\left( 2d\right)
^{p/2}\Lambda ^{p}\int_{0}^{1}\mathbb{E}\left| B(r)-B(\underline{r})\right|
^{p}\cdot \mathbb{E}\left| z_{\mathcal{P}}(\underline{r})\right| ^{p}dr \\
&=&C_{p}\left( 2d\right) ^{p/2}\Lambda ^{p}\int_{0}^{1}\left| r-\underline{r}%
\right| ^{p/2}\cdot \mathbb{E}\left| z_{\mathcal{P}}(\underline{r})\right|
^{p}dr \\
&\leq &C_{p}\left( 2d\right) ^{p/2}\Lambda ^{p}||k^{\prime }||_{1,\mathcal{P}%
}^{p}\int_{0}^{1}\left| r-\underline{r}\right| ^{p/2}\cdot \mathbb{E}e^{%
\frac{p}{2}\Lambda \sum_{j=1}^{n}|\Delta _{j}B|^{2}}dr,
\end{eqnarray*}
where equation (\ref{e.7.10}) was used in the last equality. By Lemma \ref
{l.8.5} of the Appendix, $\mathbb{E}e^{\frac{p}{2}\Lambda
\sum_{j=1}^{n}|\Delta _{j}B|^{2}}$ is bounded independent of $\mathcal{P}$
when $|\mathcal{P}|$ is sufficiently small. Hence we have shown that 
\begin{equation*}
\mathbb{E}\left[ \sup_{s}|M_{s}^{\mathcal{P}}|^{p}\right] \leq C_{p}(\Lambda
)||k^{\prime }||_{1,\mathcal{P}}^{p}\int_{0}^{1}\left| r-\underline{r}%
\right| ^{p/2}dr\leq C_{p}(\Lambda )||k^{\prime }||_{1,\mathcal{P}}^{p}|%
\mathcal{P}|^{p/2}.
\end{equation*}
So finish the proof it suffices to show that 
\begin{equation}
\mathbb{E}\left( \sum_{j=1}^{n}|\beta _{j}|\right) ^{p}\leq C||k^{\prime
}||_{1,\mathcal{P}}^{p}|\mathcal{P}|^{\gamma p}.  \label{e.7.25}
\end{equation}

By assumption, $u_{\mathcal{P}}$ solves the differential equation 
\begin{equation*}
u_{\mathcal{P}}^{\prime }(s)=\mathcal{H}_{u_{\mathcal{P}}(s)}u_{\mathcal{P}%
}(s)B_{\mathcal{P}}^{\prime }(s)\qquad u_{\mathcal{P}}(0)=u_{0},
\end{equation*}
so that for any $F\in C^{1}(O(M))$, $r\in (s_{i-1},s_{i}]$, 
\begin{equation}
|F(u_{\mathcal{P}}(r))-F(u_{\mathcal{P}}(s_{i-1}))|\leq
C|\int_{s_{i-1}}^{r}B_{\mathcal{P}}^{\prime }(s)ds|\leq C|\Delta _{i}B|,
\label{e.7.26}
\end{equation}
where $C$ bounds the horizontal derivatives of $F$. Applying this estimate
to $\Omega $ implies 
\begin{equation}
|\Omega _{u_{\mathcal{P}}(r)}-\Omega _{u_{\mathcal{P}}(s_{i-1})}|\leq
\Lambda |\Delta _{i}B|.  \label{e.7.27}
\end{equation}
Using the inequalities in (\ref{e.7.12}) and (\ref{e.7.27}) and equation (%
\ref{e.7.22}) we find that 
\begin{eqnarray}
|\beta _{i}| &\leq &\frac{1}{2}\Lambda \max_{s_{i-1}\leq s\leq s_{i}}|z_{%
\mathcal{P}}(s)-z_{\mathcal{P}}(s_{i-1})||\Delta _{i}B|^{2}+\Lambda |z_{%
\mathcal{P}}(s_{i-1})||\Delta _{i}B|^{3}  \notag \\
&\leq &\frac{1}{2}\Lambda \left( |k^{\prime }(s_{i-1}+)|\Delta _{i}s+\frac{1%
}{2}|z_{\mathcal{P}}(s_{i-1})|\Lambda \left| \Delta _{i}B\right| ^{2}\right)
\cosh \left( \sqrt{\Lambda }\left| \Delta _{i}B\right| \right) |\Delta
_{i}B|^{2}  \notag \\
&&+\Lambda |z_{\mathcal{P}}(s_{i-1})||\Delta _{i}B|^{3}.  \label{e.7.28}
\end{eqnarray}
Letting $K_{\gamma }$ denote the random variable defined in equation (\ref
{e.8.15}) of Fernique's Lemma \ref{l.8.3}, the above estimate implies that 
\begin{eqnarray*}
|\beta _{i}| &\leq &\frac{\Lambda }{2}|k^{\prime }(s_{i-1}+)|\Delta
_{i}s\cosh \left( \sqrt{\Lambda }K_{\gamma }|\mathcal{P}|\right) K_{\gamma
}^{2}|\mathcal{P}|^{2\gamma } \\
&&+\left( \frac{\Lambda ^{2}}{4}K_{\gamma }^{4}\left| \Delta _{i}s\right|
^{4\gamma }\cosh \left( \sqrt{\Lambda }K_{\gamma }|\mathcal{P}|\right)
+C_{2}K_{\gamma }^{3}\left| \Delta _{i}s\right| ^{3\gamma }\right) |z_{%
\mathcal{P}}(s_{i-1})|
\end{eqnarray*}
where $\gamma \in (0,1/2).$ We will now suppose that $\gamma $ is close to $%
1/2.$ Then by equation (\ref{e.7.10}) of Lemma \ref{l.7.10}, we find that 
\begin{eqnarray*}
\sum_{i=1}^{n}|\beta _{i}| &\leq &\frac{\Lambda }{2}||k^{\prime }||_{1,%
\mathcal{P}}\cosh \left( \sqrt{\Lambda }K_{\gamma }|\mathcal{P}|\right)
K_{\gamma }^{2}|\mathcal{P}|^{2\gamma } \\
&&+C||k^{\prime }||_{1,\mathcal{P}}|\mathcal{P}|^{3\gamma -1}\left(
K_{\gamma }^{4}\cosh \left( \sqrt{\Lambda }K_{\gamma }|\mathcal{P}|\right)
+K_{\gamma }^{3}\right) e^{\frac{1}{2}\Lambda \sum_{j=1}^{n}|\Delta
_{j}B|^{2}} \\
&\leq &C||k^{\prime }||_{1,\mathcal{P}}|\mathcal{P}|^{3\gamma -1}\left(
\left( K_{\gamma }^{4}+K_{\gamma }^{2}\right) \cosh \left( \sqrt{\Lambda }%
K_{\gamma }|\mathcal{P}|\right) +K_{\gamma }^{3}\right) e^{\frac{1}{2}%
\Lambda \sum_{j=1}^{n}|\Delta _{j}B|^{2}}.
\end{eqnarray*}
Using Lemma \ref{l.8.4} and \ref{l.8.5} of the Appendix, it follows that 
\begin{equation*}
\left( \left( K_{\gamma }^{4}+K_{\gamma }^{2}\right) \cosh \left( \sqrt{%
\Lambda }K_{\gamma }|\mathcal{P}|\right) +K_{\gamma }^{3}\right) e^{\frac{1}{%
2}\Lambda \sum_{j=1}^{n}|\Delta _{j}B|^{2}}
\end{equation*}
is bounded in all $L^{p}$ for $\left| \mathcal{P}\right| $ small. This
proves $\mathbb{E}\left( \sum_{j=1}^{n}|\beta _{j}|\right) ^{p}\leq
C||k^{\prime }||_{1,\mathcal{P}}^{p}|\mathcal{P}|^{(3\gamma -1)p}$ which
proves equation (\ref{e.7.25}) since $(3\gamma -1)$ approaches $1/2$ when $%
\gamma $ approaches $1/2$.
\end{proof}

\begin{lemma}
\label{l.7.13} Let $\epsilon _{\mathcal{P}}$ be defined by 
\begin{equation}
\epsilon _{\mathcal{P}}(s):=z_{\mathcal{P}}(s)+\int_{0}^{s}\left( \frac{1}{2}%
\mathrm{Ric}_{u_{\mathcal{P}}(r)}z_{\mathcal{P}}(r)-k^{\prime }(r)\right) dr.
\label{e.7.29}
\end{equation}
Then for all $\gamma \in (0,\frac{1}{2})$ and $p\in \lbrack 1,\infty ),$%
\begin{equation}
\mathbb{E}\left[ \max_{s}|\epsilon _{\mathcal{P}}(s)|^{p}\right] \leq
C(\gamma ,p)\left( ||k^{\prime }||_{1,\mathcal{P}}^{p}|\mathcal{P}|^{\gamma
p}+|||k^{\prime }|||_{\mathcal{P}}^{p}\right) .  \label{e.7.30}
\end{equation}
\end{lemma}

\begin{proof}
Let $\delta _{i}$ be as in the previous lemma and set $\delta _{\mathcal{P}%
}(s):=\sum_{i=1}^{n}\delta _{i}1_{(s_{i-1},s_{i}]}(s).$ By the definitions
of $\epsilon _{\mathcal{P}}$, (\ref{e.7.29}) and $\delta _{\mathcal{P}}$, we
have for $s\in (s_{i-1},s_{i}]$, 
\begin{align*}
\epsilon _{\mathcal{P}}(s)-\delta _{\mathcal{P}}(s)& =z_{\mathcal{P}}(s)-z_{%
\mathcal{P}}(s_{i}) \\
& \quad +\frac{1}{2}\left( \int_{0}^{s}\mathrm{Ric}_{u_{\mathcal{P}}(r)}z_{%
\mathcal{P}}(r)dr-\int_{0}^{s_{i}}\mathrm{Ric}_{u_{\mathcal{P}}(\underline{r}%
)}z_{\mathcal{P}}(\underline{r})dr\right) \\
& \quad +\int_{0}^{s_{i}}k^{\prime }(\underline{r}+)dr-\int_{0}^{s}k^{\prime
}(r)dr \\
& =\frac{1}{2}\int_{0}^{s_{i}}\left( \mathrm{Ric}_{u_{\mathcal{P}}(r)}z_{%
\mathcal{P}}(r)-\mathrm{Ric}_{u_{\mathcal{P}}(\underline{r})}z_{\mathcal{P}}(%
\underline{r})\right) dr \\
& \quad \,-\frac{1}{2}\int_{s}^{s_{i}}\mathrm{Ric}_{u_{\mathcal{P}}(r)}z_{%
\mathcal{P}}(r)dr \\
& \quad +(z_{\mathcal{P}}(s)-z_{\mathcal{P}}(s_{i}))+(k(s_{i})-k(s)) \\
& \quad -\int_{0}^{s_{i}}(k^{\prime }(r)-k^{\prime }(\underline{r}+))dr \\
& =:\frac{1}{2}A_{i}+\frac{1}{2}B_{i}+C_{i}(s)+E_{i},
\end{align*}
where for $r\in (s_{j-1},s_{j}],$ $\underline{r}:=s_{j-1}.$ We will now
prove the estimate 
\begin{equation*}
\mathbb{E}\left[ \sup_{s}|\epsilon _{\mathcal{P}}(s)-\delta _{\mathcal{P}%
}(s)|^{p}\right] \leq C(\gamma ,p)\left( ||k^{\prime }||_{1,\mathcal{P}}^{p}|%
\mathcal{P}|^{2\gamma p}+|||k^{\prime }|||_{\mathcal{P}}^{p}\right)
\end{equation*}
This will complete the proof (\ref{e.7.30}) in view Lemma \ref{l.7.12}.

By definition of $|||k^{\prime }|||_{\mathcal{P}}$ in equation (\ref{e.7.9}) 
\begin{equation}
\max_{i}|E_{i}|\leq |||k^{\prime }|||_{\mathcal{P}}.  \label{e.7.31}
\end{equation}
In the argument to follow let $\{K_{\mathcal{P}}\}_{\mathcal{P}}$ denote a
collection functions on $(\mathrm{W}(\mathbb{R}^{d}),\mu )$ such that $\sup_{%
\mathcal{P}}\Vert K_{\mathcal{P}}\Vert _{L^{p}(\mu )}<\infty $ for all $p\in
\lbrack 1,\infty ).$ Using equation (\ref{e.7.10}) with $b=B$ and $\sigma
=\phi (B_{\mathcal{P}})$ and Lemma \ref{l.8.5} of the Appendix, 
\begin{equation*}
|B_{i}|\leq \Vert \mathrm{Ric}\Vert _{\infty }|\mathcal{P}|\Vert z_{\mathcal{%
P}}\Vert _{\infty }\leq \Vert \mathrm{Ric}\Vert _{\infty }K_{\mathcal{P}%
}|\Vert k^{\prime }\Vert _{1,\mathcal{P}}|\mathcal{P}|.
\end{equation*}
So for $p\in \lbrack 1,\infty )$, 
\begin{equation*}
\mathbb{E}\left[ \max_{i}|B_{i}|^{p}\right] \leq \Vert \mathrm{Ric}\Vert
_{\infty }^{p}|\Vert k^{\prime }\Vert _{1,\mathcal{P}}^{p}\mathbb{E}\left[
K_{\mathcal{P}}^{p}\right] |\mathcal{P}|^{p}\leq C|\mathcal{P}|^{p}.
\end{equation*}
Next we consider $C_{i}$. We have $C_{i}(s_{i})=0$ and by (\ref{e.7.13}) and
(\ref{e.7.14}) with $b=B$ and $\sigma =\phi (B_{\mathcal{P}})$ for $s\in
(s_{i-1},s_{i}]$, 
\begin{align*}
C_{i}^{\prime }(s)& =z_{\mathcal{P}}^{\prime }(s)-k^{\prime }(s) \\
& =k^{\prime }(s_{i-1})-k^{\prime }(s)+\int_{s_{i-1}}^{s}\Omega _{u_{%
\mathcal{P}}(r)}(B^{\prime }(r),z_{\mathcal{P}}(r))drB^{\prime }(s)
\end{align*}
which implies after integrating 
\begin{equation*}
|C_{i}(s)|\leq \Lambda |\Delta _{i}B|^{2}||z_{\mathcal{P}}||_{\infty
}+|||k^{\prime }|||_{\mathcal{P}}\leq \Lambda K_{\gamma }^{2}|\Delta
_{i}s|^{2\gamma }||z_{\mathcal{P}}||_{\infty }+|||k^{\prime }|||_{\mathcal{P}%
}
\end{equation*}
where $\Lambda $ is a bound on $\Omega $ and $K_{\gamma }$ is defined in
Lemma \ref{l.8.4}. Therefore, again by (\ref{e.7.10}) and Lemma \ref{l.8.5},
if $p\in \lbrack 1,\infty )$ and $\gamma \in (0,1/2)$ then 
\begin{equation*}
\mathbb{E}\left[ \max_{i,s}|C_{i}(s)|^{p}\right] \leq C(\gamma ,p,\Lambda
)\left( ||k^{\prime }||_{1,\mathcal{P}}^{p}|\mathcal{P}|^{2\gamma
p}+|||k^{\prime }|||_{\mathcal{P}}^{p}\right) .
\end{equation*}
%By equation \eqref{e.7.29a}, 
%for $p\in [1,\infty )$ and $\gamma \in (0,1/2)$,
%\begin{align*}
%\EE\left[ \max_{i,s}|C_{i}(s)|^{p}\right] &\leq 
%C \left ( \max_i |k'(s_{i-1} +)| (\Delta_i s)^
%+ C(\gamma,p) || k'||_{1,\mathcal{P}}  |\mathcal{P}|^{2\gamma} \right )^p  \\
%&\leq
%C(\gamma,p) \left ( || k' ||_{\infty} | \mathcal{P}| 
%+ || k'||_{1,\mathcal{P}} |\mathcal{P}|^{2\gamma} \right )^p   
%%\\
%%&\leq C(\gamma,p) \left ( || k'||_{\infty,\mathcal{P}} |\mathcal{P}|^{1-\gamma} 
%%+ || k' ||_{1,\mathcal{P}} \right )^p |\mathcal{P}|^{\gamma p}
%\end{align*}
\end{proof}

%$$
%\Vert k^{\prime }\Vert_{\infty }^{p}
% \EE(1+K_{\gamma }^{2}K_{\mathcal{P}})^{p}|\mathcal{P}|^{\gamma p}\rightarrow
%0\text{ as }|\mathcal{P}|\rightarrow 0. 
%$$
So to finish the proof it only remains to consider the $A_{i}$ term.
Applying the estimate in equation (\ref{e.7.26}) with $F=\mathrm{Ric}$
gives, for $r\in (s_{j-1},s_{j}]$, 
\begin{equation*}
|\mathrm{Ric}_{u_{\mathcal{P}}(r)}-\mathrm{Ric}_{u_{\mathcal{P}}(\underline{r%
})}|\leq C|\Delta _{j}B|\leq CK_{\gamma }|\mathcal{P}|^{\gamma }.
\end{equation*}
where $C$ is a bound on the horizontal derivative of $\mathrm{Ric}$.
Therefore, 
\begin{eqnarray*}
|A_{i}| &\leq &CK_{\gamma }|\mathcal{P}|^{\gamma }\Vert z_{\mathcal{P}}\Vert
_{\infty }+\Vert \mathrm{Ric}\Vert _{\infty }\int_{0}^{1}|z_{\mathcal{P}%
}(r)-z_{\mathcal{P}}(\underline{r})|dr \\
&\leq &CK_{\gamma }|\mathcal{P}|^{\gamma }||k^{\prime }||_{1,\mathcal{P}}e^{%
\frac{1}{2}\Lambda \sum_{j=1}^{n}|\Delta _{j}B|^{2}} \\
&&+\Vert \mathrm{Ric}\Vert _{\infty }\left( ||k^{\prime }||_{1,\mathcal{P}%
}\left| \mathcal{P}\right| +\frac{1}{2}\Lambda \max_{i}|\Delta
_{i}B|^{2}\Vert z_{\mathcal{P}}\Vert _{\infty }\right) \\
&\leq &C||k^{\prime }||_{1,\mathcal{P}}\left\{ e^{\frac{1}{2}\Lambda
\sum_{j=1}^{n}|\Delta _{j}B|^{2}}\left( K_{\gamma }|\mathcal{P}|^{\gamma
}+K_{\gamma }^{2}\left| \mathcal{P}\right| ^{2\gamma }\right) +\left| 
\mathcal{P}\right| \right\} \\
&\leq &K_{\mathcal{P}}||k^{\prime }||_{1,\mathcal{P}}|\mathcal{P}|^{\gamma },
\end{eqnarray*}
wherein we have made use of equations (\ref{e.7.10}) and (\ref{e.7.12}) of
Lemma \ref{l.7.10} in the second inequality, equation (\ref{e.7.10}) and the
definition of $K_{\gamma }$ in equation (\ref{e.8.15}) in the third
inequality, and Lemmas \ref{l.8.4} and \ref{l.8.5} in the last inequality.
Thus 
\begin{equation*}
\mathbb{E}\left[ \max_{i}|A_{i}|^{p}\right] \leq C(\gamma ,p)||k^{\prime
}||_{1,\mathcal{P}}^{p}|\mathcal{P}|^{\gamma p}
\end{equation*}
for $p\in \lbrack 1,\infty )$ and $\gamma \in (0,1/2).$ %Using 
%From \eqref{e.7.29a} we get 
%$$
%\int_{s_{i-1}}^{s_i} |z_{\mathcal{P}}(r) - z_{\mathcal{P}}(s_{i-1})| dr 
%\leq |k'(s_{i-1}| (\Delta_i s)^2 + 
%%&\leq
%%%\max_i |k'(s_{i-1} +)| (\Delta_i s) 
%C (  || k' ||_{1,\mathcal{P}} K_{\mathcal{P}} K_{\gamma}^2 |\mathcal{P}|^{2\gamma}
% + || k'||_{\infty} |\mathcal{P}| +  || k'||_{1,\mathcal{P}} K_{\mathcal{P}}
%K_{\gamma}^2 |\mathcal{P}|^{2\gamma} ) 
%+ C(\gamma,p)  \left ( || k' ||_{\infty} | \mathcal{P}| 
%+ || k'||_{1,\mathcal{P}} |\mathcal{P}|^{2\gamma}\right )  
%\end{align*}
%&\le C \left (  
%\Vert k^{\prime }\Vert_{1,\mathcal{P} } K_{\mathcal{P}} K_{\gamma}
%+\Vert {\text{\rm Ric}}\Vert _{\infty }  (|| k'||_{\infty,\mathcal{P}} |\mathcal{P}|^{1-\gamma} + 
%|| k'||_{1,\mathcal{P}} K_{\mathcal{P}} K_{\gamma}^2|\mathcal{P}|^{\gamma} \right ) |\mathcal{P}|^{\gamma }
%\end{align*}
%where we have made use of equations \eqref{e.7.29a} and
%\eqref{e.7.28a}. This shows that 
%$$
%\EE\left[ \max_i|A_i|^{p}\right] \le 
%C(\gamma,p) \left ( || k' ||_{\infty } |\mathcal{P}| + 
%|| k' ||_{1,\mathcal{P}} |\mathcal{P}|^{2\gamma} \right )^p 
%\rightarrow 0\text{ as }|\mathcal{P}|\rightarrow 0 
%$$
%for $p\in [1,\infty )$ and $\gamma \in (0,1/2).$ 
This completes the proof of Lemma \ref{l.7.12}.\qed

\begin{proof}
Let $\epsilon _{\mathcal{P}}$ be defined as in equation (\ref{e.7.29}) and
let $y_{\mathcal{P}}(s)$ denote the solution to the differential equation, 
\begin{equation*}
y_{\mathcal{P}}^{\prime }(s)+\frac{1}{2}\mathrm{Ric}_{u_{\mathcal{P}}(s)}y_{%
\mathcal{P}}(s)=k^{\prime }(s)\text{\textrm{\ with }}y_{\mathcal{P}}(0)=0.
\end{equation*}
Then 
\begin{equation*}
z_{\mathcal{P}}(s)-y_{\mathcal{P}}(s)=-\int_{0}^{s}\frac{1}{2}\mathrm{Ric}%
_{u_{\mathcal{P}}(r)}(z_{\mathcal{P}}(r)-y_{\mathcal{P}}(r))dr+\epsilon _{%
\mathcal{P}}(s)
\end{equation*}
and hence 
\begin{equation*}
|z_{\mathcal{P}}(s)-y_{\mathcal{P}}(s)|\leq \int_{0}^{s}C|(z_{\mathcal{P}%
}(r)-y_{\mathcal{P}}(r)|dr+\epsilon _{\mathcal{P}}(s),
\end{equation*}
where $C$ is a bound on $\frac{1}{2}\mathrm{Ric}.$ So by Gronwall's
inequality, 
\begin{equation*}
|z_{\mathcal{P}}(s)-y_{\mathcal{P}}(s)|\leq \max_{s}\left( |\epsilon _{%
\mathcal{P}}(s)|e^{Cs}\right) \leq \max_{s}|\epsilon _{\mathcal{P}}(s)|e^{C},
\end{equation*}
which combined with equation (\ref{e.7.30}) of Lemma \ref{l.7.12} shows that 
\begin{equation*}
\mathbb{E}\left[ \max_{s}|z_{\mathcal{P}}(s)-y_{\mathcal{P}}(s)|^{p}\right]
\leq C(\gamma ,p)\left( ||k^{\prime }||_{1,\mathcal{P}}^{p}|\mathcal{P}%
|^{\gamma p}+|||k^{\prime }|||_{\mathcal{P}}^{p}\right) .
\end{equation*}
for $p\in \lbrack 1,\infty )$, $\gamma \in (0,1/2)$.

To finish the proof of the theorem it is sufficient to prove 
\begin{equation}
\mathbb{E}\left[ \max_{s}|y_{\mathcal{P}}(s)-z(s)|^{p}\right] \leq C(\gamma
,p)\left( ||k^{\prime }||_{1,\mathcal{P}}^{p}|\mathcal{P}|^{\gamma
p}+|||k^{\prime }|||_{\mathcal{P}}^{p}|\mathcal{P}|^{\gamma p}\right) .
\label{e.7.32}
\end{equation}
First note that a Gronwall estimate gives 
\begin{equation}
\max_{s}|z(s)|\leq ||k^{\prime }||_{L^{1}(ds)}e^{||\mathrm{Ric}||_{\infty
}s}\leq C||k^{\prime }||_{L^{1}(ds)}  \label{e.7.33}
\end{equation}
and similarly 
\begin{equation*}
\max_{s}|y_{\mathcal{P}}(s)|\leq C||k^{\prime }||_{L^{1}(ds)},
\end{equation*}
where $||k^{\prime }||_{L^{1}(ds)}=\int_{0}^{1}|k^{\prime }(s)|ds.$ Let $%
w=y_{\mathcal{P}}-z$. Then 
\begin{equation*}
w^{\prime }(s)=\frac{1}{2}\mathrm{Ric}_{u_{\mathcal{P}}(s)}w(s)+\frac{1}{2}%
\left( \mathrm{Ric}_{u_{\mathcal{P}}(s)}-\mathrm{Ric}_{u(s)}\right) z(s)
\end{equation*}
Letting 
\begin{equation*}
A_{\mathcal{P}}=\max_{s}\frac{1}{2}|\mathrm{Ric}_{u_{\mathcal{P}}(s)}-%
\mathrm{Ric}_{u(s)}|
\end{equation*}
the inequality (\ref{e.7.33}) and an application of Gronwall's inequality
gives 
\begin{equation}
|w(s)|\leq CA_{\mathcal{P}}||k^{\prime }||_{L^{1}}e^{Cs}  \label{e.7.34}
\end{equation}
Theorem \ref{t.4.14} implies 
\begin{equation*}
\mathbb{E}\left[ |A_{\mathcal{P}}|^{p}\right] \leq C(\gamma ,p)|\mathcal{P}%
|^{\gamma p}
\end{equation*}
and hence by (\ref{e.7.34}), 
\begin{equation*}
\mathbb{E}\left[ \max_{s}|y_{\mathcal{P}}(s)-z(s)|^{p}\right] \leq C(\gamma
,p)||k^{\prime }||_{L^{1}}^{p}|\mathcal{P}|^{\gamma p}.
\end{equation*}
This implies (\ref{e.7.32}) in view of the fact that 
\begin{equation*}
||k^{\prime }||_{L^{1}}\leq ||k^{\prime }||_{1,\mathcal{P}}+|||k^{\prime
}|||_{\mathcal{P}}.
\end{equation*}
This completes the proof of Theorem \ref{t.7.11}.
\end{proof}

\subsection{Integration by Parts for Wiener Measure\label{s.7.3}}

\begin{prop}
\label{p.7.14}Let $|\mathcal{P}|:=\max \{|\Delta _{i}s|:i=1,2,\ldots ,n\}$
denote the mesh size of the partition $\mathcal{P}$ and $f$ be a function on 
$\mathrm{H}(M)$ and $\tilde{f}$ on $\mathrm{W}(M)$ as in Theorem \ref{t.4.17}%
. Then 
\begin{equation}
\lim_{|\mathcal{P}|\rightarrow 0}\int_{H_{\mathcal{P}}(M)}f\left(
\sum_{i=1}^{n}\langle k^{\prime }(s_{i-1}+),\Delta _{i}b\rangle \right) \nu
_{\mathcal{P}}^{1}=\int_{\mathrm{W}(M)}\left( \tilde{f}\int_{0}^{1}\langle
k^{\prime },d\tilde{b}\rangle \right) d\nu ,  \label{e.7.35}
\end{equation}
where $\Delta _{i}b$ is to be interpreted as a function on $\mathrm{H}(M)$
as in equation (\ref{e.7.7}) and $\tilde{b}$ is the anti-development map.
Recall that $\tilde{b}$ is an $\mathbb{R}^{d}$--valued Brownian motion on $(%
\mathrm{W}(M),\nu )$ which was defined in Definition \ref{d.4.15}. Here $%
\int_{0}^{1}\langle k^{\prime },d\tilde{b}\rangle $ denotes the It\^{o}
integral of $k^{\prime }$ relative to $\tilde{b}.$
\end{prop}

\begin{proof}
Let $B$ denote the standard $\mathbb{R}^{d}$--valued Brownian motion in
Notation \ref{n.1.2} and $u$ denote the solution to the Stratonovich
stochastic differential equation \ref{e.7.17}. By Lemma \ref{l.4.11} and
Theorem \ref{t.4.10}, 
\begin{equation}
\int_{\mathrm{H}_{\mathcal{P}}(M)}f\left( \sum_{i=1}^{n}\langle k^{\prime
}(s_{i-1}+),\Delta _{i}b\rangle \right) \nu _{\mathcal{P}}^{1}=\mathbb{E}%
\left[ f(\phi (B_{\mathcal{P}}))\left( \sum_{i=1}^{n}\langle k^{\prime
}(s_{i-1}+),\Delta _{i}B\rangle \right) \right] .  \label{e.7.36}
\end{equation}
By the isometry property of the It\^{o} integral, we find that 
\begin{equation*}
\lim_{|\mathcal{P}|\rightarrow 0}\left( \sum_{i=1}^{n}\langle k^{\prime
}(s_{i-1}+),\Delta _{i}B\rangle \right) =\int_{0}^{1}\langle k^{\prime
},dB\rangle ,
\end{equation*}
where the convergence takes place in $L^{2}(\mathrm{W}(\mathbb{R}^{d}),\mu
). $ As in the proof of theorem \ref{t.4.17}, $f(\phi (B_{\mathcal{P}}))$
converges to $F(u)$ in $L^{2}$ as well. Therefore we may pass to the limit
in equation (\ref{e.7.36}) to conclude that 
\begin{equation*}
\lim_{|\mathcal{P}|\rightarrow 0}\int_{\mathrm{H}_{\mathcal{P}}(M)}f\left(
\sum_{i=1}^{n}\langle k^{\prime }(s_{i-1}+),\Delta _{i}b\rangle \right) \nu
_{\mathcal{P}}^{1}=\mathbb{E}\left[ F(u)\int_{0}^{1}\langle k^{\prime
},dB\rangle \right] .
\end{equation*}
Since $(B,u)$ and $(\tilde{b},/\tilde{/})$ have the same distribution, 
\begin{equation*}
\mathbb{E}\left[ F(u)\int_{0}^{1}\langle k^{\prime },dB\rangle \right]
=\int_{\mathrm{W}(M)}\left( \tilde{f}\int_{0}^{1}\langle k^{\prime },d\tilde{%
b}\rangle \right) d\nu .
\end{equation*}
The previous two displayed equations prove equation (\ref{e.7.35}).
\end{proof}

\begin{definition}
\label{d.7.15}A function $f:\mathrm{W}(M)\rightarrow \mathbb{R}$ is said to
be a \textbf{smooth cylinder function} if $f$ is of the form 
\begin{equation}
f(\sigma )=F\circ \pi _{\mathcal{P}}(\sigma )=F(\sigma _{\mathcal{P}})
\label{e.7.37}
\end{equation}
for some partition $\mathcal{P}$ and some $F\in C^{\infty }(M^{\mathcal{P}})$%
.\qed
\end{definition}

We are now prepared for the main theorem of this section.

\begin{thm}
\label{t.7.16}Let $k\in PC^{1},$ $z$ be the solution to the differential
equation (\ref{e.7.18}) of Theorem \ref{t.7.11} and $f$ be a cylinder
function on $\mathrm{W}(M).$ Then 
\begin{equation}
\int_{\mathrm{W}(M)}X^{z}f\,d\nu =\int_{\mathrm{W}(M)}f\left(
\int_{0}^{1}\langle k^{\prime },d\tilde{b}\rangle \right) \,d\nu ,
\label{e.7.38}
\end{equation}
where 
\begin{equation*}
(X^{z}f)(\sigma ):=\sum_{i=1}^{n}\langle \nabla _{i}f)(\sigma
),X_{s_{i}}^{z}(\sigma )\rangle =\sum_{i=1}^{n}\langle \nabla _{i}f)(\sigma
),/\tilde{/}_{s_{i}}(\sigma )z(s_{i},\sigma )\rangle
\end{equation*}
and $(\nabla _{i}f)(\sigma )$ denotes the gradient $F$ in the i$^{th}$
variable evaluated at $(\sigma (s_{1}),\sigma (s_{2}),\ldots ,\sigma
(s_{n})).$
\end{thm}

\begin{proof}
The proof is easily completed by passing to the limit $|\mathcal{P}%
|\rightarrow 0$ in equation (\ref{e.7.6}) of Corollary \ref{c.7.7} making
use of Proposition \ref{p.7.14}, Theorems \ref{t.7.11}, \ref{t.4.14}, and
Corollary \ref{c.4.13}
\end{proof}

\section{Appendix: Basic Estimates}

\label{s.8}

\subsection{Determinant Estimates}

\label{s.8.1}

\begin{lemma}
\label{l.8.1} Let $U$ be a $d\times d$ matrix such that $|U|<1,$ then 
\begin{equation}
\det (I-U)=\exp \left( -{\text{tr}} U+\Psi (U)\right) ,  \label{e.8.1}
\end{equation}
where $\Psi (U):=-\sum_{n=2}^{\infty }\frac{1}{n}{\text{tr}} U^{n}$.
Moreover, $\Psi (U)$ satisfies the bound, 
\begin{equation}
|\Psi (U)|\leq \sum_{n=2}^{\infty }\frac{d}{n}|U|^{n}\leq
d|U|^{2}(1-|U|)^{-1}.  \label{e.8.2}
\end{equation}
\end{lemma}

\begin{proof}
Equation (\ref{e.8.1}) is just a rewriting of the standard formula: 
\begin{equation*}
\log (\det (I-U))=-\sum_{n=0}^{\infty }\frac{1}{n+1}{\text{tr}}\left(
U^{n+1}\right) ,
\end{equation*}
which is easily deduced by integrating the following identity, 
\begin{eqnarray*}
\frac{d}{ds}\log (\det (I-sU)) &=&-{\text{tr}}((I-sU)^{-1}U) \\
&=&-{\text{tr}}\left( \sum_{n=0}^{\infty }s^{n}U^{n}U\right)
=-\sum_{n=0}^{\infty }s^{n}{\text{tr}}\left( U^{n+1}\right) .
\end{eqnarray*}
Since for any $d\times d$ matrix $|{\text{tr}} U|\leq d|U|$ and $|U^{k}|\leq
|U|^{k}$, it follows that 
\begin{equation*}
|{\text{tr}}(U^{k})|\leq d|U|^{k}
\end{equation*}
and hence 
\begin{equation*}
|\Psi (U)|\leq \sum_{n=2}^{\infty }\frac{d}{n}|U|^{n}\leq
d|U|^{2}(1-|U|)^{-1}.
\end{equation*}
\end{proof}

\subsection{Ordinary Differential Equation Estimates}

\label{s.8.2}

\begin{lemma}
\label{l.8.2} Let $A(s)$ be a $d\times d$ matrix for all $s\in \lbrack 0,1]$
and let $Z(s)$ be either a $\mathbb{R}^{d}$ valued or $d\times d$ matrix
valued solution to the second order differential equation 
\begin{equation}
Z^{\prime \prime }(s)=A(s)Z(s).  \label{e.8.3}
\end{equation}
Then 
\begin{equation}
|Z(s)-Z(0)|\leq |Z(0)|\left( \cosh \sqrt{K}s-1\right) +|Z^{\prime }(0)|\frac{%
\sinh \sqrt{K}s}{\sqrt{K}}  \label{e.8.4}
\end{equation}
and 
\begin{equation}
|Z(s)-Z(0)|\leq s|Z^{\prime }(0)|+K\frac{s^{2}}{2}Z^{\ast }(s)  \label{e.8.5}
\end{equation}
where $Z^{\ast }(s):=\max \{Z(r)|:0\leq r\leq s\},$ $K:=\sup_{s\in \lbrack
0,1]}|A(s)|$ and $|A|$ denotes the operator norm of $A.$
\end{lemma}

\begin{proof}
By the Taylor's theorem with integral remainder, 
\begin{eqnarray}
Z(s) &=&Z(0)+sZ^{\prime }(0)+\int_{0}^{s}Z^{\prime \prime }(u)(s-u)du  \notag
\\
&=&Z(0)+sZ^{\prime }(0)+\int_{0}^{s}A(u)Z(u)(s-u)du  \label{e.8.6}
\end{eqnarray}
and therefore 
\begin{eqnarray}
|Z(s)-Z(0)| &\leq &s|Z^{\prime }(0)|+K\int_{0}^{s}|Z(u)|(s-u)du  \notag \\
&\leq &s|Z^{\prime }(0)|+K\int_{0}^{s}|Z(u)-Z(0)|(s-u)du+\frac{s^{2}}{2}%
K|Z(0)|=:f(s).  \label{e.8.7}
\end{eqnarray}
One may easily deduce equation (\ref{e.8.5}) from the first inequality in
this equation.

Note that $f(0)=0,$ 
\begin{equation*}
f^{\prime }(s)=|Z^{\prime }(0)|+K\int_{0}^{s}|Z(u)-Z(0)|(s-u)du+sK|Z(0)|
\end{equation*}
$f^{\prime }(0)=|Z^{\prime }(0)|$ and 
\begin{equation*}
f^{\prime \prime }(s)=K|Z(s)-Z(0)|+K|Z(0)|\leq Kf(s)+K|Z(0)|.
\end{equation*}
That is: 
\begin{equation}
f^{\prime \prime }(s)=Kf(s)+\eta (s),\quad f(0)=0,\quad \text{and \quad }%
f^{\prime }(0)=|Z^{\prime }(0)|,  \label{e.8.8}
\end{equation}
where $\eta (s):=f^{\prime \prime }(s)-Kf(s)\leq K|Z(0)|.$ Equation (\ref
{e.8.8}) may be solved by variation of parameters to find: 
\begin{eqnarray*}
f(s) &=&|Z^{\prime }(0)|\frac{\sinh \sqrt{K}s}{\sqrt{K}}+\int_{0}^{s}\frac{%
\sinh \sqrt{K(}s-r)}{\sqrt{K}}\eta (r)dr \\
&\leq &|Z^{\prime }(0)|\frac{\sinh \sqrt{K}s}{\sqrt{K}}+|Z(0)|\int_{0}^{s}%
\sqrt{K}\sinh \sqrt{K(}s-r)dr \\
&=&|Z^{\prime }(0)|\frac{\sinh \sqrt{K}s}{\sqrt{K}}+|Z(0)|\cosh (\sqrt{K}%
s-1).
\end{eqnarray*}
Combining this equation with equation (\ref{e.8.7}) proves equation (\ref
{e.8.4}).
\end{proof}

\begin{lemma}
\label{l.8.3} Suppose that $Z$ is a $d\times d$ -- matrix valued solution to
equation (\ref{e.8.3}) with $Z(0)=0$ and $Z^{\prime }(0)=I.$ Let $K>0$, $%
K_{1}>0$ be constants so that $\sup_{s\in \lbrack 0,1]}|A(s)|\leq K$ and $%
\sup_{s\in \lbrack 0,1]}|A^{\prime }(s)|\leq K_{1}$. Then 
\begin{equation}
Z(s)=sI+\frac{s^{3}}{6}A(0)+s\mathcal{E}(s),  \label{e.8.9}
\end{equation}
where 
\begin{equation}
|\mathcal{E}(s)|=\frac{1}{6}(2K_{1}s^{3}+\frac{1}{2}K^{2}s^{4})\cosh (\sqrt{K%
}s)  \label{e.8.10}
\end{equation}
\end{lemma}

\begin{proof}
Using the definition of $Z$ in equation (\ref{e.8.3}) we have that $%
Z(0)=Z^{\prime \prime }(0)=0,$ $Z^{\prime }(0)=0,$ 
\begin{equation*}
Z^{(3)}(s):=\frac{d^{3}}{ds^{3}}Z(s)=A^{\prime }(s)Z(s)+A(s)Z^{\prime }(s),
\end{equation*}
and hence $Z^{(3)}(0)=A(0).$ By Taylor's theorem with integral remainder 
\begin{eqnarray*}
Z(s) &=&sI+\frac{1}{2}\int_{0}^{s}Z^{(3)}(\xi )(s-\xi )^{2}d\xi \\
&=&sI+\frac{s^{3}}{6}A(0)+\frac{1}{2}\int_{0}^{s}\left( Z^{(3)}(\xi
)-A(0)\right) (s-\xi )^{2}d\xi .
\end{eqnarray*}
Now using Lemma \ref{l.8.2} with $Z(0)=0,$we find 
\begin{align}
\left| Z^{(3)}(\xi )-A(0)\right| & =\left| A^{\prime }(\xi )Z(\xi )+A(\xi
)\left( I+\int_{0}^{\xi }A(r)Z(r)dr\right) -A(0)\right|  \label{e.8.11} \\
& \leq K_{1}\frac{\sinh (\sqrt{K}\xi )}{\sqrt{K}}+K(\cosh (\sqrt{K}\xi
)-1)+K_{1}\xi \\
& \leq K_{1}\xi (\cosh (\sqrt{K}\xi )+1)+\frac{1}{2}K^{2}\xi ^{2}\cosh (%
\sqrt{K}\xi )  \label{e.8.13} \\
& \leq (2K_{1}s+\frac{1}{2}K^{2}s^{2})\cosh (\sqrt{K}s)  \label{e.8.14}
\end{align}
where we used the elementary inequalities $\sinh (a)/a\leq \cosh (a)$ and $%
\cosh (a)-1\leq \frac{1}{2}a^{2}\cosh (a)$ valid for all $a\in $. Using $%
Z^{(3)}(0)=A(0)$ and the definition of $\mathcal{E}$ completes the proof.
\end{proof}

\subsection{Gaussian Bounds\label{s.8.3}}

In this subsection, $B(s)$ will always denote the standard $\mathbb{R}^{d}$
--valued Brownian motion defined in Notation \ref{n.1.2}.

\begin{lemma}[Fernique]
\label{l.8.4} For $\gamma \in (0,1/2)$ let $K_{\gamma }$ be the random
variable, 
\begin{equation}
K_{\gamma }:=\sup \left\{ \frac{|B(s)-B(r)|}{|s-r|^{\gamma }}:0\leq s<r\leq
1\right\} .  \label{e.8.15}
\end{equation}
Then there exists an $\epsilon =\epsilon (\gamma )>0$ such that $\mathbb{E}%
\left[ e^{\epsilon K_{\gamma }^{2}}\right] <\infty $
\end{lemma}

\begin{proof}
Since $K_{\gamma }$ as a functional of $B$ is a ``measurable'' semi-norm,
equation (\ref{e.8.15}) is a direct consequence of Fernique's theorem 
\cite[Theorem 3.2]{Kuo75}.
\end{proof}

\begin{lemma}
\label{l.8.5}For $p\in \lbrack 1,\infty )$, 
\begin{equation}
\mathbb{E}e^{\frac{p}{2}C\sum_{j=1}^{n}|\Delta
_{j}B|^{2}}=\prod_{j=1}^{n}(1-pC\Delta _{j}s)^{-d/2}  \label{e.8.16}
\end{equation}
provided that $pC\Delta _{j}s<1$ for all $j$. Furthermore, 
\begin{equation}
\lim_{|\mathcal{P}|\rightarrow 0}\mathbb{E}e^{\frac{p}{2}C\sum_{j=1}^{n}|%
\Delta _{j}B|^{2}}=e^{dpC/2}.  \label{e.8.17}
\end{equation}
\end{lemma}

\begin{proof}
By the independence of increments and scaling properties of $B$ we have 
\begin{equation*}
\mathbb{E}\left[ e^{\frac{p}{2}C\sum_{j=1}^{n}|\Delta _{j}B|^{2}}\right]
=\prod_{j=1}^{n}\mathbb{E}\left[ e^{pC|\Delta _{j}B|^{2}/2}\right]
=\prod_{j=1}^{n}\left( \mathbb{E}\left[ e^{pC\Delta _{j}sN^{2}/2}\right]
\right) ^{d},
\end{equation*}
where $N$ is an standard normal random variable. This proves equation (\ref
{e.8.16}), since an elementary Gaussian integration gives 
\begin{equation*}
\mathbb{E}\left[ e^{pC\Delta _{j}sN^{2}/2}\right] =(1-pC\Delta _{j}s)^{-1/2}
\end{equation*}
provided that $pC\Delta _{j}s<1.$ Equation \ref{e.8.17} is an elementary
consequence of (\ref{e.8.16}).
\end{proof}

\begin{lemma}[Gaussian Bound]
\label{l.8.6} For every $k\geq 0$ there is a constant $C=C(k,d)$ which is
increasing in $k$ such that 
\begin{equation}
\mathbb{E}[e^{k|B(1)|}:|B(1)|\geq \rho ]\leq Ce^{-\frac{1}{4}\rho ^{2}}/\rho
^{2}\text{ for all }\rho \geq 1.  \label{e.8.18}
\end{equation}
\end{lemma}

\begin{proof}
A compactness argument shows that there is a constant $\tilde{C}(k,d)$ such
that $r^{d-1}e^{kr}e^{-\frac{1}{2}r^{2}}\leq \tilde{C}(k,d)e^{-\frac{3}{8}%
r^{2}}$ for all $r\geq 0.$ Passing to polar coordinates and using this
inequality shows that 
\begin{align*}
\mathbb{E}[e^{k|B(1)|}& :|B(1)|\geq \rho ]=\omega _{d-1}(2\pi
)^{d/2}\int_{\rho }^{\infty }r^{d-1}e^{kr}e^{-\frac{1}{2}r^{2}}dr \\
& \leq \omega _{d-1}(2\pi )^{d/2}\tilde{C}(k,d)\int_{\rho }^{\infty }\frac{r%
}{\rho }e^{-\frac{3}{8}r^{2}}dr \\
& =\omega _{d-1}(2\pi )^{d/2}\tilde{C}(k,d)\frac{4}{3\rho }e^{-\frac{3}{8}%
\rho ^{2}} \\
& \leq Ce^{-\frac{1}{4}\rho ^{2}}/\rho ^{2},
\end{align*}
where $\omega _{d-1}$ is the volume of the $d-1$ sphere in $\mathbb{R}^{d}.$
\end{proof}

\begin{lemma}
\label{l.8.7} Fix $\epsilon >0$ and $K\geq 0$. Let $\chi _{\epsilon
}(r)=1_{r\geq \epsilon }$, let $B$ be a standard $\mathbb{R}^{d}$--valued
Brownian motion and let $\mathcal{P}=\{0=s_{0}<s_{1}<\cdots <s_{n}=1\}$ be a
partition of $[0,1]$.

Define the function $\psi :\mathbb{R}_{+}\rightarrow \mathbb{R}_{+}$ by 
\begin{equation}
\psi (u):=\mathbb{E}\left[ \left( \frac{\sinh (\sqrt{K}|B(u^{2})|)}{\sqrt{K}%
|B(u^{2})|}\right) ^{d-1}\right] =\mathbb{E}\left[ \left( \frac{\sinh (\sqrt{%
K}u|B(1)|)}{\sqrt{K}u|B(1)|}\right) ^{d-1}\right] .  \label{e.8.19}
\end{equation}
Then there is a constants $C=C(K,d)<\infty $ such that 
\begin{equation}
\sum_{i=1}^{n}\mathbb{E}\left[ \chi _{\epsilon }(|\Delta _{i}B|)\left( \frac{%
\sinh (\sqrt{K}|\Delta _{i}B|)}{\sqrt{K}|\Delta _{i}B|}\right) ^{d-1}\right]
\prod_{j\neq i}\psi (\sqrt{\Delta _{j}s})\leq C\epsilon ^{-2}\exp \left( -%
\frac{\epsilon ^{2}}{4|\mathcal{P}|}\right) .  \label{e.8.20}
\end{equation}
\end{lemma}

\begin{proof}
It is easily checked that $\psi $ is an even smooth (in fact analytic)
function and that $\psi (u)=1+\frac{d(d-1)}{6}u^{2}+O(u^{4})$ and hence
there is a constant $C<\infty $ such that $\psi (u)\leq e^{Cu^{2}}$for $%
0\leq u\leq 1.$ Thus 
\begin{equation*}
\prod_{j\neq i}\psi (\sqrt{\Delta _{j}s})\leq e^{C\sum_{j\neq i}\Delta
_{j}s}\leq e^{C}.
\end{equation*}

Recall the elementary inequality $\sinh (a)/a\leq \cosh (a)\leq e^{|a|}$
which is valid for all $a\in \mathbb{R}.$ Using this inequality and the
scaling properties of $B$ and Lemma \ref{l.8.6}, 
\begin{eqnarray*}
\mathbb{E}\left[ \chi _{\epsilon }(|\Delta _{i}B|)\left( \frac{\sinh (\sqrt{K%
}|\Delta _{i}B|)}{\sqrt{|K|}|\Delta _{i}B|}\right) ^{d-1}\right] &=&\mathbb{E%
}\left[ \chi _{\epsilon }(\sqrt{\Delta _{i}s}|B(1)|)\left( \frac{\sinh (%
\sqrt{K\Delta _{i}s}|B(1)|)}{\sqrt{K\Delta _{i}s}|B(1)|}\right) ^{d-1}\right]
\\
&\leq &\mathbb{E}\left[ \chi _{\epsilon \Delta _{i}s^{-1/2}}(|B(1)|)\exp
\left( (d-1)\sqrt{K\Delta _{i}s}|B(1)|\right) \right] \\
&\leq &C(d,K|\mathcal{P}|)\frac{\Delta _{i}s}{\epsilon ^{2}}\exp \left( -%
\frac{\epsilon ^{2}}{4\Delta _{i}s}\right) \\
&\leq &C(d,K|\mathcal{P}|)\frac{\Delta _{i}s}{\epsilon ^{2}}\exp \left( -%
\frac{\epsilon ^{2}}{4|\mathcal{P}|}\right) .
\end{eqnarray*}
Combining the above estimates completes the proof of Lemma \ref{l.8.7}.
\end{proof}

\begin{prop}
\label{p.8.8} Let $B$ be the $\mathbb{R}^{d}$--valued Brownian motion
defined on $\left( \mathrm{W}(\mathbb{R}^{d}),\mu \right) $ as in Notation 
\ref{n.1.2} above and let Let $\mathbf{R}_{i}$ for $i=0,1,,\ldots ,n$ be a
random symmetric $d\times d$ matrices which are $\sigma (B_{s}:s\leq
s_{i-1}) $--measurable for each $i$. Note that $\mathbf{R}_{0}$ is
non-random. Further assume there is a non-random constant $K<\infty $ such
that $|\mathbf{R}_{i}|\leq K$ for all $i$. Then for all $p\in \mathbb{R}$
there is an $\epsilon =\epsilon (K,d,p)>0$ 
\begin{equation}
1\leq \mathbb{E}\left[ e^{p\sum_{i=1}^{n}(\langle \mathbf{R}_{i}\Delta
_{i}B,\Delta _{i}B\rangle -{\text{tr}}(\mathbf{R}_{i})\Delta _{i}s)}\right]
\leq e^{dp^{2}K^{2}|\mathcal{P}|}  \label{e.8.21}
\end{equation}
whenever $|\mathcal{P}|\leq \epsilon .$
\end{prop}

\begin{proof}
By It\^{o}'s Lemma, 
\begin{equation*}
\langle \mathbf{R}_{i}\Delta _{i}B,\Delta _{i}B\rangle -{\text{tr}}(\mathbf{R%
}_{i})\Delta _{i}s=2\int_{s_{i-1}}^{s_{i}}\langle \mathbf{R}_{i}(B(s)-B(%
\underline{s})),dB(s)\rangle ,
\end{equation*}
and hence $\sum_{i=1}^{n}(\langle \mathbf{R}_{i}\Delta _{i}B,\Delta
_{i}B\rangle -{\text{tr}}(\mathbf{R}_{i})\Delta _{i}s)=M_{1}$, where $M_{t}$
is the continuous square integrable martingale 
\begin{equation*}
M_{t}:=2\int_{0}^{t}\langle \mathbf{R}_{s}(B(s)-B(\underline{s}%
)),dB(s)\rangle
\end{equation*}
and $\mathbf{R}_{s}:=\mathbf{R}_{i}$ if $s\in (s_{i-1},s_{i}].$ The
quadratic variation of this martingale is 
\begin{equation*}
\langle M\rangle _{t}=4\int_{0}^{t}|\mathbf{R}_{s}(B(s)-B(\underline{s}%
))|^{2}ds\leq 4K^{2}\int_{0}^{t}|B(s)-B(\underline{s})|^{2}ds.
\end{equation*}
Let $p\in (1,\infty ).$ Then by the independent increment property of the
Brownian motion $B,$ it follows that 
\begin{eqnarray}
\mathbb{E}\left[ e^{p^{2}\langle M\rangle _{1}}\right] &\leq &\mathbb{E}%
\left[ \exp \left( 4p^{2}K^{2}\int_{0}^{1}|B(s)-B(\underline{s}%
)|^{2}ds\right) \right]  \notag \\
&=&\prod_{i=1}^{n}\mathbb{E}\left[ \exp \left(
4p^{2}K^{2}\int_{s_{i-1}}^{s_{i}}|B(s)-B(\underline{s})|^{2}ds\right) \right]
\notag \\
&=&\prod_{i=1}^{n}\mathbb{E}\left[ \exp \left( 4p^{2}K^{2}\Delta
_{i}s^{2}\int_{0}^{1}|B(s)|^{2}ds\right) \right] ,  \label{e.8.22}
\end{eqnarray}
wherein the last equality we haves used scaling and independence properties
of $B$ to conclude that $\int_{s_{1-1}}^{s_{i}}|B(s)-B(\underline{(}%
s)|^{2}ds,$ $\int_{0}^{\Delta _{i}s}|B(s)|^{2}ds$ and $\int_{0}^{\Delta
_{i}s}\Delta _{i}s|B(\frac{s}{\Delta _{i}s})|^{2}ds=\Delta
_{i}s^{2}\int_{0}^{1}|B(s)|^{2}ds$ all have the same distribution.

Fernique's theorem \cite[Theorem 3.2]{Kuo75} implies that 
\begin{equation*}
\psi (\lambda ):=\mathbb{E}\left[ \exp \left( \frac{\lambda }{2}%
\int_{0}^{1}|B(s)|^{2}ds\right) \right]
\end{equation*}
is a well defined analytic function of $\lambda $ in a neighborhood of $0.$
Because $\psi (0)=1$ and 
\begin{equation*}
\psi ^{\prime }(0)=\frac{1}{2}\mathbb{E}\int_{0}^{1}|B(s)|^{2}ds=\frac{d}{4}
\end{equation*}
it follows that $\psi (\lambda )\leq e^{d\lambda /2}$ for all positive $%
\lambda $ sufficiently near $0.$ Using this fact in equation \ref{e.8.22}
gives the bound 
\begin{eqnarray}
\mathbb{E}\left[ e^{p^{2}\langle M\rangle _{1}}\right] &\leq
&\prod_{i=1}^{n}\exp \left( 2dp^{2}K^{2}\Delta _{i}s^{2}\right) =\exp \left(
2dK^{2}p^{2}\sum_{i=1}^{n}\Delta _{i}s^{2}\right)  \notag \\
&\leq &\exp \left( 2dK^{2}p^{2}|\mathcal{P}|\right) <\infty ,  \label{e.8.23}
\end{eqnarray}
which is valid when the mesh of $\mathcal{P}$ is sufficiently small.

By It\^{o}'s Lemma, 
\begin{equation*}
Z_{t}^{(p)}=\exp \left( pM_{t}-\frac{p^{2}}{2}\langle M\rangle _{t}\right)
\end{equation*}
is a positive local martingale. Because of the bound in equation \ref{e.8.23}%
, Novikov's criterion \cite[Prop. 1.15, p.308]{Revuz:Yor:ed2} implies that $%
Z_{t}^{(p)}$ is in fact a martingale and hence in particular $\mathbb{E}%
\left[ Z_{1}^{(p)}\right] =1$. Therefore, 
\begin{equation*}
\mathbb{E}\left[ e^{pM_{1}}\right] =\mathbb{E}\left[ e^{pM_{1}-\frac{p^{2}}{2%
}\langle M\rangle _{1}}e^{\frac{p^{2}}{2}\langle M\rangle _{1}}\right] \geq 
\mathbb{E}\left[ e^{pM_{1}-\frac{p^{2}}{2}\langle M\rangle _{1}}\right] =1
\end{equation*}
and 
\begin{eqnarray*}
\mathbb{E}\left[ e^{pM_{1}}\right] &=&\mathbb{E}\left[ \exp \left(
pM_{1}-p^{2}\langle M\rangle _{1}\right) \exp \left( p^{2}\langle M\rangle
_{1}\right) \right] \\
&\leq &\sqrt{\mathbb{E}\left[ \exp \left( 2pM_{1}-2p^{2}\langle M\rangle
_{1}\right) \right] }\sqrt{\mathbb{E}\left[ \exp \left( p^{2}\langle
M\rangle _{1}\right) \right] } \\
&=&\sqrt{\mathbb{E}Z_{t}^{(2p)}}\sqrt{\mathbb{E}\left[ \exp \left(
p^{2}\langle M\rangle _{1}\right) \right] }=\sqrt{\mathbb{E}\left[ \exp
\left( p^{2}\langle M\rangle _{1}\right) \right] } \\
&\leq &\exp \left( dK^{2}p^{2}|\mathcal{P}|\right) .
\end{eqnarray*}
This completes the proof of Proposition \ref{p.8.8}.
\end{proof}

%\bibliographystyle{amsplain}
%\bibliography{density2}

\providecommand{\bysame}{\leavevmode\hbox to3em{\hrulefill}\thinspace}

\end{document}